\documentclass[12pt,leqno]{article}

\usepackage[francais]{babel}
\usepackage{amsmath,amssymb,amscd,a4wide}
\usepackage[latin1]{inputenc}

\newenvironment{paragr}[1][]{\refstepcounter{subsection} \noindent \textbf{\thesubsection . \ #1}}{\medskip}

\newenvironment{theoreme}{ \medskip\refstepcounter{theo}  \noindent\textbf{Th\'eor\`eme \thetheo}. ---\em}{\em \medskip}
\newenvironment{proposition}{\medskip\refstepcounter{theo}   \noindent\textbf{Proposition \thetheo}. ---\em}{\em\medskip}
\newenvironment{corollaire}{\medskip\refstepcounter{theo}  \noindent\textbf{Corollaire \thetheo}. ---\em}{\em\medskip}
\newenvironment{remarque}{\medskip \noindent \textbf{Remarque}. --- }{}

\newenvironment{lemme}{\medskip\refstepcounter{theo}   \noindent\textbf{Lemme \thetheo}. ---\em}{\em\medskip}
\newenvironment{conjecture}{ \medskip\refstepcounter{theo}  \noindent\textbf{Conjecture \thetheo}. ---\em}{\em \medskip}

\newenvironment{definition}{\medskip\refstepcounter{theo}  \noindent\textbf{D\'efinition \thetheo}. ---}{\medskip}

\newenvironment{preuve}[1][]{\noindent \textbf{Démonstration.} #1 --- }{\hfill
  \ensuremath{\square} \medskip}

\DeclareMathOperator{\nr}{nr}

\DeclareMathOperator{\Ad}{Ad}

\DeclareMathOperator{\mes}{mes}
\DeclareMathOperator{\der}{der}
\DeclareMathOperator{\scnx}{sc}

\DeclareMathOperator{\Norm}{Norm}

\DeclareMathOperator{\Gal}{Gal}
\DeclareMathOperator{\Hom}{Hom}

\DeclareMathOperator{\Id}{Id}

\DeclareMathOperator{\Ker}{Ker}
\DeclareMathOperator{\Cok}{Cok}

\DeclareMathOperator{\val}{val}
\DeclareMathOperator{\SL}{SL}
\DeclareMathOperator{\Spec}{Spec}
\DeclareMathOperator{\Div}{Div}
\DeclareMathOperator{\Tor}{Tor}

\DeclareMathOperator{\trace}{trace}
\DeclareMathOperator{\Sym}{Sym}
\newcommand{\ZZ}{\mathbb{Z}}
\newcommand{\Gm}{\mathbb{G}_m}

\newcommand{\NN}{\mathbb{N}}
\newcommand{\RR}{\mathbb{R}}

\newcommand{\CC}{\mathbb{C}}
\newcommand{\QQ}{\mathbb{Q}}
\newcommand{\FF}{\mathbb{F}}
\newcommand{\Fq}{\mathbb{F}_q}

\newcommand{\Ql}{\mathbb{Q}_{\ell}}
\newcommand{\Qlb}{\bar{\mathbb{Q}}_{\ell}}

\newcommand{\oc}{\mathcal{O}}

\newcommand{\Sc}{\mathcal{S}}

\newcommand{\ec}{\mathcal{E}}

\newcommand{\lc}{\mathcal{L}}
\newcommand{\dc}{\mathcal{D}}
\newcommand{\fc}{\mathcal{F}}
\newcommand{\pc}{\mathcal{P}}

\newcommand{\Gc}{\widehat{G}}
\newcommand{\Tc}{\widehat{T}}
\newcommand{\Pc}{\widehat{P}}

\newcommand{\Mc}{\widehat{M}}

\newcommand{\Uc}{\widehat{U}}
\newcommand{\Hc}{\widehat{H}}

\newcommand{\PP}{\mathbf{P}}

\newcommand{\ggo}{\mathfrak{g}}

\newcommand{\of}{\mathfrak{o}}

\newcommand{\mgo}{\mathfrak{m}}
\newcommand{\ngo}{\mathfrak{n}}
\newcommand{\ago}{\mathfrak{a}}

\newcommand{\hgo}{\mathfrak{h}}
\newcommand{\tgo}{\mathfrak{t}}

\newcommand{\Xgo}{\mathfrak{X}}
\newcommand{\Fgo}{\mathfrak{F}}

\newcommand{\al}{\alpha}

\newcommand{\la}{\lambda}

\newcommand{\back}{\backslash}

\newcommand{\bg}{\langle}
\newcommand{\bd}{\rangle}

\newcommand{\eps}{\varepsilon}

\renewcommand{\leq}{\leqslant}
\renewcommand{\geq}{\geqslant}

\title{Sur l'homologie des fibres de Springer affines tronquées}
\author{Pierre-Henri Chaudouard et Gérard Laumon}
\date{ }
\begin{document}
\maketitle

\makeatletter
\@addtoreset{equation}{section}         
\def\theequation{\thesection.\arabic{equation}}
\makeatother

\begin{abstract} Following Goresky, Kottwitz and MacPherson, we compute the homology of truncated
affine Springer fibers in the unramified case but under a purity assumption. We prove this assumption in the equivalued case. The truncation parameter
is viewed as a divisor on an $\ell$-adic toric variety. For each fiber, we introduce a graded quasi-coherent sheaf on the toric variety
whose space of global sections is precisely the $\ell$-adic homology of the truncated affine Springer fiber. Moreover, for some
families of endoscopic groups, these sheaves show up in an exact sequence. As a consequence we prove Arthur's weighted fundamental lemma
in the unramified equivalued case.
\end{abstract}

\tableofcontents

\section{Introduction}

\begin{paragr}
Le lemme fondamental de Langlands et Shelstad est une famille
d'identités combinatoires entre int\'{e}grales orbitales sur les
groupes r\'{e}ductifs sur un corps local non archim\'{e}dien $F$.  Ici
nous ne consid\'{e}rons que le cas o\`{u} $F$ est de caract\'{e}ristiques égales et nous nous res\-treignons \`{a} la variante du
lemme fondamental pour les alg\`{e}bres de Lie.  Des r\'{e}sultats de
Waldspurger (cf. \cite{W1}) permettent de ramener le lemme fondamental pour les
algèbres de Lie sur un corps local non archim\'{e}dien $F$ arbitraire (mais de caractéristique résiduelle assez grande) \`{a}
cette variante pour $F$ de caract\'{e}ristiques égales. D'autre part, en caractéristique nulle, le lemme fondamental sur les groupes est essentiellement équivalent au lemme fondamental sur les algèbres de Lie (pour des résultats plus généraux, cf. \cite{W2}). 
\bigskip

Les int\'{e}grales orbitales qui interviennent dans le lemme 
fondamental admettent alors une interpr\'{e}tation
cohomologique gr\^{a}ce au dictionnaire fonctions-faisceaux de
Grothendieck.  Dans \cite{GKM} et \cite{GKM-purete}, Goresky, Kottwitz et MacPherson ont
utilis\'{e} cette interpr\'{e}tation cohomologique et le
th\'{e}or\`{e}me de localisation d'Atiyah-Borel-Segal en cohomologie
\'{e}quivariante pour d\'{e}montrer le lemme fondamental de
Langlands-Shelstad dans de nombreux cas.
\bigskip

Le lemme fondamental pond\'{e}r\'{e} d'Arthur est une extension du
lemme fondamental de Langlands-Shelstad aux int\'{e}grales orbitales
pond\'{e}r\'{e}es (cf. conjecture  5.1 de \cite{STF1}).  Rappelons que le c\^{o}t\'{e} g\'{e}om\'{e}trique
de la formule des traces d'Arthur-Selberg est une combinaison
lin\'{e}aire d'int\'{e}grales orbitales pond\'{e}r\'{e}es globales qui se
d\'{e}vissent en int\'{e}grales orbitales pond\'{e}r\'{e}es locales (cf. \cite{inv-trace}), et que pour stabiliser cette formule des traces, et donc obtenir
le transfert endoscopique du programme de Langlands, il faut disposer
du lemme fondamental pond\'{e}r\'{e}.
\bigskip

Dans ce travail nous \'{e}tendons l'approche de Goresky, Kottwitz et
MacPherson au cas du lemme fondamental pond\'{e}r\'{e} d'Arthur.  Les
limitations de notre travail sont les m\^{e}mes que celles de \cite{GKM} et \cite{GKM-purete} : nous ne pouvons atteindre que les int\'{e}grales orbitales pondérées en
les \'{e}l\'{e}ments semi-simples de valuations égales qui sont
contenus dans les tores non ramifi\'{e}s.
\end{paragr}

\begin{paragr} \label{par:intro2}Plus pr\'{e}cis\'{e}ment, soit $k$ une clôture algébrique du corps fini $\Fq$ et $\tau$ le Frobenius. soit $G$ un groupe r\'{e}ductif connexe, $T$ un sous-tore maximal de
$G$ et $M$ un sous-groupe de Levi de $G$ contenant $T$, tous ces groupes étant définis sur $\mathbb{F}_{q}$ . Dans cette introduction, pour simplifier, on suppose de plus que $G$ est semi-simple, simplement connexe et d\'{e}ploy\'{e} sur $\mathbb{F}_{q}$. \\

Soit $\Gc$ le groupe complexe dual de $G$ au sens de Langlands. Soit $s\in \Tc^\tau$ un élément, fixé par le Frobenius $\tau$, du tore complexe dual de $T$. La composante neutre du  centralisateur de $s$ dans $\Mc$ est le dual d'un groupe réductif connexe défini sur $\Fq$ qui admet $T$ comme tore maximal. Appelons ce groupe $M'$ : c'est un groupe endoscopique de $M$.

Lorsque $s_{1}$ parcourt
$sZ_{\hat{M}}$ o\`{u} $Z_{\hat{M}}$ est le centre de $\hat{M}$, la famille formée des composantes neutres $\hat{G}_{s_{1}}$ des centralisateurs de $s_{1}$ dans
$\hat{G}$ est finie. Indexons la sous-famille des centralisateurs maximaux (pour l'inclusion) par l'ensemble $[n]=\{1,\ldots ,n\}$ pour un certain entier $n$. Pour toute partie $I$ non vide de $[n]$, on pose 
$$\Gc_I= (\bigcap_{i\in I}\hat{G}_{i})^{0}.$$
Par dualité de Langlands, on obtient une famille de groupes réductifs connexes $G_I$ sur $\Fq$ qui admettent $M'$ comme sous-groupe de Lévi. On complète cette famille en posant $G_\emptyset=G$. Lorsque la différence $J-I$ est un singleton le groupe $G_J$ est un groupe endoscopique de $G_I$. 

\bigskip

Soit $\mathfrak{t}\subset \mathfrak{m}$ les alg\`{e}bres de Lie de
$T\subset M$. Soit $\of=\Fq[[\eps]]\subset F=\Fq((\eps))$  et $1_{\mathfrak{k}}$ la fonction caract\'{e}ristique du
$\of$-r\'{e}seau $\mathfrak{k}=\mathfrak{g}(\of)$ dans $\mathfrak{g}(F)$.
Soit $\gamma\in
\mathfrak{t}(\of)\subset\mathfrak{g}(F)$ un
\'{e}l\'{e}ment $G$-r\'{e}gulier.

Pour définir des intégrales orbitales pondérés, nous avons besoin d'introduire un certain paramètre de troncature adapté à $M$ que nous notons $D$. On lui associe une fonction poids à valeurs entières 
$$x \in G(F) \mapsto v^G_D(x)$$
qui est invariante à gauche par $M(F)$ et à droite par $G(\of)$. 

 Pour tout $\gamma'\in \mgo(F)$ stablement conjugué à $\gamma$, on définit à la suite de Langlands un accouplement $\bg s, \gamma'\bd$ qui est une racine de l'unité. Celui-ci ne dépend que de la classe de $M(F)$-conjugaison de $\gamma'$. Pour un tel $\gamma'$, on note $T_{\gamma'}$ le centralisateur de $\gamma'$ dans $G$ : c'est un tore maximal.
 
On peut alors former la $s$-int\'{e}grale orbitale pond\'{e}r\'{e}e
$$
\mathrm{J}_{D}^{G,s}(\gamma )=\sum_{\gamma'}\langle s,\gamma'\rangle
\int_{T_{\gamma'}(F)\backslash G(F)}1_{\mathfrak{k}}(\mathrm{ad}(g^{-1})(\gamma'))
v_{D}^G(g)\frac{\mathrm{d}g}{\mathrm{d}t}
$$
o\`{u} $\gamma'$ parcourt un système de représentants des classes de $M(F)$-conjugaison dans la  classe stable de $\gamma$. 

Pour chaque partie $I\subset [n]$ non vide, le paramètre $D$ est adapté à $M'$. On a donc de même un poids $v_D^{G_I}$ sur $G_I(F)$ invariant à gauche par $M'(F)$ et à droite par $G_I(\of)$ et une $s$-int\'{e}grale orbitale pond\'{e}r\'{e}e $\mathrm{J}_{D}^{G_I,s}(\gamma )$.

Notre r\'{e}sultat principal est alors le suivant:

\begin{theoreme}
Si la caract\'{e}ristique de $\mathbb{F}_{q}$ est assez grande par
rapport \`{a} $G$,  pour des éléments  $\gamma$ de valuations égales et pour des paramètres de troncature $D$ assez ``réguliers'',  on a la 
relation
$$
\Delta_{s}(\gamma )\mathrm{J}_{D}^{G,s}(\gamma )=
\sum_{I}(-1)^{|I|-1}q^{d_{I}(\gamma)}\mathrm{J}_{D}^{G_{I},s}(\gamma ).
$$
o\`{u} $\Delta_{s}(\gamma )$ est un facteur de transfert de 
Langlands-Shelstad et $d_{I}(\gamma)$ est la demi-somme des 
valuations de $\alpha (\gamma)$ pour $\gamma$ parcourant les racines 
de $G$ qui ne sont pas dans $G_{I}$.
\end{theoreme}

Goresky, Kottwitz et MacPherson ont montr\'{e} que les
$s$-int\'{e}grales orbitales s'interpr\`{e}tent
g\'{e}om\'{e}triquement comme des sommes finies index\'{e}es par les
points rationnels de certains quotients des fibres de Springer
affines.  De m\^{e}me les $s$-int\'{e}grales orbitales
pond\'{e}r\'{e}es s'interpr\`{e}tent comme des sommes finies
index\'{e}es par les points rationnels de certains quotients des
fibres de Springer affines tronqu\'{e}es.  Elles admettent donc une
interpr\'{e}tation en cohomologique $\ell$-adique  gr\^{a}ce \`{a} la formule des points fixes de Grothendieck-Lefschetz.

Une grande partie de ce travail porte sur l'\'{e}tude de la 
cohomologie $\ell$-adique des fibres de Springer affines 
tronqu\'{e}es.
\end{paragr}

\begin{paragr} Rappelons que $k$ est une clôture algébrique de $\Fq$. On note $\Xgo^G$ 
 la grassmannienne affine de $G$. C'est un ind-$k$-sch\'{e}ma
projectif dont l'ensemble des $k$-points est $G(k((\eps)))/G(k[[\varepsilon]])$.  La fibre de Springer affine en
$\gamma\in\mathfrak{g}(k((\varepsilon )))$ est le ferm\'{e}
$\mathfrak{X}_{\gamma}^{G}$ de $\mathfrak{X}^{G}$ d\'{e}fini par la
condition $\mathrm{Ad}(g)^{-1}\gamma\in
\mathfrak{g}(k[[\varepsilon]])$.  Si $\gamma$ est
r\'{e}gulier semi-simple, la fibre de Springer affine
$\mathfrak{X}_{\gamma}^{G}$ est un $k$-sch\'{e}ma localement de type
fini et de dimension finie (cf. \cite{KL}).  Dans la suite on se limite aux
$\gamma\in\mathfrak{t}(k[[\varepsilon ]])\subset
\mathfrak{g}(k((\varepsilon)))$ r\'{e}guliers.
\bigskip

Les tronqu\'{e}s de $\mathfrak{X}_{\gamma}^{G}$ que nous
consid\'{e}rons ici sont les intersections de
$\mathfrak{X}_{\gamma}^{G}$ avec des tronqu\'{e}s de la grassmannienne
affine.  Soit $a_M=\Hom_{\ZZ}(X^\ast(M),\ZZ)$. Les tronqu\'{e}s de $\mathfrak{X}^{G}$ ont
\'{e}t\'{e} d\'{e}finis par Arthur; ils sont relatifs au choix d'un
sous-groupe de Levi $M$ de $G$ contenant $T$ et d'un param\`{e}tre de
troncature qui est une famille $\la =(\la_{P})_{P}$ de points $\la_P\in a_M$, indexée par les sous-groupes
paraboliques $P$ de $G$ de Levi $M$.  Pour tout tel $P$, en utilisant la décomposition d'Iwasawa $G(k((\eps)))=P(k((\eps))) G(k[[\eps]])$, on définit une fonction 
$$H_{P}\ :\ G(k((\eps)))/G(k[[\varepsilon]]) \longrightarrow a_M.$$
Le tronqué par $\la$ de la grassmannienne affine $\Xgo^G$ est l'ensemble des $x\in \Xgo^G$ qui vérifie : pour tout sous-groupe parabolique $P$ de $G$ de Lévi $M$, le point $H_P(x)$ se trouve dans l'enveloppe convexe de la famille $\la$ dans $a_M\otimes \RR$.

Lorsque de plus la famille $\la$  est orthogonale au sens d'Arthur, elle s'interprète naturellement comme un diviseur $D$ sur une compactification torique partielle $Y$ du tore dual $\hat{T}$ de $T$. Pour notre calcul de cohomologie $\ell$-adique, il est naturel de
consid\'{e}rer $\hat{T}$ comme un tore sur $\Qlb$, une clôture algébrique de $\Ql$ avec $\ell$ premier à $q$.  La compactification torique partielle $Y$ de $\hat{T}$ est associ\'{e}e
\`{a} l'\'{e}ventail $\Sigma_{M}^{G}$ form\'{e} des c\^{o}nes
$a_{P}^{G,+}$ pour $P$ parcourant les sous-groupes paraboliques de $G$
contenant $M$.  Par diviseur on entend un diviseur de Weil $\hat{T}$-\'{e}quivariant : c'est donc une combinaison lin\'{e}aire \`{a} coefficients entiers des composantes irr\'{e}ductibles du bord de $Y-\hat{T}$.

Pour chaque diviseur $D$, on a donc un tronqu\'{e}
$\mathfrak{X}^{G}(D)$ de la grassmannienne affine $\mathfrak{X}^{G}$
qui est un sous-ind-sch\'{e}ma ferm\'{e}, et par intersection avec la
fibre de Springer affine on obtient un tronqu\'{e}
$\mathfrak{X}_{\gamma}^{G}(D)$ qui est un ferm\'{e} de
$\mathfrak{X}_{\gamma}^{G}$.  Le groupe $X_{\ast}(T\cap
M_{\mathrm{der}})$ agit librement sur $\mathfrak{X}_{\gamma}^{G}(D)$
et le quotient $X_{\ast}(T\cap
M_{\mathrm{der}}) \back \mathfrak{X}_{\gamma}^{G}(D)$ est $k$-sch\'{e}ma projectif.
\bigskip

Une propri\'{e}t\'{e} conjecturale essentielle des fibres de Springer 
affine  est la 
propri\'{e}t\'{e} de puret\'{e} de leur homologie.

\begin{conjecture} (Goresky, Kottwitz et MacPherson)
Pour tout entier $n\geq 0$, le groupe d'homologie $\ell$-adique
$H_{n}(\mathfrak{X}_{\gamma}^{G},\Qlb)$ est pur de 
poids $n$ au sens de Grothendieck et Deligne.
\end{conjecture}

De m\^{e}me pour les fibres de Springer affines tronqu\'{e}es on 
conjecture:

\begin{conjecture} Soit $\gamma\in\mathfrak{t}(k[[\varepsilon ]])$ un élément régulier. Pour tout diviseur $D$ ``assez régulier'' sur la compactification partielle $Y$ de $\hat{T}$, la fibre de Springer affine tronquée $\Xgo_\gamma^G(D)$ est pure au sens de la définition \ref{def:purete}. En particulier, pour tout entier $n\geq 0$, le groupe d'homologie
$\ell$-adique $H_{n}(\mathfrak{X}_{\gamma}^{G}(D),\Qlb)$
est lui aussi pur de poids $n$ au sens de Grothendieck et Deligne. 
\end{conjecture}

On a alors les résultats suivants.

\begin{theoreme} Soit $\gamma\in\mathfrak{t}(k[[\varepsilon ]])$ un élément r\'{e}gulier  et équivalué (au sens de la définition \ref{def:equivalue}) alors la fibre de Springer affine tronquée $\Xgo_\gamma^G(D)$ est pure au sens de la définition \ref{def:purete} pour tout diviseur $D$ sur $Y$ assez ``régulier''.
\end{theoreme}

\begin{theoreme} Il existe un $\mathcal{O}_{Y}$-Module quasi-coh\'{e}rent gradu\'{e} $\mathfrak{F}^{G}=\bigoplus_{n\geq
0}\mathfrak{F}_{n}^{G}$ qui vérifie la propriété suivante : pour tout entier $n$, le $\mathcal{O}_{Y}$-Module $\mathfrak{F}_{n}^{G}$ est cohérent et  pour tout diviseur $D$ assez ``régulier'' tel que la fibre de Springer affine tronquée $\Xgo_\gamma^G(D)$ soit pure au sens de la définition \ref{def:purete}, on a
$$
H_{2n}(\mathfrak{X}_{\gamma}^{G}(D),\Qlb)=
H^{0}(Y,\mathfrak{F}_{n}^{G}(D)),
$$
et les $H^{i}(Y,\mathfrak{F}_{n}^{G}(D))$ pour $i>0$ sont tous nuls. 

Par ailleurs, dès que la fibre de Springer affine tronquée $\Xgo_\gamma^G(D)$ est pure, on a 
$$
H_{2n+1}(\mathfrak{X}_{\gamma}^{G}(D),\Qlb)=0.$$
\end{theoreme}

En guise d'illustration, considérons le cas o\`{u} $G=\mathrm{SL}(2)$, $M=T$ est le tore diagonal et  $\gamma =\mathrm{diag}(\varepsilon^d,-\varepsilon^d)$ où $d\geq 1$ est un entier. 

Soit $\PP^d$ la fibre de Springer classique formée des sous-espaces de dimension $d$ dans un espace de dimension $2d$ stables par le carré d'un nilpotent régulier. Dans ce cas, $\mathfrak{X}_{\gamma}^{G}$ est une cha\^{\i}ne infinie indexée par $\ZZ$ de $\PP^d$ : deux $\PP^d$ successifs indexés par $n$ et $n+1$ se recollent suivant un $\PP^{d-1}$. Cette chaîne  est pavée par des espaces affines standard et son homologie est donc pure. Le tore dual $\Tc$ est le groupe multiplicatif $\mathbb{G}_m$ sur $\Qlb$ et la variété torique est la droite projective $Y=\mathbb{P}_{\Qlb}^{1}$. Soit $(x,y)$ des coordonnées homogènes sur $Y$. Le $\mathcal{O}_{Y}$-Module quasi-coh\'{e}rent $\mathfrak{F}^{G}$ est défini par la donnée du $\Qlb[x,y]$-module gradué
$$\frac{\Qlb[x,y][z]}{\sum_{i=1}^d (x-y)^i \Ker(\partial_z^i)}\cap \Ker \partial_z.$$
La graduation sur $\mathfrak{F}^{G}$ est donnée par la graduation sur l'indéterminée $z$ qui s'explicite de la manière suivante 
$$\sum_{i=0 }^{d-1} (x-y)^i \frac{\Qlb[x,y]}{(x-y)\Qlb[x,y]} z^i + (x-y)^d\Qlb[x,y]z^d.$$

Soit un diviseur $D=n_{0}[0]+n_{\infty}[\infty]$ sur $Y$ avec $n_0$ et $n_\infty$ deux entiers. Ici $D$ est régulier si  $n_0+n_\infty\geq 0$. Si $0\leq n_0+n_\infty \leq d$, $\Xgo_\gamma^G(D)$ s'identifie à $\PP^{n_0+n_\infty}$. Si  $n_0+n_\infty \geq d$, $\Xgo_\gamma^G(D)$ est une sous-chaîne dans $\Xgo_\gamma^G$ formée de $n_0+n_\infty-d +1$ copies de $\PP^d$. D'après le théorème précèdent, l'homologie en degré $2i$ de $\Xgo_\gamma^G(D)$ s'identifie à la partie homogène de degré $n_0+n_\infty$ de
 $$(x-y)^i \frac{\Qlb[x,y]}{(x-y)\Qlb[x,y]}$$
si $i>d$ et de
$$(x-y)^d\Qlb[x,y]$$
si $i=d$.
\end{paragr}

\begin{paragr} Reprenons les notations du \S \ref{par:intro2} : l'élément $s\in \hat{T}^\tau$ définit un  groupe endoscopique $M'$ de $M$. Pour toute partie $I\subset[n]$, on a défini un groupe réductif connexe $G_I$ sur $\Fq$. Soit $Y_s$ l'adhérence dans $Y$ de la classe $sZ_{\hat{M}}\subset \hat{T}$. Pour tout $\mathcal{O}_{Y}$-Module coh\'{e}rent $\mathfrak{F}$, soit $\hat{\Fgo}$ le complété de $\Fgo$ le long de $Y_s$. Le lemme fondamental pond\'{e}r\'{e} r\'{e}sulte alors d'une suite  exacte longue
$$
0\longrightarrow \hat{\mathfrak{F}}^{G}\longrightarrow 
\bigoplus_{|I|=1}\hat{\mathfrak{F}}^{G_{I}}\longrightarrow 
\bigoplus_{|I|=2}\hat{\mathfrak{F}}^{G_{I}}\longrightarrow \ldots \longrightarrow 
\hat{\mathfrak{F}}^{G_{[n]}}\longrightarrow 0
$$
de $\mathcal{O}_Y$-Modules coh\'{e}rents sur le compl\'{e}t\'{e} formel 
de $Y$ le long de $Y_{s}$.
\end{paragr}

\begin{paragr} Cet article est organisé ainsi. Les sections \ref{sec:notations} à \ref{sec:suite-spec} sont consacrées à quelques notations et préliminaires cohomologiques. Dans la section \ref{sec:vartorique}, on définit les variétés toriques associées à $\Tc$. On introduit les grassmanniennes affines tronquées et on étudie leurs orbites sous l'action de certains tores aux sections \ref{sec:gat} et \ref{sec:orbitesga}. Dans la section \ref{sec:fsat}, on introduit les fibres de Springer affines tronquées et on prouve leur pureté dans la situation mentionnée plus haut. La démonstration est en deux étapes. On commence par le cas où la troncature est complète et régulière. On montre alors que  la fibre de Springer affine tronquée est une réunion d'orbites sous $G(k[[\eps]])$ intersectées avec la fibre de Springer affine : c'est réminiscent d'un lemme dû à Arthur dans sa preuve de la formule des traces locales. Des résultats antérieurs de Goresky-Kottwitz-MacPherson donnent alors la pureté voulue. On élargit ensuite la troncature dans certaines directions et on montre que la différence est également pure. En suivant la méthode de Goresky-Kottwitz-MacPherson, on procède au calcul de l'homologie équivariante des fibres de Springer affines tronquées d'abord pour $\SL(2)$ à la section  \ref{sec:SL2} puis, à la section \ref{sect:Hom_orb}, pour un groupe quelconque mais sous une hypothèse de pureté. On introduit les faisceaux $\Fgo^G$ sur $Y$ et on prouve dans les cas purs la suite exacte ci-dessus à la section \ref{sec:cplx-fx}. Dans la section \ref{sec:IOP}  finale, on définit les intégrales orbitales pondérées. On donne leur interprétation cohomologique et on déduit des résultats précédents le lemme fondamental pondéré dans le cas équivalué et non ramifié.

\end{paragr}

\begin{paragr}[Remerciements.] --- Nous remercions L. Fargues, L. Illusie, B.C. Ngô et  J.-L. Waldspurger pour l'aide qu'ils nous ont apportée durant la préparation de cet article. Nous remercions également les organisateurs des différents séminaires ou conférences où nous avons pu exposer une partie des résultats présentés ici. Enfin le premier auteur nommé remercie  l'A.C.I. ``Réalisations géométriques de correspondances de Langlands'' pour son aide matérielle.

\end{paragr}

\section{Notations}\label{sec:notations}

\begin{paragr}
Soit $k$ une clôture algébrique d'un corps fini de caractéristique assez grande (cf. la remarque qui suit le théorème \ref{thm:purete}). Soit $F=k((\eps))$ le corps des séries formelles de Laurent et $\of=k[[\eps]]$ l'anneau des séries formelles. Soit $\val$ la valuation de $F$ déterminée par $\val(\eps)=1$.
\end{paragr}

\begin{paragr}
  Pour tout ensemble fini $E$, on note $|E|$ son cardinal.
\end{paragr}

\begin{paragr} Pour tout groupe réductif connexe $G$ défini sur $k$, on note $X^\ast(G)$ le groupe de ses caractères algébriques et $Z_G$ son centre. On note $\ggo$ son algèbre de Lie.
Pour tout tore $T$, on note $X_\ast(T)$ le groupe de ses cocaractères. 
\end{paragr}

\begin{paragr} Soit $T$ un sous-tore maximal de $G$. Pour tout sous-groupe fermé $H$ de $G$ stable par $T$ pour l'action adjointe, on note $\Phi^H(T)$ l'ensemble des racines de $T$ dans $H$.  On note $\Phi^G_+(T)$ un système de représentants du quotient de  $\Phi^G(T)$ par la relation d'équivalence $\al\sim \beta$ si $\al=\pm \beta$. Si $\Phi\subset \Phi^G(T)$ est une partie stable par multiplication par $\pm 1$, on pose $\Phi_+=\Phi\cap\Phi^G_+(T)$.

À toute racine $\al \in \Phi^G(T)$, on associe de la manière usuelle une coracine $\al^\vee \in X_\ast(T)$. On définit alors $\Phi^{H,\vee}(T)$ comme l'ensemble des coracines associées aux racines dans $\Phi^H(T)$.

 Pour tout sous-groupe de Borel $B$ de $G$ qui contient $T$, on note $\Delta_B$ l'ensemble des racines simples dans $\Phi^B(T)$.
\end{paragr}

\begin{paragr}
On note $\fc^G(T)$ l'ensemble des sous-groupes paraboliques de $G$ qui contiennent $T$. Pour tout $P$ dans $\fc^G(T) $, on note $N_P$ le radical unipotent de $P$ et $M_P$ l'unique sous-groupe de Lévi de $P$ qui contient $T$.  On note $\lc^G(T)$ l'ensemble des $M_P$  pour $P\in \fc^G$. Dans la suite, lorsque le tore $T$ est clairement sous-entendu, on appelle simplement ``sous-groupes de Lévi'' de tels sous-groupes.

On omet le $T$ dans les notations $\Phi^G$, $\fc^G$ et $\lc^G$ dès que le contexte est clair.
\end{paragr}

\begin{paragr}  Dans toute la suite, on fixe un nombre premier $\ell$ premier à la caractéristique de $k$, ainsi que $\Qlb$ une clôture algébrique de $\QQ_{\ell}$ le corps des nombres $\ell$-adiques. Considérons la donnée radicielle $(X^\ast(T),\Phi^G,X_\ast(T),\Phi^{G,\vee})$ associée au couple $(G,T)$. On fixe alors un couple $(\Gc,\Tc)$ formé d'un groupe réductif connexe $\Gc$ et d'un sous-tore maximal $\Tc$ tous deux définis sur $\Qlb$ de sorte que la donnée radicielle associée à $(\Gc,\Tc)$ soit duale de celle de $(G,T)$ au sens où
$$(X^\ast(\Tc),\Phi^{\Gc},X_\ast(\Tc),\Phi^{\Gc,\vee})=(X_\ast(T),\Phi^{G,\vee},X^\ast(T),\Phi^G).$$
L'application qui, à une racine, associe une coracine, induit une bijection de $\Phi^G$ sur  $\Phi^{G,\vee}=\Phi^{\Gc}$. On en déduit des bijections $\lc^G \to \lc^{\Gc}$ et $\fc^{G}\to \fc^{\Gc}$, notées $H\mapsto \Hc$. Si $M\in \lc^G$, le groupe réductif $\Mc$ est un groupe dual pour $M$. 
\end{paragr}

\begin{paragr}
Soit $s\in \Tc$. Pour tout sous-groupe fermé $H$ de $\Gc$ qui contient $\Tc$, on note 
$$H_s=Z_H(s)^0$$ 
la composante neutre du centralisateur de $s$ dans $H$. C'est un groupe réductif connexe qui contient $\Tc$.  

 Un groupe endoscopique associé au triplet $(G,T,s)$ est un groupe $G'$ réductif, connexe, défini sur $k$ et muni d'un plongement de $T$ dans $G'$ de sorte que la donnée  radicielle de $(G',T)$ soit duale de celle de $(\Gc_s,\Tc)$. 

Si $M\in \lc^G$, on a $\Mc\in \lc^{\Gc}$ et on peut trouver $s\in \Tc$ de sorte que $\Mc=\Gc_s$. Par conséquent, $M$ est un groupe endoscopique associé au triplet $(G,T,s)$.
\end{paragr}

\begin{paragr}
Pour tout $M\in \lc^G$, on définit le $\QQ$-espace vectoriel $\ago_M^\ast=X^\ast(M)\otimes \QQ$. Le morphisme de restriction $X^\ast(M)\to X^\ast(T)$ induit une injection $\ago_M^\ast\to \ago_T^\ast$. On note $(\ago_T^M)^\ast$ le sous-espace de $\ago_T^\ast$  engendré par $\Phi^M$. On a une décomposition en somme directe 
$$\ago_T^\ast=(\ago_T^M)^\ast\oplus \ago_M^\ast.$$
Plus généralement, si $L\in \lc^G$ contient $M$, on pose $(\ago_M^L)^\ast=\ago_M^\ast\cap (\ago_T^L)^\ast$. On a une inclusion canonique $\ago_L^\ast\subset \ago_M^\ast$ et une décomposition
$$\ago_M^\ast=(\ago_M^L)^\ast\oplus \ago_L^\ast.$$

L'accouplement canonique
$$X^\ast(T) \times X_\ast(T) \to \ZZ$$
se prolonge linéairement en un  accouplement parfait $\ago_T^\ast\times \ago_T \to \QQ,$ avec $\ago_T=X_\ast(T)\otimes\QQ$. Plus généralement, pour $M\in \lc^G$, on note $\ago_M$ l'espace vectoriel sur $\QQ$ dual de $\ago_M^\ast$. On identifie alors $\ago_M$ au sous-espace de $\ago_T$ annulé par $(\ago_T^M)^\ast$. On note $\ago_T^M$ le sous-espace de $\ago_T$ engendré par $\Phi^{M,\vee}$. C'est aussi l'orthogonal du sous-espace $\ago_M^\ast$ de $\ago_T^\ast$. On a donc la décomposition $\ago_T=\ago_T^M\oplus \ago_M$.
\end{paragr}

\begin{paragr} Les espaces vectoriels obtenus à partir des $\QQ$-espaces du paragraphe précédent par extension des scalaires à $\RR$ sont désignés par un $a$ en lieu et place du $\ago$. Les espaces $a$ et $a^\ast$ affublés des mêmes indice et exposant sont alors en dualité. Pour tout $P\in \fc^G$, on définit le cône dans $a_T^\ast$
\begin{equation}
  \label{eq:coneap}
a_P^{G,+}=\{\chi \in a_T^{G,\ast} \ | \forall \al \in \Phi^{N_P}, \, \chi(\al^\vee)\geq 0 \text{  et  }  \forall \al  \in \Phi^{M_P}, \, \chi(\al^\vee)=0\}.\end{equation}
On notera que $a_P^{G,+}\subset a_{M_P}^{G,\ast}$.
\end{paragr}

\section{Cohomologie équivariante}\label{sec:cohequiv}

\begin{paragr} \label{S:coh-equiv1} Soit $X$ un $k$-schéma séparé et de type fini, muni d'une action algébrique d'un tore $T$. On peut alors définir la cohomologie $\ell$-adique $T$-équivariante de $X$
$$H^\bullet_T(X)=H^\bullet_T(X,\Qlb)=\bigoplus_{n\in\NN}H_{T}^{n}(X,\Qlb)$$
à savoir la cohomologie du champ algébrique $[X/T]$ à valeurs dans $\Qlb$.
La cohomologie équivariante $H_T^\bullet(X)$ est une $\Qlb$-algèbre graduée pour le cup-produit; en particulier, $H_{T}^{\bullet}(X)$ a une structure naturelle de module gradué sur $H_{T}^{\bullet}(\mathrm{Spec}(k))$, toujours pour le cup-produit.
On note pour tout entier $n$
$$\dc^n=\dc^n(T)=\mathrm{Sym}^n (X^\ast(T)\otimes \Qlb(-1))$$
et $\dc=\dc^\bullet=\oplus_{n\in \NN}\dc^n$. On a un isomorphisme, dit de Chern-Weil,
$$\dc^\bullet \to H_{T}^{\bullet}(\Spec(k))$$
d'algèbres graduées, qui double les degrés. En particulier, la cohomologie équivariante du point s'annule en degrés impairs.

La suite spectrale de Leray pour le morphisme structural $[X/T]\to [\Spec(k)/T]$ s'écrit
$$E_2^{p,q}=H^p_T(\Spec(k))   \otimes H^{q}(X)\Rightarrow H_{T}^{p+q}(X),$$ 
où $H^\bullet(X)=H^\bullet(X,\Qlb)$ désigne la cohomologie $\ell$-adique ordinaire de $X$.
\end{paragr}

\begin{paragr} \label{S:coh-equiv2} Il existe un sous-corps fini $\FF_q$ de $k$ et un schéma $X_0$ sur $\FF_q$ tels que $X=X_0\otimes_{\FF_q} k$. Soit $F_q$ le Frobenius géométrique qui est l'inverse du $\FF_q$-automorphisme de $k$ défini  par l'élévation à la puissance $q$-ième. Il agit sur $X$ via $\Id_{X_0} \times F_q$ et sur la cohomologie $\ell$-adique ordinaire de $X$
$$H^\bullet(X)=X^\bullet(X,\Qlb).$$
Cette cohomologie est pure si, par définition, pour tout entier $n$, l'espace  $H^{n}(X)$ est pur de poids $n$ au sens de Grothendieck et Deligne, c'est-à-dire si pour tout plongement 
$$\iota \ : \ \Qlb \to \CC$$ 
et toute valeur propre $\la$ de $F_q$ dans $H^{n}(X)$ on a 
$$|\iota(\la)| =q^{n/2},$$
où $|.|$ est la valeur absolue usuelle de $\CC$.
Cette propriété ne dépend pas du choix de $q$ et $X_0$.

La pureté de la cohomologie de $X$ implique que la suite spectrale de Leray ci-dessus dégénère en $E_2$. En particulier, on a un isomorphisme non-canonique en tout degré $n$
$$H^n_T(X) \simeq \bigoplus_{p+2q=n}  H^p(X)\otimes \dc^{q},$$
et la cohomologie ordinaire de $X$ s'identifie à la réduction de la cohomologie équivariante de $X$ modulo l'idéal d'augmentation $\dc^+=\oplus_{n \geq 1} \dc^n$.
\end{paragr}

\begin{paragr} \label{S:coh-equiv3}Dans la suite, nous aurons à travailler avec des $ind$-schémas projectifs $X$ (cf. par exemple \cite{Kumar}) qui sont des réunions croissantes de schémas projectifs $\cup_{n\in \NN}X_n$. Dans ce cas, la cohomologie de $X$ est définie comme une limite projective
$$H^\bullet(X)=\underset{n}{\underleftarrow{\lim}} \,  H^\bullet(X_n).$$

Il est cependant plus agréable de travailler en homologie, de façon à manipuler des limites inductives. On définit alors l'homologie et l'homologie équivariante par dualité. Par exemple, si $X$ vérifie les hypothèses du \S \ref{S:coh-equiv1}, on pose
$$H^T_\bullet(X)=\mathrm{Hom}_{\Qlb}(H_{T}^{\bullet}(X),\Qlb).$$
On pose pour tout entier $n$
$$\Sc_n(T)=\Sc_n=\mathrm{Sym}^n (X_\ast(T)\otimes \Qlb(1)),$$
et $\Sc=\Sc_\bullet=\oplus_{n\in \NN}\Sc_n$. On peut alors identifier $\Sc_\bullet$ à $H^T_\bullet(\Spec(k))$ par un isomorphisme qui double les degrés et l'action de $\dc$ sur $\Sc$ par dérivation s'identifie à l'action de $H_T^\bullet(\Spec(k))$ sur $H^T_\bullet(\Spec(k))$ par cap-produit. Pour tout $\chi\in X^\ast(T)$, on note $\partial_\chi\in \dc$ la dérivation associée.

Sous les hypothèses de pureté du \S\ref{S:coh-equiv2}, l'homologie ordinaire s'identifie au sous-module $H^T_\bullet(X)\{\dc^+\}$ de $H^T_\bullet(X)$ annulé par tous les éléments de l'idéal $\dc^+$.

\end{paragr}

\begin{paragr} \label{S:coh-equiv4} Soit $V$ une représentation algébrique de dimension finie de $T$. Soit $X$ un sous-$k$-schéma fermé $T$-stable de l'espace projectif des droites dans $V$. On suppose que la cohomologie de $X$ est pure. Le lemme suivant est une version  en homologie $\ell$-adique d'un lemme dû à Chang-Skjelbred (cf. \cite{Chang}). Notons $X_0$ l'ensemble des points fixes de $X$ sous $T$ et $X_1$ la réunion des orbites de $T$ de dimension $\leq 1$ : ce sont des ferm\'{e}s de $X$.

\begin{lemme} \label{lem:CS} On a une suite exacte
$$ H^{T}_{\bullet}(X_{1},X_{0}) \rightarrow H^{T}_{\bullet}(X_{0})\rightarrow H^{T}_{\bullet}(X)\rightarrow 0.$$
\end{lemme}

Supposons de plus que les ensembles  $\{t_1,\ldots,t_n\}$ des points fixes et $\{O_1,\ldots,O_d\}$ des orbites de dimension $1$ soient tous deux finis. Pour toute orbite $O_i$, l'algèbre de Lie du stabilisateur dans $T$ de tout point de $O_i$ est le noyau d'un caractère $\chi_i\in X^\ast(T)$, bien défini à un multiple près. L'adhérence de cette orbite dans $X$ est la réunion disjointe de $O_i$ et de deux points fixes $t_{i_0}$ et $t_{i_\infty}$. Alors la suite exacte ci-dessous s'explicite de la manière suivante
$$\begin{CD} \bigoplus_{i=1}^d \ker(\partial_{\chi_i}) @>{\beta}>> \bigoplus_{i=1}^n \Sc @>>> H^{T}_{\bullet}(X)\rightarrow 0.$$
\end{CD}$$
où la flèche de gauche $\beta$ est définie ainsi : pour tout $f\in \ker(\partial_{\chi_i})$, on pose $\beta(f)_{j_0}=f$, $\beta(f)_{j_\infty} =-f$ et $\beta(f)_{i}=0$ si $i\notin\{j_0,j_\infty\}$ (cf. \cite{GKM2} \S6,\S7, \cite{GKM}\S4).
\end{paragr}

\section{Une suite spectrale}\label{sec:suite-spec}

\begin{paragr} Soit $k$ un corps alg\'{e}briquement clos.  Par sch\'{e}ma on entend
ici un sch\'{e}ma $X$ s\'{e}par\'{e}, localement de type fini et de
dimension finie sur $k$ qui est admet un recouvrement ouvert
d\'{e}nombrable croissant $(U_{m})_{m}$ pour lequel chaque $U_{m}$ est
de type fini sur $k$.

Soit $\Lambda$ un anneau commutatif unitaire n{\oe}th\'{e}rien 
annul\'{e} par une puissance de $\ell$ et $\mathcal{F}$ un 
faisceau constructible de $\Lambda$-modules sur
un tel sch\'{e}ma $X$. On d\'{e}finit ses groupes d'homologie
$H_{q}(X,\mathcal{F})$, $q\geq 0$, par
$$H_{q}(X,\mathcal{F})=\lim_{\underset{m}{\longrightarrow}}
H_{\mathrm{c}}^{-q}(U_{m},\mathcal{F}\otimes K_{U_{m}})
$$
o\`{u} $(U_{m})_{m\geq 0}$ est n'importe quel recouvrement comme
ci-dessus et o\`{u} pour tout $k$-sch\'{e}ma s\'{e}par\'{e} de type
fini $U$ on a not\'{e} $K_{U}\in D_{\mathrm{c}}^{\mathrm{b}}
(U,\Lambda)$ son complexe dualisant.  On pose la m\^{e}me
d\'{e}finition pour un $\mathbb{Z}_{\ell}$-faisceau ou
$\mathbb{Q}_{\ell}$-faisceau $\mathcal{F}$.  Si $\mathcal{F}$ est un
$\mathbb{Z}_{\ell}$-faisceau on a comme d'habitude
$$
H_{q}(X,\mathbb{Q}_{\ell}\otimes_{\mathbb{Z}_{\ell}}\mathcal{F})
=\mathbb{Q}_{\ell}\otimes_{\mathbb{Z}_{\ell}}H_{q}(X,\mathcal{F}),
$$
par contre la fl\`{e}che canonique
$$
H_{q}(X,\mathcal{F})\rightarrow\lim_{\underset{n}{\longleftarrow}}
H_{q}(X,\mathcal{F}/\ell^{n}\mathcal{F})
$$
n'est pas en g\'{e}n\'{e}ral un isomorphisme.
\bigskip

Soit $G$ un groupe ab\'{e}lien libre de type fini de rang $d$ et $X$
un sch\'{e}ma muni d'une action propre et libre de $G$ (le morphisme
$G\times X\rightarrow X\times_{k}X,~(g ,x)\mapsto (g\cdot x,x)$, est
une immersion ferm\'{e}e).  On suppose que pour tout sous-groupe $H$
de $G$ le quotient $X/H$ existe dans la cat\'{e}gorie des sch\'{e}mas
consid\'{e}r\'{e}s ici et que le morphisme quotient $X\rightarrow X/H$
est \'{e}tale.  On note $Y=X/G$ et $f:X\rightarrow Y$ le morphisme
quotient.  On suppose de plus que pour tout ouvert $V$ de $Y$ qui est
de type fini sur $k$ il existe un ouvert $U$ de $f^{-1}(V)$ de type
fini sur $k$ tel que les $g\cdot U$ recouvrent $f^{-1}(V)$.  Il est
possible que certaines de ces hypoth\`{e}ses soient redondantes ou ne
soient pas n\'{e}cessaires.

Soit $\kappa :G\rightarrow \mathbb{Q}_{\ell}^{\times}$ un
caract\`{e}re d'ordre fini.  Le rev\^{e}tement fini \'{e}tale
galoisien $X/\mathrm{Ker}(\kappa)\rightarrow Y$ de groupe de Galois
$G/\mathrm{Ker}(\kappa)$ et le caract\`{e}re $\kappa$ d\'{e}finissent
un syst\`{e}me local $\mathcal{L}_{\kappa}$ sur $Y$.

Le but de cette section est de prouver la proposition suivante.

\begin{proposition} \label{prop:suitespec}Il existe une suite spectrale canonique
\begin{equation}\label{SuiteSpectrale}
E_{pq}^{2}=\mathrm{Tor}_{p}^{\mathbb{Q}_{\ell}[G]}(\mathbb{Q}_{\ell 
,\kappa},H_{q}(X,\mathbb{Q}_{\ell}))\Rightarrow H_{p+q}(Y,\mathcal{L}_{\kappa})
\end{equation}
o\`{u} $\mathbb{Q}_{\ell ,\kappa}$ est le 
$\mathbb{Q}_{\ell}[G]$-module de rang $1$ sur $\mathbb{Q}_{\ell}$ 
d\'{e}fini par $\kappa$.
\end{proposition}

Le reste de cette section est consacrée à la preuve de cette proposition.
\end{paragr}

\begin{paragr}
Le cas g\'{e}n\'{e}ral se d\'{e}duit du cas particulier o\`{u} $Y$ 
est de type fini par passage \`{a} la limite inductive sur les 
ouverts de type fini sur $k$ de $Y$. On peut donc supposer et on 
supposera dans la suite que $Y$ est de type finie sur $k$.
\end{paragr}

\begin{paragr}
Commen\c{c}ons par construire une suite spectrale analogue pour des
coefficients de torsion.  Soit $\Lambda$ un anneau commutatif unitaire
n{\oe}th\'{e}rien annul\'{e} par une puissance de $\ell$.

Pour chaque ouvert $U$ de $X$ on a le faisceau constructible de
$\Lambda$-modules $f_{U,!}\Lambda_{U}$ sur $Y$ o\`{u}
$f_{U}:U\rightarrow Y$ est la restriction de $f$ \`{a} $U$.  Les
ouverts $U$ consid\'{e}r\'{e}s forment un syst\`{e}me inductif pour
l'inclusion et il en est de m\^{e}me des $f_{U,!}\Lambda_{U}$.  La
limite inductive de ces faisceaux de $\Lambda$-modules existe dans la
cat\'{e}gorie de $\Lambda$-faisceaux de $\Lambda$-modules et n'est
autre que le faisceau $f_{!}\Lambda_{X}$ dont les sections sur un
ouvert \'{e}tale $V$ de $Y$ sont les fonctions $f^{-1}(V)\rightarrow
\Lambda$ localement constantes \`{a} support propre sur $V$. 
L'action de $G$ sur $X$ munit le $\Lambda$-faisceau $f_{!}\Lambda_{X}$
d'une structure de faisceau de $\Lambda [G]$-modules.  Comme
$f:X\rightarrow Y$ est localement trivial pour la topologie \'{e}tale
sur $Y$, $f_{!}\Lambda_{X}$ est localement isomorphe au faisceau
constant $\Lambda_{Y}[G]$.

On peut former le complexe $R\Gamma_{\mathrm{c}}(Y,K_{Y}
\otimes_{\Lambda} f_{!}\Lambda_{X})$ dans la cat\'{e}gorie
d\'{e}riv\'{e}e des complexes born\'{e}s de $\Lambda [G]$-modules
\`{a} cohomologie de type fini.  Comme $f$ est \'{e}tale, on a
$K_{X}=f^{!}K_{Y}=f^{\ast}K_{Y}$ et la formule des projections donne
l'\'{e}galit\'{e}
$$
H_{q}(X,\Lambda )=H_{\mathrm{c}}^{-q}(Y,K_{Y}
\otimes_{\Lambda}^{\mathrm{L}}f_{!}\Lambda_{X})
$$
pour tout entier $q$.

Soit maintenant $\chi :G\rightarrow\Lambda^{\times}$ un 
caract\`{e}re d'ordre fini.  Si $\Lambda_{\chi}$ est le
$\Lambda_{\ell}[G]$-module de rang $1$ sur $\Lambda$ d\'{e}fini par
$\chi$, le faisceau localement constant de $\Lambda$-modules libres
de rang $1$
$$
\mathcal{L}_{\chi}^{\Lambda}=\Lambda_{\chi}\otimes_{\Lambda 
[G]}f_{!}\Lambda_{X}=
\Lambda_{\chi}\otimes_{\Lambda [G]}^{\mathrm{L}}f_{!}\Lambda_{X}
$$
est la version faisceau constructible de $\Lambda$-modules du
$\mathbb{Q}_{\ell}$-faisceau $\mathcal{L}_{\chi}$.  Appliquant le
foncteur $R\Gamma_{\mathrm{c}}(Y,K_{Y}
\otimes_{\Lambda}^{\mathrm{L}}(-))$ \`{a} l'\'{e}galit\'{e} ci-dessus
on obtient un isomorphisme canonique
$$
R\Gamma_{\mathrm{c}}(Y,K_{Y}\otimes_{\Lambda}
\mathcal{L}_{\chi}^{\Lambda})=\Lambda_{\chi}\otimes_{\Lambda 
[G]}^{\mathrm{L}}R\Gamma_{\mathrm{c}}(Y,K_{Y}\otimes_{\Lambda}
f_{!}\Lambda_{X})
$$
dans la cat\'{e}gorie d\'{e}riv\'{e}e des complexes born\'{e}s de
$\Lambda$-modules \`{a} cohomologie de type fini. On en d\'{e}duit la 
suite spectrale cherch\'{e}e
\begin{equation}\label{SuiteSpectraleLambda}
E_{pq}^{2}=\mathrm{Tor}_{p}^{\mathbb{Q}_{\ell}[G]}(\Lambda_{\chi},
H_{q}(X,\Lambda))
\Rightarrow H_{p+q}(Y,\mathcal{L}_{\chi}^{\Lambda}).
\end{equation}
Tous les termes $E_{pq}^{r}$ sont des $\Lambda$-modules de type fini.
\end{paragr}

\begin{paragr}
Revenons maintenant \`{a} la situation de d\'{e}part \`{a}
coefficients dans $\mathbb{Q}_{\ell}$.  Comme $\kappa$ est d'ordre
fini, $\kappa$ est \`{a} valeurs dans $\mathbb{Z}_{\ell}^{\times}$ et,
pour chaque entier $n>0$ on peut le r\'{e}duire modulo $\ell^{n}$ en
un caract\`{e}re $\kappa_{n}:G\rightarrow
(\mathbb{Z}/\ell^{n}\mathbb{Z})^{\times}$.

On v\'{e}rifie que pour $n$ variable les suites spectrales
\ref{SuiteSpectraleLambda} avec
$\Lambda=\mathbb{Z}/\ell^{n}\mathbb{Z}$ et $\chi =\kappa_{n}$ forment
un syst\`{e}me projectif $\ell$-adique.  En passant \`{a} la limite
sur $n$ on obtient une suite spectrales de $\mathbb{Z}_{\ell}$-modules
de type fini.  Nous allons voir que cette suite spectrale
tensoris\'{e}e par $\mathbb{Q}_{\ell}$ au-dessus de
$\mathbb{Z}_{\ell}$ r\'{e}pond \`{a} la question.

Par construction de la cohomologie $\ell$-adique on a
$$
H_{p+q}(Y,\mathcal{L}_{\kappa})=\mathbb{Q}_{\ell}
\otimes_{\mathbb{Z}_{\ell}}\lim_{\longleftarrow\atop n}
H_{p+q}(Y,\mathcal{L}_{\kappa}^{\mathbb{Z}/\ell^{n}\mathbb{Z}}).
$$
Pour terminer il ne reste donc plus qu'\`{a} montrer que la 
fl\`{e}che naturelle
$$
\mathrm{Tor}_{p}^{\mathbb{Z}_{\ell}[G]}(\mathbb{Z}_{\ell
,\kappa},H_{q}(X,\mathbb{Z}_{\ell}))\rightarrow \lim_{\longleftarrow\atop n}
\mathrm{Tor}_{p}^{(\mathbb{Z}/\ell^{n}\mathbb{Z})[G]}
((\mathbb{Z}/\ell^{n}\mathbb{Z})_{\kappa},H_{q}(X,\mathbb{Z}/\ell^{n}\mathbb{Z}))
$$
est un isomorphisme pour tous entiers $p$ et $q$.

Choisissons une base $(e_{1},\dots ,e_{d})$ du $\mathbb{Z}$-module $G$
et notons $u:(\mathbb{Z}/\ell^{n}\mathbb{Z})[G]^{d}\rightarrow
(\mathbb{Z}/\ell^{n}\mathbb{Z})[G]$ la forme
$(\mathbb{Z}/\ell^{n}\mathbb{Z})[G]$-lin\'{e}aire d\'{e}finie par
$x\mapsto \sum_{i=1}^{d}([e_{i}]-\kappa (e_{i})[0])x_{i}$ o\`{u} on a
not\'{e} comme d'habitude $\sum_{g\in G}a_{g}[g]$ avec
$a_{g}\in\mathbb{Z}/\ell^{n}\mathbb{Z}$ les \'{e}l\'{e}ments de
$(\mathbb{Z}/\ell^{n}\mathbb{Z})[G]$.  Consid\'{e}rons alors le
complexe de Koszul
\begin{equation*}
\begin{split}
&K^{\mathbb{Z}/\ell^{n}\mathbb{Z}}(e_{1},\ldots ,e_{d})\\
&=(0\rightarrow \wedge^{d}(\mathbb{Z}/\ell^{n}\mathbb{Z})[G]^{d}
\rightarrow \cdots\rightarrow \wedge^{2}(\mathbb{Z}/\ell^{n}\mathbb{Z})[G]^{d}
\rightarrow (\mathbb{Z}/\ell^{n}\mathbb{Z})[G]^{d}
\rightarrow (\mathbb{Z}/\ell^{n}\mathbb{Z})[G]\rightarrow 0)
\end{split}
\end{equation*}
o\`{u} la diff\'{e}rentielle $\wedge^{r}(\mathbb{Z}/\ell^{n}\mathbb{Z})[G]^{d}
\rightarrow\wedge^{r-1}(\mathbb{Z}/\ell^{n}\mathbb{Z})[G]$ envoie
$x_{1}\wedge\cdots\wedge x_{r}$ sur
$$
\sum_{\rho =1}^{r}(-1)^{\rho -1}u(x_{\rho})x_{1}\wedge\cdots \wedge x_{\rho
-1}\wedge x_{\rho +1}\wedge\cdots \wedge x_{r}
$$
pour $r=1,\ldots ,d$. 

La fl\`{e}che $(\mathbb{Z}/\ell^{n}\mathbb{Z})[G]\rightarrow
\mathbb{Z}/\ell^{n}\mathbb{Z}_{\kappa_{n}}$ qui envoie
$\sum_{g}a_{g}[g]$ sur $\sum_{g}a_{g}\kappa (g)^{-1}$ est une
augmentation de $K^{\mathbb{Z}/\ell^{n}\mathbb{Z}}(e_{1},\ldots
,e_{d})$ et, comme $([e_{1}]-\kappa (e_{1})[0],\ldots ,[e_{d}]-\kappa
(e_{d})[0])$ est une suite r\'{e}guli\`{e}re dans
$(\mathbb{Z}/\ell^{n}\mathbb{Z})[G]$ le complexe
$K^{\mathbb{Z}/\ell^{n}\mathbb{Z}}(e_{1},\ldots ,e_{d})$ augment\'{e}
est acyclique.  Par suite
$$
\mathrm{Tor}_{p}^{\mathbb{Z}/\ell^{n}\mathbb{Z} 
[G]}((\mathbb{Z}/\ell^{n}\mathbb{Z})_{\kappa},H_{q}(X,\mathbb{Z}/\ell^{n}\mathbb{Z} ))
=H_{p}(K^{\mathbb{Z}/\ell^{n}\mathbb{Z}}(e_{1},\ldots 
,e_{d})\otimes_{\mathbb{Z}/\ell^{n}\mathbb{Z} [G]}H_{q}(X,\mathbb{Z}/\ell^{n}\mathbb{Z}))
$$
pour tout $p$.

En rempla\c{c}ant dans ce qui pr\'{e}c\`{e}de 
$\mathbb{Z}/\ell^{n}\mathbb{Z}$ par $\mathbb{Z}_{\ell}$ on a de 
m\^{e}me un complexe de Koszul $K^{\mathbb{Z}_{\ell}}(e_{1},\ldots 
,e_{d})$ et la relation
$$
\mathrm{Tor}_{p}^{\mathbb{Z}_{\ell}[G]}(\mathbb{Z}_{\ell, 
\kappa},H_{q}(X,\mathbb{Z}_{\ell}))
=H_{p}(K^{\mathbb{Z}_{\ell}}(e_{1},\ldots 
,e_{d})\otimes_{\mathbb{Z}_{\ell}[G]}H_{q}(X,\mathbb{Z}_{\ell})).
$$
On veut donc montrer que les fl\`{e}ches naturelles de complexes
$$
K^{\mathbb{Z}_{\ell}}(e_{1},\ldots 
,e_{d})\otimes_{\mathbb{Z}_{\ell}[G]}H_{q}(X,\mathbb{Z}_{\ell})
\rightarrow K^{\mathbb{Z}/\ell^{n}\mathbb{Z}}(e_{1},\ldots 
,e_{d})\otimes_{\mathbb{Z}/\ell^{n}\mathbb{Z} 
[G]}H_{q}(X,\mathbb{Z}/\ell^{n}\mathbb{Z}),
$$
o\`{u} plus explicitement les fl\`{e}ches naturelles de 
$$
C=[H_{q}(X,\mathbb{Z}_{\ell})^{d\choose d}\rightarrow 
H_{q}(X,\mathbb{Z}_{\ell})^{d\choose d-1}\rightarrow\cdots \rightarrow
H_{q}(X,\mathbb{Z}_{\ell})^{d}\rightarrow H_{q}(X,\mathbb{Z}_{\ell})]
$$
vers les
$$
C_{n}=[H_{q}(X,\mathbb{Z}/\ell^{n}\mathbb{Z})^{d\choose d}\rightarrow 
H_{q}(X,\mathbb{Z}/\ell^{n}\mathbb{Z})^{d\choose d-1}\rightarrow\cdots \rightarrow
H_{q}(X,\mathbb{Z}/\ell^{n}\mathbb{Z})^{d}\rightarrow 
H_{q}(X,\mathbb{Z}/\ell^{n}\mathbb{Z})],
$$
induisent un isomorphisme de l'homologie de la source vers la limite 
projective sur $n$ de l'homologie du but.

Soit $U$ un ouvert de $X$ de type fini sur $k$ tel que les $g\cdot U$
recouvrent $X$.  On peut calculer \`{a} la \v{C}ech l'homologie des
$X$ \`{a} partir des homologies des intersections $g_{1}\cdot U\cap
g_{2}\cdots U\cap\cdots g_{q}\cdot U$ pour $g_{1},\ldots ,g_{q}\in G$.
Comme $X$ est de dimension finie, on peut se limiter aux $q$
inf\'{e}rieurs ou \'{e}gaux \`{a} 2 fois la dimension de $X$.  Comme
chaque $H_{q}(g_{1}\cdot U\cap g_{2}\cdots U\cap\cdots g_{q}\cdot
U,\mathbb{Z}_{\ell})$ est un $\mathbb{Z}_{\ell}$-module de type fini,
il r\'{e}sulte de ce calcul \`{a} la \v{C}ech que pour chaque entier
$q$, $H_{q}(X,\mathbb{Z}_{\ell})$ est un $\mathbb{Z}_{\ell}[G]$-module
de type fini et il existe un entier $N\geq 0$ tel le noyau
$H_{q}(X,\mathbb{Z}_{\ell})[\ell^{N}]$ de la multiplication par
$\ell^{N}$ dans$H_{q}(X,\mathbb{Z}_{\ell})$ soit le sous-module de
torsion de $H_{q}(X,\mathbb{Z}_{\ell})$ tout entier.
Comme l'homologie du complexe $C$ est donc de type fini sur 
$\mathbb{Z}_{\ell}$ et que la $\ell$-torsion dans $C$ est born\'{e}e, 
le syst\`{e}me projectif de morphismes quotients de complexes
$$
C\rightarrow C/\ell^{n}C
$$
est induit des isomorphismes
\begin{equation}\label{Isomorphisme1}
H_{q}(C)\buildrel\sim\over\longrightarrow \lim_{{\longleftarrow\atop 
n}}H_{q}(C/\ell^{n}C).
\end{equation}
En effet, remarquons que, si $0\rightarrow M'\rightarrow 
M\rightarrow M''\rightarrow 0$ est une suite exacte courte de 
$\mathbb{Z}_{\ell}$-modules et si $M''[\ell^{n}]$ ne d\'{e}pend pas 
de $n$ pour $n$ assez grand, alors la suite
$$
0\rightarrow \lim_{{\longleftarrow\atop n}}M'/\ell^{n}M'\rightarrow 
\lim_{{\longleftarrow\atop n}}M/\ell^{n}M\rightarrow 
\lim_{{\longleftarrow\atop n}}M''/\ell^{n}M''\rightarrow 0
$$
est aussi exacte; pour d\'{e}montrer \ref{Isomorphisme1} il suffit 
donc d'appliquer cette remarque aux suites exactes
$$
0\rightarrow Z_{q}\rightarrow C_{q}\rightarrow B_{q}\rightarrow 0
$$
et
$$
O\rightarrow B_{q+1}\rightarrow Z_{q}\rightarrow H_{q}(C)\rightarrow 0
$$
o\`{u} $Z_{q}=\mathrm{Ker}(C_{q}\rightarrow C_{q-1})$ et
$B_{q}=\mathrm{Im}(C_{q}\rightarrow C_{q-1})$ (la $\ell$-torsion de 
$B_{q}$ est born\'{e}e puisque $B_{q}\subset C_{q}$ et que celle de 
$C_{q}$ l'est, et la $\ell$-torsion de $H_{q}(C)$ est born\'{e} 
puisque $H_{q}(C)$ est de type fini sur $\mathbb{Z}_{\ell}$).

Maintenant on a un syst\`{e}me projectif de suites exactes
$$
0\rightarrow H_{q}(X,\mathbb{Z}_{\ell})/\ell^{n}H_{q}(X,\mathbb{Z}_{\ell})
\rightarrow H_{q}(X,\mathbb{Z}/\ell^{n}\mathbb{Z})\rightarrow
H_{q-1}(X,\mathbb{Z}_{\ell})[\ell^{n}]\rightarrow 0
$$
obtenues par passage \`{a} la limite inductive \`{a} partir des suites
exactes analogues pour les ouverts de $X$ de type fini sur $k$, et 
donc on a un syst\`{e}me projectif de suites exactes de complexes
$$
0\rightarrow C/\ell^{n}C\rightarrow C_{n}\rightarrow D_{n}\rightarrow 
0
$$
o\`{u} le syst\`{e}me projectif des
$$
D_{n}=[H_{q-1}(X,\mathbb{Z}_{\ell})[\ell^{n}]^{d\choose d}\rightarrow 
H_{q-1}(X,\mathbb{Z}_{\ell})[\ell^{n}]^{d\choose d-1}\rightarrow\cdots 
\rightarrow H_{q-1}(X,\mathbb{Z}_{\ell})[\ell^{n}]]
$$
est essentiellement nul (pour tout $n$ la fl\`{e}che de transition 
it\'{e}r\'{e}e $D_{n+N}\rightarrow \cdots D_{n+1}\rightarrow D_{n}$ 
est nulle). 

On a donc un syst\`{e}me projectif de suites exactes longues
$$
\cdots\rightarrow H_{q}(C/\ell^{n}C)\rightarrow H_{q}(C_{n})\rightarrow 
H_{q}(D_{n})\rightarrow H_{q-1}(C/\ell^{n}C)\rightarrow\cdots
$$
d'o\`{u} la conclusion puisque le syst\`{e}me projectif des 
$H_{q}(D_{n})$ est essentiellement nul.
\end{paragr}

\section{Variétés toriques.}\label{sec:vartorique}

\begin{paragr} \label{S:vartor-intro} Soit $T$ un tore défini sur $k$ et $\Sigma$ un éventail dans $X^\ast(T)$. On entend par là un ensemble $\Sigma$ de cônes saillants dans $a_{T}^\ast$, de la forme 
$$\{\sum_{i=1}^N x_i \varpi_i \ | \ (x_i)_{1\leq i\leq N} \in \RR_+^N\},$$
pour une famille $(\varpi_i)_{1\leq i\leq N}\in (X^\ast(T))^N$ et un entier $N$, qui vérifie les deux propriétés suivantes :
\begin{itemize}
\item toute face d'un cône de $\Sigma$ est dans $\Sigma$ ;
\item toute intersection  de deux cônes de $\Sigma$ est dans $\Sigma$.
\end{itemize}
Soit $Y=Y_\Sigma$ la variété torique sur $\Qlb$ associée au tore $\Tc$ et à l'éventail $\Sigma$. 

À chaque cône $\sigma\in \Sigma$ est associé un ouvert $\Tc$-invariant $U_\sigma$ qui est isomorphe à $\Spec(\Qlb[X_\ast(T)\cap\sigma^\vee])$ où $\sigma^\vee$ est le cône dans $a_T$ dual de $\sigma$ et l'application $\sigma \mapsto U_\sigma$ préserve les inclusions. Les ouverts $U_\sigma$ pour $\sigma\in \Sigma$ recouvrent $Y$.
On a également une bijection qui renverse les inclusions $\sigma  \mapsto D_\sigma$ de $\Sigma$ sur l'ensemble des sous-variétés fermées de $Y$ qui sont irréductibles et $\Tc$-invariantes. La variété $D_\sigma$ est recouverte par $U_\tau \cap D_\sigma$ pour $\tau\in \Sigma$ tel que $\sigma\subset \tau$. L'intersection $U_\tau \cap D_\sigma$ est défini dans l'ouvert affine $U_\tau$ par l'idéal engendré par 
$$\{\la\in X_\ast(T)\ |\ \forall \al\in \sigma,  \ \al(\la)>0 \}.$$
La codimension de $D_\sigma$ dans $Y$ est égale à la dimension de l'espace vectoriel engendré par $\sigma$. On appelle rayon un cône de $\Sigma$ qui engendre une droite. Soit $\Sigma(1)$ le sous-ensemble de $\Sigma$ formé des rayons. Soit $\Div_{\Tc}(Y)$ le groupe des diviseurs de Weil $\Tc$-invariants : c'est le groupe abélien libre de base $D_\sigma$ pour $\sigma\in \Sigma(1)$. 

 Pour tout rayon $\sigma\in \Sigma(1)$, on note $\varpi_\sigma\in X^\ast(T)$ l'unique élément qui engendre le monoïde $\sigma\cap  X^\ast(T)$. Soit $\Qlb(Y)$ le corps des fonctions de $Y$.  Puisque la variété $Y$ est normale, toute fonction $f\in \Qlb(Y)$ définit un diviseur $(f)$. Si $f$ est la fonction définie par $\la\in  X_\ast(T)$, on pose $(\la)=(f)$ et ce diviseur est donné par 
$$(\la)=\sum_{\sigma \in \Sigma(1)} \varpi_\sigma(\la) D_\sigma.$$
Le morphisme $\la \mapsto (\la)$ s'insère dans une suite exacte de groupes
\begin{equation}
  \label{eq:suite-classe}
  X_\ast(T) \longrightarrow \Div_{\Tc}(Y) \longrightarrow \mathrm{Cl}(Y) \longrightarrow 0,
\end{equation}
où  $\mathrm{Cl}(Y)$ désigne le groupe des classes de diviseurs de $Y$ : c'est un groupe abélien de type fini. 

Soit $D\in \Div_{\Tc}(Y)$, on note $[D]$ la classe $D$ dans $ \mathrm{Cl}(Y)$. On dit que $D$ est effectif, ce que l'on note $D\geq 0$, si $D=\sum_{\sigma\in\Sigma(1)} n_\sigma D_\sigma$ avec $n_\sigma\in \NN$. On écrit aussi $D\geq D'$ pour $D-D'\geq 0$.
\end{paragr}

\begin{paragr} \label{S:coord-hom} Nous allons introduire l'anneau de Cox des coordonnées homogènes de $Y$  (cf. \cite{Cox}) qui permet de réaliser $Y$ comme un quotient catégorique d'un ouvert d'un espace affine standard par un groupe diagonalisable et qui donne en sus une description agréable des faisceaux quasi-cohérents sur $Y$.

Soit $S$ le sous-tore de $T$ défini par
$$X_\ast(S) = \ker(X_\ast(T) \to \Div_{\Tc}(Y)).$$
Le sous-$\ZZ$-module $X_\ast(S)$ est un facteur direct de $X_\ast(T)$. Soit
\begin{equation}
  \label{eq:projecteur}
  \begin{array}{lll} X_\ast(T) & \to & X_\ast(S) \\ \la &\mapsto& \la_S \end{array}
\end{equation}
un projecteur. Les constructions qui suivent dépendent de ce choix non-canonique. 

La suite (\ref{eq:suite-classe}) se complète en une suite exacte 
\begin{equation}
  \label{eq:suite-classe1}
 0 \longrightarrow X_\ast(S) \longrightarrow X_\ast(T) \longrightarrow \Div_{\Tc}(Y) \longrightarrow \mathrm{Cl}(Y) \longrightarrow 0,
\end{equation}
qui donne par dualité une suite exacte de schémas en groupes diagonalisables sur $\Qlb$
\begin{equation}
  \label{eq:suite-classe2}
1 \longrightarrow \widehat{ \mathrm{Cl}}(Y)   \longrightarrow \mathbb{G}_m^{\Sigma(1)} \longrightarrow \Tc \longrightarrow \hat{S} \longrightarrow 1,
\end{equation}
où l'on a posé
$$\widehat{ \mathrm{Cl}}(Y)=\Hom_\ZZ( \mathrm{Cl}(Y),\Qlb^\times).$$

Soit $(y_\sigma)_{\sigma\in\Sigma(1)}$ une famille d'indéterminées. On introduit ``l'anneau des coordonnées homogènes''
\begin{equation}
  \label{eq:coord-hom}
 A=\Qlb[X_\ast(S)][(y_\sigma)_{\sigma\in\Sigma(1)}]
\end{equation}
et $A'$ le localisé de  $A$ par la partie multiplicative engendrée par $y_\sigma$, pour $\sigma\in \Sigma(1)$,
$$A'=\Qlb[X_\ast(S)][(y_\sigma^{\pm 1})_{\sigma\in\Sigma(1)}].$$
On définit un morphisme injectif de $\Qlb[X_\ast(T)]$ dans $A'$ par 
\begin{equation}
  \label{eq:mor}
  \lambda  \in X_\ast(T) \mapsto y^\lambda= \la_S \prod_{\sigma \in \Sigma(1)} y_\sigma^{\varpi_\sigma(\lambda)}
\end{equation}
Pour $D=\sum_{\sigma\in \Sigma(1)}n_\sigma D_\sigma \in \Div_{\Tc}(Y)$,  on pose 
$$y^D=\prod_{\sigma\in \Sigma(1)} y_\sigma^{n_\sigma} \in A'.$$
L'anneau $A$ est alors gradué par le groupe $\mathrm{Cl}(Y)$ :
$$A=\bigoplus_{[D]\in \mathrm{Cl}(Y)} A[D]$$
où $ A[D]$  est le sous-$\Qlb$-espace de $A$ engendré par les monômes 
$$y^{\la+D}=y^\la y^D$$
lorsque $\la\in X_\ast(T)$ vérifie $(\la)+D \geq 0$.

Pour tout $\sigma\in \Sigma$, on définit un monôme
$$u_\sigma=\prod_{\tau \in \Sigma(1), \, \tau \not\subset \sigma} y_\tau,$$
et l'on note $I$ l'idéal de $A$ engendré par  les $u_\sigma$ lorsque $\sigma$ décrit $\Sigma$. L'anneau $A$ est l'anneau des fonctions régulières sur le produit $\mathbb{A}^{\Sigma(1)}\times \hat{S}$. Le schéma en groupes $\widehat{\mathrm{Cl}}(Y)$ opère sur le premier facteur via le morphisme $\widehat{\mathrm{Cl}}(Y) \to \mathbb{G}_m^{\Sigma(1)}$ et l'action évidente de $\mathbb{G}_m^{\Sigma(1)}$ sur $\mathbb{A}^{\Sigma(1)}$. L'ouvert complémentaire du fermé défini par $I$ est de la forme $V \times \hat{S}$ et le quotient $V//\widehat{\mathrm{Cl}}(Y) \times \hat{S}$ s'identifie non canoniquement à $Y$.

Pour $\sigma\in \Sigma$, l'ouvert affine $U_\sigma$ défini au paragraphe précédent s'identifie à $V_\sigma // \widehat{\mathrm{Cl}}(Y) \times \hat{S}$ où $V_\sigma\times \hat{S}$ est le complémentaire dans $\mathbb{A}^{\Sigma(1)}\times \mathbb{G}_m^{\mathbb{B}}$ du fermé défini par l'idéal engendré par $u_\sigma$. En effet, par définition, le quotient $V_\sigma //\widehat{\mathrm{Cl}}(Y) \times \hat{S}$ est isomorphe à $\Spec(A_\sigma[0])$ où $A_\sigma$ est l'anneau obtenu à partir de $A$ par localisation par la  partie multiplicative engendrée par $u_\sigma$. On vérifie que le morphisme $\la \mapsto y^\la$ induit un isomorphisme de $\Qlb[X_ \ast(T)\cap \sigma^\vee]$ sur $A_\sigma[0]$. 

 On a un foncteur exact $L \mapsto \tilde{L}$ de la catégorie des $A$-modules gradués vers la catégorie des faisceaux quasi-cohérents sur $Y$. Ce foncteur envoie $A$ sur un faisceau isomorphe au faisceau structural et tout $A$-module de type fini sur un faisceau cohérent. De plus, tout faisceau quasi-cohérent (resp. cohérent) sur  $Y$ est isomorphe à un faisceau  $\tilde{L}$ pour un $A$-module gradué $L$ (resp. et de type fini) (cf. \cite{Cox} proposition 3.1 et théorème 3.2). 

Pour un $A$-module gradué $L$, le faisceau $\tilde{L}$ sur $U_\sigma$ est le faisceau quasi-cohérent associé au $A_\sigma[0]$-module $L_\sigma[0]$ avec $L_\sigma=L\otimes_A A_\sigma$. Pour tout $D \in \Div_{\Tc}(Y)$, on note $L(D)$ le $A$-module dont la graduation est donnée par  $L(D)[D']=L[D+D']$. Dans ces conditions, on note $\tilde{L}(D)$ le faisceau qui devrait être noté $\widetilde{L(D)}$. 
\end{paragr}

\begin{paragr}[Familles $(G,M)$-orthogonales.]\  \label{S:GMfamille} Soit $G$ un groupe défini sur $k$ réductif et connexe, $T$ un sous-tore maximal de $G$ et $M$ un sous-groupe de Lévi dans $\lc^G(T)$. On considère l'éventail  dans $X^\ast(T)$ défini par
  \begin{equation}
    \label{eq:SGM}
\Sigma=\Sigma_M^G=\{ a_P^{G,+} \ | \ P\in \fc^G(M)\}.
  \end{equation}
Les cônes maximaux sont indexés par les éléments de $\pc^G(M)$ et les rayons par les sous-groupes paraboliques propres maximaux dans $\fc^G(M)$. 

Soit $P\in \fc(M)$. On pose
$$\Pi_P=\{  \varpi_\sigma \ | \ \sigma \in \Sigma(1) \text{ et } \sigma\subset a_P^{G,+}  \}.$$
Les éléments de $\Pi_P$ forment une base de l'espace vectoriel $(a_M^G)^\ast$ et on note 
$$\Pi_P^\vee=\{\alpha_\sigma^\vee  \ | \ \sigma \in \Sigma(1) \text{ et } \sigma\subset a_P^{G,+}  \}$$ 
la base duale dans $a_M^G$. Soit $B$ un sous-groupe de Borel inclus dans $P$ et $\Delta_B$ l'ensemble des racines simples dans $B$. L'ensemble $\Delta_B\cap \Phi^{N_P}$ est en bijection avec l'ensemble des sous-groupes paraboliques maximaux contenant $P$. Soit $\al\in \Delta_B\cap \Phi^{N_P}$ ; le sous-groupe parabolique maximal $Q$ associé est défini de manière unique par $Q\supset P$ et $\Delta_B \cap \Phi^{N_Q}=\{\al \}$.  Soit $\sigma=a_Q^{G,+}$. Avec ces notations, le vecteur $\al_\sigma^\vee$ est égal à la projection de $\al^\vee$ sur $a_M^G$ à un coefficient rationnel positif près. Soit $(\varpi_\beta)_{\beta \in \Delta_B}$ la base de $(a_T^G)^\ast$ duale de la base $(\beta^\vee)_{\beta\in \Delta_B}$ de $a_T^G$. Les vecteurs $\varpi_\al$ et $\varpi_\sigma$ sont égaux à un coefficient  rationnel positif près.

Deux éléments $P$ et $P'$ de $\pc(M)$ sont adjacents si les chambres associées ont un mur commun. Ce dernier est alors porté par un hyperplan de $(a_M^G)^\ast$ d'équation $\chi(\al^\vee_{P,P'})=0$ pour un unique élément  $\al^\vee_{P,P'} \in \Pi_P^\vee$. 

Soit une famille $(\la_P)_{P\in \pc(M)}$ d'éléments de $a_M^G$. Suivant la terminologie d'Arthur (cf. \cite{dis_series} en particulier), on dit que cette famille est $(G,M)$-\emph{orthogonale} (resp.   $(G,M)$-\emph{orthogonale positive}) si pour toute paire $(P,P')$ de sous-groupes paraboliques adjacents de $\pc(M)$, il existe $x\in \RR$ (resp.  $x\in \RR_+$) tel que
$$\la_P -\la_{P'}= x\, \al^\vee_{P,P'}.$$ 

Toute famille $(G,M)$-orthogonale  $(\la_P)_{P\in \pc(M)}$ s'étend en une famille  $(\la_Q)_{Q\in \fc(M)}$ où, par définition, $\la_Q$ est le projeté sur $a_{M_Q}^G$, relativement à la décomposition 
$$a_T=a_T^{M_Q}\oplus a_{M_Q}^G \oplus a_G,$$
de $\la_P$ pour tout $P\in \pc(M)$ tel que $P\subset Q$.

\end{paragr}

\begin{paragr} \label{S:dla} Soit $B\in \pc^G(T)$ et $P=MN_P$ un sous-groupe parabolique de $G$ qui contient $B$. La projection de $a_{T}^\ast$ sur $a_{M}^\ast$ parallèlement à $(a_T^M)^\ast$ envoie bijectivement l'ensemble $\Delta_B\cap \Phi^{N_P}$ sur un ensemble noté $\Delta_P$. Ce dernier ne dépend pas du choix de $B$. 

Soit $\la$ une famille $(G,M)$-orthogonale. On pose
$$d(\la)= \min \{\al(\la_P) \ | \    P\in \pc^G(M),\al\in \Delta_P\}$$
et
$$\delta(\la)=\min \{ \al(\la_P-\la_{P'}) \ | \    P, P'\in \pc^G(M) \text{ adjacents  et }\al\in \Delta_P\cap (-\Delta_{P'})\}.$$
Il résulte de la définition que $\delta(\la)$ est positif si et seulement si la famille $\la$ est positive.  On dit que la famille $\la$ est très positive si $d(\la)\geq 0$. Il résulte de la majoration évidente 
$$\delta(\la)\geq 2d(\la)$$
que très positif implique positif.

 Soit $\la\in X_\ast(T)$ et $P\in \fc^G(T)$.  Soit $\la_P$ l'élément de $a_{M_P}^G$ défini de la façon suivante. Soit $B\in \pc^G(T)$ inclus dans $P$ et $c_B^+$ la chambre de Weyl correspondante dans $a_T^G$. Soit $\la'$ l'unique élément de $c_B^++a_G$ qui soit dans l'orbite de $\la$ sous le groupe de Weyl $W^G(T)$ ; la projection de $\la'$ sur $a_{M_P}^G$ parallèlement à $a_T^{M_P}\oplus a_G$ ne dépend pas du choix de $B$ tel que $B\subset P$. Par définition, cette projection est $\la_P$. 

Pour tout $M\in \lc^G(T)$, la famille $(\la_P)_{P\in \pc(M)}$ est $(G,M)$-orthogonale. Posons
\begin{equation} 
  \label{eq:d}
 d(\la)= d((\la_B)_{B\in \pc(T)}).
\end{equation}
On a l'inégalité (cf. \cite{localtrace} l.(3.2) p. 20)
$$d(\la)\leq d((\la_P)_{P\in \pc(M)}).$$
Énonçons le lemme suivant qui résulte immédiatement de l'inégalité ci-dessus.
 
\begin{lemme} \label{lem:positivite} 
Pour toute famille $(G,M)$-orthogonale $(\mu_P)_{P\in \pc(M)}$, il existe un entier $d\in \NN$ de sorte que pour tout $\la \in X^\ast(T)$ tel que $d(\la)\geq d$, la famille $(\la_P+\mu_P)_{P\in \pc(M)}$ est  $(G,M)$-orthogonale très positive donc positive.
\end{lemme}
\end{paragr}

\begin{paragr}[Diviseurs et familles $(G,M)$-orthogonales.]\ \label{S:div-GM} --- Soit  $\Sigma_M^G$ l'éventail défini en (\ref{eq:SGM}).  

\begin{definition} \label{def:adapte} On appelle éventail $(G,M)$-adapté un éventail $\Sigma$ dans $X^\ast(T)$ qui vérifie la condition suivante : pour tout cône $\sigma \in \Sigma$, il existe $\sigma' \in \Sigma_M^G$ tel que  $\sigma\subset  \sigma'+a_G^\ast$ c'est-à-dire il existe $P\in \fc^G(M)$ tel que $\sigma\subset  a_P^{G,+}+a_G^\ast$.
\end{definition}

Soit  $\Sigma$ un éventail $(G,M)$-adapté. 

\begin{lemme} Soit $\la$ une famille $(G,M)$-orthogonale. Soit $\sigma \in \Sigma(1)$ un rayon et $P$ et $Q$ deux sous-groupes paraboliques dans $\fc^G(M)$ tels que $\sigma \subset a_P^{G,+} \cap a_Q^{G,+}+ a_G^\ast$. Alors on  a 
$$\varpi_\sigma(\la_P)=\varpi_\sigma(\la_Q).$$
\end{lemme}

\begin{preuve}  Puisqu'il existe $R\in \fc^G(M)$ tel que  $a_R^{G,+}=a_P^{G,+} \cap a_Q^{G,+}$, il suffit de prouver le lemme dans le cas où $P\subset Q$, cas que nous considérons désormais. En écrivant  $\la_P$ et $\la_Q$ comme projetés du même élément $\la_{P_0}$ pour $P_0\in\pc(M)$ tel que $P_0\subset P$ (cf. la fin du \S  \ref{S:GMfamille}), on voit que $\la_P -\la_Q \in a_T^{M_Q}$. Comme $\varpi_\sigma\in  a_{Q}^{G,+}+a_G^\ast$, on a $\varpi_\sigma\in a_{M_Q}^\ast$. L'égalité cherchée est alors évidente.
\end{preuve}

Soit $Y$ la variété torique associée à l'éventail $\Sigma$. On pose 
$$  \Div_{\Tc}(Y)_\RR=\Div_{\Tc}(Y)\otimes_\ZZ \RR.$$
On définit alors une application linéaire du $\RR$-espace vectoriel des familles $(G,M)$-orthogonales dans $  \Div_{\Tc}(Y)_\RR$ par 
\begin{equation}
  \label{eq:appliGM}
\la \mapsto D_\la=D_\la^\Sigma=\sum_{\sigma \in \Sigma(1)} \varpi_\sigma(\la_{Q_\sigma}) D_\sigma
\end{equation}
où la famille $(Q_\sigma)_{\sigma \in \Sigma(1)}$ est une famille d'éléments de $\fc(M)$ qui vérifie $\sigma\subset  a_{Q_\sigma}^{G,+}+a_G^\ast$. Le lemme précédent montre que $D_\la$ ne dépend pas du choix de cette famille.

Lorsque $\Sigma=\Sigma_M^G$, l'application (\ref{eq:appliGM}) est bijective et son inverse envoie le diviseur  
$$D=\sum_{\sigma\in \Sigma(1)} n_\sigma D_\sigma \in \Div_{\Tc}(Y)_\RR$$
sur la   famille  $(G,M)$-orthogonale $\mu^D$ définie pour tout $P\in \pc(M)$ par 
$$\mu_P^D= \sum_{ \alpha^\vee_\sigma \in \Pi_P^\vee} n_\sigma \alpha^\vee_\sigma.$$

Si $\la\in X_\ast(T)$, on pose $D_\la=D_{(\la_P)_{P\in \pc^G(M)}}$ où $(\la_P)_{P\in \pc^G(M)}$ est la famille $(G,M)$-orthogonale associée à $\la$ (cf. la paragraphe précédent).
\end{paragr}

\begin{paragr}\label{S:polyedre}  Désormais soit $\Sigma=\Sigma_M^G$ l'éventail défini en (\ref{eq:SGM}) et $Y=Y_\Sigma$ la variété torique associée à $\Sigma$. Soit $D\in \Div_{\Tc}(Y)_\RR$. On définit une fonction $\psi^D$ continue sur $(a_M^G)^\ast$ et linéaire par morceaux de la façon suivante : pour tout $P\in \pc(M)$ et tout $\chi \in a_P^{G,+}$, on pose
$$\psi^D(\chi)=\chi(\mu_P^D).$$  

Finalement, on définit un polyèdre dans $(a_M^G)^\ast$ par 
$$\mathfrak{P}(D)=\{ \la \in a_M^G \ | \ \forall \chi\in (a_M^G)^\ast\   \chi(\la)\leq \psi^D(\chi) \}.$$

\begin{lemme} \label{lem:cvx}Soit $D\in  \Div_{\Tc}(Y)_\RR$. Les assertions suivantes sont équivalentes :
  \begin{enumerate}
  \item la fonction $\psi^D$ est convexe ;
 \item le polyèdre $\mathfrak{P}(D)$ est l'enveloppe convexe des points $\mu_P^D$ pour $P\in\pc(M)$.
  \end{enumerate}
 \end{lemme}

\begin{preuve}  Pour alléger les notations, on omet la lettre $D$. Notons $\mathfrak{C} $ l'enveloppe convexe des points $\mu_P $ pour $P\in\pc(M)$. On a toujours $\mathfrak{P} \subset \mathfrak{C} $. En effet, soit $\la \in a_M^G$ tel que $\la\notin\mathfrak{C} $ ; il existe donc $\chi\in (a_M^G)^\ast$ tel que 
$$\chi(\la) > \sup_{\mu\in \mathfrak{C} } \chi(\mu).$$
Il est clair alors que $\chi(\la)> \psi (\chi)$ et donc $\la \notin \mathfrak{P} $.
On est donc ramené à prouver que l'assertion 1 est équivalente à
\begin{equation}
  \label{eq:2'}
\mathfrak{C} \subset \mathfrak{P} .
\end{equation}

La fonction $\psi$ est convexe si, par définition, pour tout $t\in[0;1]$ et tous $\chi$ et $\chi'$ dans $(a_M^G)^\ast$, on a
\begin{equation}
  \label{eq:cvx}
  \psi (t\chi + (1-t)\chi')\leq t\psi (\chi)+(1-t) \psi (\chi').
\end{equation}
Soit $Q$ dans $\pc(M)$ de sorte que $t\chi + (1-t)\chi'$ appartienne à $a_Q^{G,+}$. L'inégalité (\ref{eq:cvx}) se récrit 
\begin{equation}
  \label{eq:cvx2}
  t (\chi(\mu_Q)-\psi(\chi)) +(1-t) (\chi(\mu_Q)-\psi(\chi'))\leq 0.
\end{equation}
On voit alors que l'assertion 1 est équivalente à l'inégalité $\chi(\mu_Q)\leq \psi(\chi)$ pour tout $Q\in \pc(M)$ elle-même équivalente à l'inclusion $\mathfrak{C} \subset \mathfrak{P}$, puisque $\mathfrak{P}$ est convexe.
\end{preuve}

\begin{lemme} \label{lem:famille-positive} Soit $D\in  \Div_{\Tc}(Y)$. Les assertions du lemme \ref{lem:cvx} sont vérifiées dès que la famille $(\mu_P^D)_{P\in \pc(M)}$ est positive.
\end{lemme}

\begin{preuve} Pour alléger, on omet la lettre $D$. Reprenons les notations de la preuve du  lemme \ref{lem:cvx} : on y a vu qu'on a toujours $\mathfrak{P}\subset \mathfrak{C}$. Il s'agit donc de voir que  $\mu_P$ appartient à $\mathfrak{P}$ pour tout $P\in \pc(M)$. Or $\mu_P$ appartient à $\mathfrak{P}$ si et seulement si on a
$$\chi(\mu_{P'} - \mu_{P})\geq 0$$
pour tous $P'\in \pc(M)$ et $\chi\in a_{P'}^{G,+}$. On montre facilement cette dernière inégalité en partant du cas (évident) où $P$ et $P'$ sont adjacents (cf. \cite{dis_series} l.(3.2) p.217).

\end{preuve}

\end{paragr}

\begin{paragr}\label{S:cohosup} Soit $\la\in X_\ast(T)$ et $D_\la$ le diviseur de $\Div_{\Tc}(Y)$ associé par l'application (\ref{eq:appliGM}) à la famille $(G,M)$-orthogonale $(\la_P)_{P\in\pc(M)}$ définie au \S\ref{S:dla}. Dans ce même paragraphe, on a défini un entier $d(\la)$.

\begin{theoreme} \label{thm:cohosup} Pour tout $A$-module gradué $L$ de type fini, il existe un entier $d\in \NN$ tel que pour tout $\la\in X_\ast(T)$ tel que $d(\lambda)\geq d$, on ait
  \begin{enumerate}
  \item $H^0(Y,\tilde{L}(D_\la))=L[D_\la]$ ;
\item pour tout entier $i>0$, $H^i(Y,\tilde{L}(D_\la))=0$.
  \end{enumerate}
\end{theoreme}

\begin{preuve} On commence par traiter le cas du $A$-module $A(D)$ pour $D\in \Div_{\Tc}(Y)$. On vérifie directement que le morphisme naturel
$$A[D] \to H^0(Y,\tilde{A}(D))$$
est un isomorphisme. 

 Les éléments $y^{\la+D}$, pour $\la \in X_\ast(T)$, forment une base du $\Qlb$-espace vectoriel $A'[D]$. On munit $A'[D]$ d'une action de $\Tc$ en posant
$$t.y^{\la+D}=\la(t) y^{\la+D}$$
pour tout $t\in \Tc$ et $\la\in X_\ast(T)$. Pour tout $P\in \pc(M)$, soit $U_P$ l'ouvert affine de $Y$ associé au cône $a_P^{G,+}$. Soit  $\{P_1,\ldots,P_k\}\subset \pc(M)$ une partie non vide. L'espace vectoriel 
$$ H^0(U_{P_1}\cap \ldots \cap U_{P_k}, \tilde{A}(D))$$
est un sous-espace de $A'[D]$ stable par $\Tc$. Soit 
$$ H^0(U_{P_1}\cap \ldots \cap U_{P_k}, \tilde{A}(D))_\la$$
le sous-espace de  poids $\la\in X_\ast(T)$. Cet espace est non  nul si et seulement si pour tout $\sigma\in \Sigma(1)$ tel que $\sigma\subset a_{P_1}^{G,+} \cap \ldots \cap a_{P_k}^{G,+}$, 
$$\varpi_\sigma(\la)+n_\sigma \geq 0,$$
c'est-à-dire si et seulement si 
$$b_\la \cap a_{P_1}^{G,+} \cap \ldots \cap a_{P_k}^{G,+}=a_{P_1}^{G,+} \cap \ldots \cap a_{P_k}^{G,+}$$
où l'on a posé
$$b_\la=b_\la^D=\{\chi\in (a_M^G)^\ast \ | \ \chi(\la)+\psi^D(\chi)\geq 0\}.$$
Dans ce cas, $ H^0(U_{P_1}\cap \ldots \cap U_{P_k}, \tilde{A}(D))_\la$ est de dimension $1$ engendré par $y^{\la+D}$. Pour tout $\la\in X_\ast(T)$, on a donc un isomorphisme
\begin{equation}
  \label{eq:isom-cohom}
  H^0(U_{P_1}\cap \ldots \cap U_{P_k}, \tilde{A}(D))_\la \simeq H^0_{b_\la \cap a_{P_1}^{G,+} \cap \ldots \cap a_{P_k}^{G,+}}(a_{P_1}^{G,+} \cap \ldots \cap a_{P_k}^{G,+},\Qlb).
\end{equation}
La cohomologie du faisceau $\tilde{A}(D)$ se calcule comme la  cohomologie du complexe de Cech associé au faisceau  $\tilde{A}(D)$ et au recouvrement de $Y$ par les ouverts $U_P$ pour $P\in \pc(M)$. Pour tout entier $i$, le $\Qlb$-espace vectoriel  $H^i(Y,\tilde{A}(D))$ est donc muni d'une action de $\Tc$ et il résulte de (\ref{eq:isom-cohom}) que le sous-espace de poids $\la$  vérifie
 \begin{equation} \label{eq:iso-Fulton}
 H^i(Y,\tilde{A}(D))_{\la}\simeq H^i_{b_\la}((a_M^G)^\ast,\Qlb).
 \end{equation}
L'argument, qu'on ne reproduit pas ici, est celui donné par Fulton dans la preuve de la proposition du \S3.5 du livre  \cite{Fulton} (le lecteur notera que notre fonction $\psi$ est notée $-\psi$ dans ce livre.)

On en déduit alors le fait suivant : si la famille $(\mu_P^D)_{P\in \pc(M)}$ est positive, on a pour tout entier $i>0$
$$H^i(Y,\tilde{A}(D))=0.$$  
En effet, la positivité de la famille entraîne que la fonction $\psi^D$ est convexe (cf. lemmes \ref{lem:cvx} et \ref{lem:famille-positive}) et donc que l'ouvert $a_M^G-b_\la$ est convexe pour tout $\la\in X_\ast(T)$. L'isomorphisme (\ref{eq:iso-Fulton}) ci-dessus donne l'annulation voulue.

Soit $d$ un entier qui vérifie l'assertion suivante : pour tout $\la$ tel que $d(\la)\geq d$, la famille $(G,M)$-orthogonale $(\mu_P^{D+D_\la})_{P\in\pc(M)}$ est positive. Comme $\mu_P^{D+D_\la}=\mu_P^D +\la_P$, un tel entier existe (cf. lemme \ref{lem:positivite}). On a donc 
$$H^i(Y,\tilde{A}(D+D_\la))=0.$$
qui est l'assertion 2 pour le $A$-module $A(D)$.

Passons ensuite au cas général d'un $A$-module $L$ gradué de type fini. On commence par prouver l'assertion 2. Soit $r+1$ le nombre de cônes maximaux dans $\Sigma$. Alors $Y$ est recouvert par $r+1$ ouverts affines et donc l'assertion 2 est vraie pour tout $i>r$ et tout faisceau cohérent. On raisonne alors par récurrence. Soit $i>1$ telle que l'assertion 2 soit vraie pour tout  $A$-module gradué de type fini. Soit $L$ un $A$-module gradué de type fini. Il existe donc une suite exacte courte
$$\begin{CD} 0 @>>> L_2 @>>> L_1 @>>> L @>>>0
\end{CD}
$$
où $L_1$ et $L_2$ sont des $A$-modules gradués de type fini et $L_1$ est libre. En considérant la suite exacte de faisceaux associée et la suite exacte longue de cohomologie qui en résulte, on extrait la suite exacte
\begin{equation}\label{eq:presentation}
\begin{CD} H^{i-1}(Y,\tilde{L}_1(D_\la))@>>>   H^{i-1}(Y,\tilde{L}(D_\la))@>>>  H^i(Y,\tilde{L}_2(D_\la)).
\end{CD}
\end{equation}

Par hypothèse de récurrence, il existe $d_2$ tel que si $d(\la)\geq d_2$ le troisième terme est nul. Comme $L_1$ est isomorphe à une somme directe finie de modules $A(D)$ pour $D\in \Div_{\Tc}(Y)$, on sait, d'après ce qui précède, qu'il existe un entier $d_1$ tel que $d(\la)\geq d_1$ implique $H^{i-1}(Y,\tilde{L}_1(D_\la))=0$. On en déduit la nullité de  $H^{i-1}(Y,\tilde{L}(D_\la))$ pour $d(\la)\geq \max(d_1,d_2)$ ce qui prouve l'assertion 2 pour $i-1$ et achève la récurrence.

\medskip

Prouvons maintenant l'assertion $1$ pour le $A$-module $L$. Il existe une filtration finie 
$$L_0=\{0\} \subset L_1 \subset \ldots \subset L_n=L$$
telle que  chaque quotient $L_{i+1}/L_i$ soit isomorphe à un $A$-module de la forme $A/\mathfrak{p}(D)$ où $\mathfrak{p}$ est un idéal premier homogène de $A$ et $D\in \Div_{\Tc}(Y)$. On est alors immédiatement ramené au cas où $L=A/\mathfrak{p}(D)$. Considérons alors le diagramme commutatif 
$$\begin{CD}
0 @>>> \mathfrak{p}[D+D_\la] @>>> A[D+D_\la] @>>> A/ \mathfrak{p}[D+D_\la] @>>> 0\\
@. @VVV @VVV @VVV \\
0 @>>> H^0(Y,\tilde{ \mathfrak{p}}(D+D_\la)) @>>> H^0(Y,\tilde{A}(D+D_\la)) @>>> H^0(Y,\widetilde{A/ \mathfrak{p}}(D+D_\la)) @>>>0
\end{CD}$$
Les flèches verticales sont les morphismes naturels. Comme $A/ \mathfrak{p}$ est intègre, la dernière flèche verticale est injective. Il s'agit de prouver qu'elle est surjective. D'après ce qui précède, il existe $d$ tel que $d(\la)\geq d$ implique la nullité de $H^1(Y,\tilde{ \mathfrak{p}}(D+D_\la))$ et donc l'exactitude de la suite horizontale du bas. On a vu que la flèche $A[D+D_\la]\to   H^0(Y,\tilde{A}(D+D_\la))$ est toujours bijective. On en déduit la surjectivité voulue.
\end{preuve}

\end{paragr}

\section{Grassmanniennes affines tronquées}\label{sec:gat}

\begin{paragr} 
 Soit $G$ un groupe réductif connexe  défini sur $k$ et $K=G(\of)$. Soit $\Xgo^G=G(F)/K$ la grassmannienne affine. On rappelle que $\Xgo^G$ possède une structure naturelle de ind-schéma projectif sur $k$ (cf. par exemple \cite{Kumar}) i.e. il existe une suite croissante $\Xgo_0\subset \Xgo_1 \subset \ldots$ de schémas projectifs sur $k$ de sorte que $\Xgo^G=\bigcup_{n\in \NN}\Xgo_n$ et les inclusions $\Xgo_n\subset \Xgo_{n+1}$ sont des immersions fermées. La grassmannienne affine est alors munie de la topologie induite de sorte qu'une partie de $\Xgo^G$ est fermée si et seulement sa trace sur $ \Xgo_n$ est fermée dans $\Xgo_n$ pour la topologie de Zariski. 
\end{paragr}

\begin{paragr}
Soit $H_G$ l'unique application de $G(F)$ dans $\ago_G$ qui vérifie : pour tous $\chi\in X^\ast(G)$ et  $g\in G(F)$,
$$\chi(H_G(g))=\val(\chi(g)).$$
L'application $H_G$ est invariante à droite par $K$. Elle induit donc une application $\Xgo^G \to \ago_G$ encore notée $H_G$.

Il est bien connu que les  composantes connexes de $\Xgo^G$ sont les translatés sous $G(F)$ de l'image du morphisme naturel $\Xgo^{G_{\scnx}} \to \Xgo^G$. La proposition suivante montre que l'application $H_G$ est constante sur les composantes connexes de $\Xgo^G$.

\begin{lemme} \label{lem:HG} Soit $g\in G(F)$. On a $H_G(g)=0$ si et seulement si $gK$ appartient à l'image du morphisme naturel $\Xgo^{G_{\der}} \to \Xgo^G$.
\end{lemme}

\begin{preuve} Soit $g\in G(F)$ et $T$ un sous-tore maximal de $G$ défini sur $k$. Soit $L$ une extension finie de $F$,  $t\in T(L)$ et $g_1\in G_{\der}(L)$ tels que $g=g_1 t$. Soit $T_{\der}=T\cap G_{\der}$ : c'est un sous-tore maximal de $G_{\der}$. Comme $H^1(L/F,T_{\der})= 1$, on peut prendre $t\in T(F)$ et $g_1\in G_{\der}(F)$. Soit $\la\in X_\ast(T)$ et $u\in T(\of)$ tels que $t=\eps^\la u$. Si $H_G(g)=0$, on a $\la\in X_\ast(T_{\der})$ et $gK$ appartient à l'image de  $\Xgo^{G_{\der}} \to \Xgo^G$. La réciproque est évidente.
\end{preuve}
\end{paragr}

\begin{paragr} Soit $P$ un sous-groupe parabolique de $G$ défini sur $k$. Soit $P=M N$ une décomposition de Lévi définie sur $k$. L'injection $M(F) \to G(F)$ induit une immersion fermée $\Xgo^M \to \Xgo^G$. On définit une rétraction (non-algébrique)
$$\pi_P \ : \ \Xgo^G \to \Xgo^M$$
en posant $\pi_P(xk)= m M(\of)$ pour tout $x\in G(F)$ qui s'écrit $x=nmk$ suivant la décomposition d'Iwasawa $G(F)=N(F)M(F)K$.  

Soit $H_P=H_M\circ \pi_P$. C'est une application  de $\Xgo^G$ dans $\ago_M$. On note encore $H_P$ l'application définie sur $G(F)$ obtenue par composition avec l'application canonique $G(F)\to \Xgo^G$.

\begin{lemme}Soit $X\in \ago_M$. La partie $H_P^{-1}(X)$ de $\Xgo^G$ est localement fermée. La rétraction $\pi_P$ induit  un morphisme de ind-schéma de  $H_P^{-1}(X)$ dans  $\Xgo^M$.
\end{lemme}

\begin{preuve} Soit $X\in \ago_M$. Si $H_P^{-1}(X)$ est vide, il n'y a rien à prouver. Sinon il existe $g\in G(F)$ tel que $H_P(g)=X$ et par  translation à gauche par $g$ on se ramène au cas où $X=0$. Soit $P_{\der}$ le groupe dérivé de $P$. On note $\Xgo^{P_{\der}}=P_{\der}(F)/P_{\der}(\of)$ la grassmannienne affine de $P_{\der}$ : c'est un ind-schéma, limite inductive de schémas de type fini. Comme le quotient $G/ P_{\der}$ est quasi-affine, l'injection canonique $P_{\der}\to G$ induit un plongement localement fermé  $\Xgo^{P_{\der}}\to \Xgo^G$ (cf. \cite{BB} théorème 4.5.1). D'après le lemme \ref{lem:HG}, l'image de ce plongement est précisément $H_P^{-1}(X)$.  Le groupe $M_{\der}(\of)$ opère sur le produit $M_{\der}(F)\times N(F)/N(\of)$ par $\gamma.(m,n)=(m\gamma^{-1},\gamma n)$. La structure naturelle de ind-schéma du ``quotient'' $M_{\der}(F)\times^{M_{\der}(\of)} N(F)/N(\of)$ (cf. \cite{Gaitsgory} A.3) donne un isomorphisme de ce quotient sur $\Xgo^{P_{\der}}$. Via cet isomorphisme, l'application  $\Xgo^{P_{\der}} \to \Xgo^{M}$ obtenue par composition de $\pi_P$ avec le plongement $\Xgo^{P_{\der}}\to \Xgo^G$ n'est autre que le morphisme canonique   $M_{\der}(F)\times^{M_{\der}(\of)} N(F)/N(\of)\to \Xgo^{M}$. Cela permet de conclure.
\end{preuve}

\end{paragr}

\begin{paragr} Soit $T\subset G$ un sous-tore maximal. Soit $M\in \lc^G(T)$ et $\Sigma$ un éventail $(G,M)$-adapté (cf. définition \ref{def:adapte}). Soit $Y=Y_\Sigma$ la variété torique associée. Pour tout $x\in G(F)$, la famille $H(x)$ obtenue par projection sur $\ago_T^G$ (parallèlement à $\ago_G$) de la famille $(H_P(x))_{P\in \pc(M)}$ est une famille $(G,M)$-orthogonale positive (cf. \cite{dis_series} lemme 3.6). En utilisant la définition (\ref{eq:appliGM}), on obtient  un  élément de $\Div_{\Tc}(Y)$ 
 \begin{equation}
   \label{eq:DMG}
   D_M^G(x)=D_{M}^{G,\Sigma}(x)=D_{H(x)}^\Sigma.
\end{equation}
Pour $Q\in \fc(M)$ et $P\in \pc(M)$ tels que $P \subset Q$, l'élément $H_Q(x)$ est la projection de $H_P(x)$ sur $\ago_{M_Q}^G$ parallèlement à $\ago_T^{M_Q}\oplus \ago_G$. On obtient alors l'expression suivante
  \begin{equation}
   \label{eq:DMGII}
   D_M^G(x)=\sum_{\sigma \in \Sigma(1)} \varpi_\sigma(x)  D_\sigma
\end{equation}
où l'on a posé
$$\varpi_\sigma(x)=  \varpi_\sigma(H_{Q_\sigma}(x))$$
et où la famille $(Q_\sigma)_{\sigma\in \Sigma(1)}$ vérifie les conditions qui suivent (\ref{eq:appliGM}).

On remarque que l'application $x \mapsto D_M^G(x)$ est invariante à droite par $K$ et à gauche par $T$ ainsi que par $S(F)$ où $S$ est le sous-tore de $T$ défini au \S\ref{S:coord-hom}.

\begin{definition} \label{def:grafftr} Soit $D\in \Div_{\Tc}(Y_\Sigma)$. La grassmannienne affine tronquée par le diviseur $D$ est définie par
$$\Xgo^G(D)=\Xgo^G(\Sigma,D)=\{x\in \Xgo^G \ | \  D_M^G(x)\leq D\}.$$
\end{definition}

\begin{proposition} Pour tout $D\in \Div_{\Tc}(Y_\Sigma)$, la grassmannienne affine tronquée $\Xgo^{G}(D)$ est un sous-ind-schéma fermé de $\Xgo^G$.
\end{proposition}

\begin{preuve} Pour tout $\sigma \in \Sigma(1)$, soit $P_\sigma$ un sous-groupe parabolique  tel que $\sigma\subset a_{P_\sigma}^{G,+}+a^\ast_G$ ; cela existe  puisque $\Sigma$ est $(G,M)$-adapté. Par conséquent $\varpi_\sigma$ est un caractère de $M_{P_\sigma}$ qui est $P_\sigma$-dominant. Soit $(\rho_\sigma,V_\sigma)$ une représentation algébrique de $G$ sur le $k$-espace vectoriel de dimension finie $V_\sigma$ de plus haut poids $\varpi_\sigma$ (relatif à l'ordre défini par $P_\sigma$). Soit $v_\sigma \in V_\sigma$ un vecteur de plus haut poids $\varpi_\sigma$.  Notons que si $x \in G(F)$ s'écrit $x=nmk$ suivant la décomposition d'Iwasawa $N_{P_\sigma} M_{P_\sigma}K$ on a 

\begin{equation}
  \label{eq:actionsurpoids}
  \rho(x^{-1})v_\sigma=\rho(k)^{-1} (\varpi_\sigma (m)^{-1}v_\sigma).
\end{equation}
On pose alors $(\rho,V)=\bigoplus_{\sigma\in \Sigma(1)} (\rho_\sigma,V_\sigma)$. On définit $V(F)=V\otimes_k F$ et $V(\of)=V\otimes_k\of$. Notons que $V(\of)$ est stable par $K$.  On associe au diviseur $D=\sum_{\sigma\in \Sigma(1)}n_\sigma D_\sigma$ un vecteur $v\in V(F)$ défini par 
$$v=\bigoplus_{\sigma\in\Sigma(1) } \eps^{n_\sigma} v_\sigma.$$
Considérons alors la fibre de Springer affine ``généralisée''
$$\Xgo^G_v=\{x\in G(F)/ K  | \rho(x)^{-1}v \in V(\of) \}.$$
On voit que c'est un sous-ind-schéma fermé de $\Xgo^G$. Comme $V(\of)$ est stable par $K$, il résulte de la ligne (\ref{eq:actionsurpoids}) que $\Xgo^{G}(D)=\Xgo^G_v$. Cela donne le résultat.
 \end{preuve}

\end{paragr}

\begin{paragr} On continue avec les hypothèse du paragraphe précédent. Soit $H$ un groupe réductif connexe défini sur $F$ de sorte que son groupe dual $\Hc$ s'identifie à un sous-groupe fermé de $\Gc$ qui contient $\Tc$. Le groupe $\Mc\cap \Hc$ est un sous-groupe de Lévi de $\Hc$ qui contient $\Tc$. Dualement, on obtient un sous-groupe de Lévi $M_H$ de $H$.

\begin{proposition} \label{eventail-adapte}
L'éventail $\Sigma_M^G$ défini en (\ref{eq:SGM}) est $(H,M_H)$-adapté. 
\end{proposition}

\begin{preuve} Soit $\sigma \in \Sigma_M^G$. Alors $\sigma= a_P^{G,+}$ pour un sous-groupe parabolique $P\in \fc^G(M)$. Dualement, on a un sous-groupe parabolique $\Pc\in \fc^{\Gc}(\Mc)$. Soit $P_H \in \fc^H(M_H)$ défini dualement par $\Pc_H=\Pc \cap \Hc$. Soit $\chi \in  a_P^{G,+}$. Écrivons $\chi=\chi^H +  \chi_H$ suivant $a_T^\ast=a_T^{H,\ast}+a_H^\ast$. Toute racine $\al \in \Phi^{P_H}$ vérifie $\chi(\al^\vee)=\chi^H(\al^\vee)$ et donc  $\chi^H\in a_{P_H}^{H,+}$. Il s'ensuit que $ a_P^{G,+}\subset  a_{P_H}^{H,+} + a_H^\ast$ d'où le lemme.
\end{preuve}

Soit $L\in \lc^G$ et $R=M\cap L \in \lc^L$. La proposition précédente implique  que l'éventail $\Sigma_M^G$ est $(L,R)$-adapté. On peut donc considérer la grassmannienne affine tronquée $\Xgo^L(D)$ mais on peut aussi regarder l'intersection de $\Xgo^G(D)$ avec $\Xgo^L\subset \Xgo^G$. Cela revient au même.

\begin{proposition}  Pour tout $x\in L(F)$, on a 
$$D_R^L(x)=D_M^G(x).$$
Pour tout $D\in \Div_{\Tc}(Y_\Sigma)$,
$$\Xgo^L(D)=\Xgo^G(D)\cap \Xgo^L.$$
\end{proposition}

\begin{preuve} Soit $\sigma\in \Sigma(1)$ et $P\in \fc^G(M)$ tels que $\sigma\subset a_P^{G,+}+a_G^\ast$. Il résulte de la preuve de la proposition \ref{eventail-adapte} que 
$$ a_P^{G,+}\subset  a_{Q}^{L,+} + a_L^\ast$$
où $Q\in \fc^L(R)$ est défini par $Q=P\cap L$. Il s'agit donc simplement de vérifier que pour tout $x\in L(F)$
\begin{equation}
  \label{eq:varpi}
  \varpi_\sigma(H_P(x))=\varpi_\sigma(H_Q(x)).
\end{equation}
Or un tel $x$ s'écrit $nmk$ suivant la décomposition d'Iwasawa $L(F)=N_{Q}(F) M_{Q}(F) (K\cap L)$. Comme $N_Q\subset N_P$, on a les égalités $H_P(x)=H_{M_{P}}(m)$ et $H_Q(x)=H_{M_Q}(m)$. Le caractère $\varpi_\sigma$ est naturellement un caractère de $M_P$, donc de $M_Q$, puisque $\varpi_\sigma \in a_{M_P}^\ast$. Il résulte alors de la définition des fonctions $H$ que $  \varpi_\sigma(H_{M_{P}}(m))=\varpi_\sigma(H_{M_Q}(m))$, d'où (\ref{eq:varpi}). Cela démontre la première égalité de la proposition. Le seconde en est une conséquence immédiate.
\end{preuve}
\end{paragr}

\section{Les orbites de $T$ dans la grassmannienne affine}\label{sec:orbitesga}
\label{sec:orbitesdeT}

\begin{paragr}
Soit $G$ un groupe réductif connexe défini sur $k$ et $T$ un sous-tore maximal. Le tore $T$ agit par translation à gauche sur la grassmannienne affine $\Xgo^G$. L'appartement associé à $T$ dans l'immeuble de Bruhat-Tits de $G(F)$ est un espace affine sous l'action de $a_T$. On choisit comme origine le sommet hyperspécial dont $K$ est le fixateur dans $G(F)$ et on identifie dans la suite cet appartement et $a_T$. Soit $x_0$ un point général dans une alcôve de sommet $0$. Le fixateur de $x_0$ dans $G(F)$ est alors un sous-groupe d'Iwahori de $G(F)$ noté $I$. Pour tout $\lambda\in X_\ast(T)$, soit $C_\lambda$ la cellule de Bruhat définie par
\begin{equation}
  \label{eq:Bruhat-cell}
  C_\lambda= I  \eps^\lambda K.
\end{equation}
La décomposition de Bruhat donne la partition suivante de la grassmannienne affine $\Xgo^G$ 
$$\Xgo^G=\bigcup_{\lambda\in X_\ast(T)} I \eps^\lambda K.$$
Notons que les cellules de Bruhat sont stables par $T$.

Les sous-groupes de Moy-Prasad $G_{x_0,r}$ (cf. \cite{MP}), pour $r\in \RR_+$, forment une filtration décroissante de $I$. Soit $(r_i)_{i\in \NN}$ la suite croissante des indices $r$ pour lesquels le terme $G_{x_0,r}/ G_{x_0,r+}$ du gradué associé est non trivial. Pour tout $i\in \NN$, soit $I_i=G_{x_0,r_i}$. En particulier $I_0=I$. Pour tout $\la\in X_\ast(T)$, soit  $K_\lambda = \eps^\lambda K  \eps^{-\lambda}$. Pour tout entier $i$, le tore $T$  agit par conjugaison sur les groupes $I_i$ et $I_i\cap  K_\lambda$ ; il agit donc sur le quotient $I_i/I_i\cap  K_\lambda$. L'isomorphisme évident $I_{1}/I_1\cap  K_\lambda \to C_\lambda$ est $T$-équivariant.

\begin{lemme} Soit $i\in \NN^\ast$ et $x\in I_i/I_i\cap  K_\lambda$. Soit $A$ un sous-tore du stabilisateur de $x$ dans $T$  et $M=Z_G(A)^0$. Alors $x$ appartient à l'image de 
$$(M(F)\cap I_i)I_{i+1}(I_i\cap  K_\lambda).$$
\end{lemme}

\begin{preuve} Soit $\ggo(i)$ le $k$-espace vectoriel défini par
$$\ggo(i)=\oplus \, \ggo_\alpha \eps^m$$
où $(\alpha,m)$ parcourt les couples de $X^\ast(T)\times \ZZ$ qui vérifient $\alpha(x_0)+m=r_i$ et où $\ggo_\alpha$ est l'espace propre de $T$ dans $\ggo$ associé à $\alpha$. On a un isomorphisme $T$-équivariant 
\begin{equation}
  \label{eq:composee}
 I_{i}/I_{i+1} \to \ggo(i).
\end{equation}
L'image de $(I_i\cap  K_\lambda) I_{i+1}$ par cette application est le sous-$k$-espace vectoriel
$$\ggo(i,\lambda)=\oplus \, \ggo_\alpha \eps^m,$$
où $(\alpha,m)$ parcourt le même ensemble que ci-dessus avec en sus la condition $\alpha(\lambda)+m\geq 0$.
Relevons $x$ en un élément de $I_i/I_{i+1}$,  par abus encore noté $x$. Soit $v\in  \ggo(i)$ l'image de $x$ par (\ref{eq:composee}). Soit $t\in T$. Par définition $t$ stabilise $x$ si $txt^{-1}\in (I_i\cap  K_\lambda) I_{i+1}$ ce qui se traduit par la condition $\Ad(t)(v)-v \in \ggo(i,\lambda)$. 

Soit $v=\sum v_{\al,m}$ la décomposition de $v$ suivant $\ggo(i)=\oplus \ggo_{\al}\eps^m$. Soit  $(\alpha,m)$ tel que $ v_{\alpha,m}$ soit non nul et n'appartienne pas à $\ggo(i,\lambda)$. Pour tout $t\in T$, on a  
$$\Ad(t)(v_{\alpha,m})-v_{\alpha,m}=(\alpha(t)-1)v_{\alpha,m}.$$ 
Si $t$ stabilise $x$ cet élément doit être dans $\ggo(i,\lambda)$. Cela n'est possible que si $\alpha(t)$ vaut $1$. On en déduit que $\al$ doit être une racine de $T$ dans $M$ et qu'on a $v\in (\mgo\cap\ggo(i))+ \ggo(i,\lambda)$. On conclut alors aisément.
\end{preuve}

\end{paragr}

\begin{paragr} Pour tout entier $r$, on note $\Xgo_r^G$ la réunion des orbites de dimension $r$ de $T$ dans $\Xgo^G$. Le rang semi-simple d'un sous-groupe de Lévi $M\in \lc^G(T)$ est par définition le rang du tore $T/Z_M$. On note $\lc^G_r(T) \subset \lc^G(T)$ le sous-ensemble des sous-groupes de Lévi de rang semi-simple $r$.

\begin{theoreme} \label{thm:orbites}
 Pour tout $r \in \NN$, on a la réunion disjointe 
$$\Xgo^G_r=\bigcup_{M\in \lc^G_r}\Xgo^M_r.$$
En particulier, l'ensemble des points fixes de $T$ dans $\Xgo^G$ est la grassmannienne affine du tore $T$. 
Si $r$ est le rang semi-simple de $T$, on a 
$$\Xgo^G_r= \Xgo^G -  \bigcup_{M\in \lc^G, M\not=G} \Xgo^M.$$
\end{theoreme}

\begin{preuve} Soit $n$ le rang de $T$. Commençons par prouver la première égalité. Seule l'inclusion $\subset$ n'est pas évidente. Soit $\lambda\in X_\ast(T)$ et $x\in C_\lambda\cap \Xgo_r^G$. La composante neutre du stabilisateur de $x$ est un tore $A$ de rang $n-r$. La composante neutre du centralisateur de $A$ dans $G$ est un sous-groupe de Lévi $M$ de rang semi-simple inférieur ou égal à $r$. 

Nous allons prouver que $x\in \Xgo^M$. Rappelons que la cellule $C_\la$ est isomorphe à $I_1/I_1\cap  K_\lambda$. Il suffit donc de prouver que $x$ appartient à l'image de $M(F)\cap I_1$ dans $I_1/I_1\cap  K_\lambda$. On va prouver par récurrence sur l'entier $i$ que $x$ appartient à l'image de $(M(F)\cap I_1)I_i (I_1\cap K_\lambda)$ dans $I_1/I_1\cap  K_\lambda$. Cela permettra de conclure puisque pour $i$ assez grand $I_i\subset I_1\cap K_\lambda$. Le cas $i=1$ amorce la récurrence. Passons ensuite du cas $i$ au cas $i+1$. Soit $m\in M(F)\cap I_1$ et $v\in I_i$ tels que l'image de $mv$ dans $I_1/I_1\cap  K_\lambda$ soit $x$. Pour tout $t\in A$, l'élément $x$ est fixe par $t$. Comme $t$ est central dans $M$, on a  $tvt^{-1}\in I_i\cap K_\la$.  Donc $A$ est un sous-tore du stabilisateur dans $T$ de l'image de $v$ dans le quotient  $I_i/ I_i\cap K_\la$. Le lemme précédent montre que $v$ appartient à $(M(F)\cap I_i)I_{i+1} (I_i\cap K_\lambda)$. On en déduit l'assertion pour $i+1$.

L'action de $T$ sur $\Xgo^M$ se factorise par le quotient $T/Z_M$. Si le rang semi-simple de $M$ est strictement inférieur à $r$, les orbites de $T$ dans $\Xgo^M$ sont toutes de dimension inférieure strictement à $r$. On en déduit que le rang semi-simple de $M$ est $r$. D'où la première égalité. On vérifie avec ce qui précède que la réunion est disjointe. Les autres assertions sont alors évidentes. 
\end{preuve}

\end{paragr}

\section{Fibres de Springer affines tronquées.}\label{sec:fsat}
\begin{paragr} Soit $(G,M,T)$ un triplet formé de $G$ un groupe réductif connexe et défini sur $k$, de $M$  un sous-groupe de Lévi de $G$ et de $T$ un sous-tore maximal vérifiant $T\subset M\subset G$. Soit $\Sigma$ un éventail $(G,M)$-adapté, $Y=Y_\Sigma$ la variété torique associée et $D\in  \Div_{\Tc}(Y)$. 

  Soit $\gamma\in \ggo(F)$ un élément semi-simple et régulier. À la suite de Kazhdan-Lusztig, on introduit la fibre de Springer affine associée à $\gamma$ 
$$\Xgo^G_\gamma=\{x\in \Xgo^G=G(F)/G(\of) \ | \ \Ad(x^{-1}   )\gamma \in \ggo(\of)  \}.$$

\begin{definition} \label{def:fsafftr} La fibre de Springer affine associée à $\gamma$ tronquée par le diviseur $D$ est définie par
$$\Xgo^G_\gamma(D)=\Xgo^G_\gamma(\Sigma,D)=\Xgo^G(\Sigma,D)\cap \Xgo^G_\gamma$$
où $\Xgo^G(\Sigma,D)$ est la grassmannienne affine tronquée par $D$ (cf. définition \ref{def:grafftr}).
\end{definition}\\

On obtient ainsi un sous-ind-schéma fermé de la grassmannienne affine.

\end{paragr}

\begin{paragr} \label{S:GH}Soit $H$ un groupe réductif connexe défini sur $k$ muni d'un plongement de $T$ dans $H$. On suppose que le groupe dual $\Gc$ de $G$ s'identifie à un sous-groupe réductif du dual de $\Hc$ qui contient $ \Tc$. On a donc les inclusions $\Phi^{G}(T)\subset \Phi^{H}(T)$ et $\Phi^{G,\vee}(T)\subset \Phi^{H,\vee}(T)$.

Soit $\Sigma=\Sigma^{H}_{T}$ l'éventail défini en  (\ref{eq:SGM}) et $Y=Y_\Sigma$ la variété torique associée au tore $\Tc$ et à l'éventail $\Sigma$. L'éventail $\Sigma$ est $(G,T)$-adapté. 

Soit $\gamma\in \ggo(F)$ un élément semi-simple et régulier. Par abus, on note encore $G$ le groupe réductif sur $F$ obtenu par extension des scalaires. Le centralisateur de $\gamma$ est un sous-tore maximal $T_\gamma$ défini sur $F$. Soit $A_\gamma$ le plus grand sous-$F$-tore déployé de $T_\gamma$.  Quitte à conjuguer $\gamma$  par un élément de $G(k)$, on peut et on va supposer que $A_\gamma\subset T$. Le centralisateur connexe $M_0$ de $A_\gamma$ dans $H$ est un sous-groupe de Lévi $M_0$ de $H$ qui contient $T$. Soit l'éventail $\Sigma'=\Sigma^{H}_{M_0}$ et $Y'=Y_{\Sigma'}$ la variété torique  associée. \\

Soit $\la\in X_\ast(T)$. On a vu au paragraphe \ref{S:dla} comment définir une $(H,T)$-famille $(\la_{B})_{B\in \pc^{H}(T)}$ et une $(H,M_0)$-famille  $(\la_{P})_{P\in \pc^{H}(M_0)}$. Soit $D_\la=D_\la^{\Sigma}\in \Div_{\Tc_0}(Y)$ et $D_\la'=D_\la^{\Sigma'}\in \Div_{\Tc_0}(Y')$ les diviseurs qui leur sont associés par l'application (\ref{eq:appliGM}).\\
\end{paragr}

\begin{paragr} \label{S:les2prop} Les deux propositions suivantes s'inspirent très fortement du ``lemme géométrique principal'' d'Arthur qui est au c\oe ur du développement géométrique de la formule des traces locale (cf. \cite{localtrace} lemme 4.4 et \S5).

\begin{proposition} \label{prop:tronc-ell}Il existe une constante $c\geq 0$ telle que pour tout $\la\in X_\ast(T)$ tel que $d(\la)\geq c$, on ait
$$\Xgo_\gamma^G(\Sigma, D_\la)=\Xgo_\gamma^G(\Sigma',D'_{\la}).$$
\end{proposition}

Pour tout $\la\in X_\ast(T)$, soit 
$$\la=\la^H+\la_H$$
la décomposition de $\la$ suivant $a_{T}=a_{T}^H \oplus a_H$. On pose
$$K=G(\of).$$
\begin{proposition}  \label{prop:cartan} Il existe une constante $c\geq 0$ qui vérifie : pour tout $\la\in X_\ast(T)$ tel que $d(\la)\geq c$, on a 
$$ \Xgo_\gamma^G(D_\la)= \Xgo_\gamma^G(\Sigma,D_\la)= \bigcup_\mu  \big( \Xgo_\gamma^G \cap K\eps^{\mu} K / K\big)$$
où la réunion est prise sur les cocaractères $\mu\in X_\ast(T)$ tels que $\mu^H$ appartienne au polyèdre $\mathfrak{P}(D_\la)$ défini au \S \ref{S:polyedre}.
\end{proposition}

La démonstration de ces propositions occupe respectivement les paragraphes \ref{S:tronc-ell} et \ref{S:Cartan}. Elle s'inspire bien évidemment des travaux d'Arthur mais aussi de la présentation qu'en a donnée Kottwitz (cf. \cite{K-LTF} et plus particulièrement le paragraphe 22). Auparavant, nous aurons besoin de quelques notions auxiliaires.

\end{paragr}

\begin{paragr} Soit $V$ un schéma affine de type fini sur $F$ et $\ZZ_-$ l'ensemble des entiers relatifs négatifs ou nuls. Soient $v_1$ et $v_2$ deux applications
$$V(F) \to \ZZ_-.$$
On dit que $v_2$ domine $v_1$ s'il existe $a\in\NN^\ast$ et $b\in \ZZ$ tels que pour tout $x\in  V(F)$
$$av_2(x) \leq b+v_1(x).$$
On dit que $v_1$ et $v_2$ sont équivalentes si $v_1$ et $v_2$ se dominent mutuellement.

Soit $(f_1,\ldots,f_n)$ un ensemble fini de générateurs de l'algèbre des fonctions régulières sur $V$. Pour tout $x\in V(F)$, on pose
$$v(x)= \min(0,\val(f_1(x)),\ldots, \val(f_n(x))).$$
On appelle valuation sur $V$ toute application $V(F)\to \ZZ_-$ équivalente à $v$. Cette notion ne dépend pas du choix des générateurs $f_1,\ldots,f_n$.

\begin{lemme} \label{lem:TG}Soit $v_G$ une valuation sur $G(F)$. Il existe $c_1\geq0$ tel tout $x\in G(F)$ qui vérifie
$$\Ad(x^{-1})\gamma \in \ggo(\of)$$
appartient à l'ensemble
$$\{\eps^\chi y \ |\  \chi\in X_\ast(A_\gamma) \text{ et } v_G(y)\geq -c_1 \}.$$
\end{lemme}

\begin{preuve} Par abus, on ne distingue pas dans les notations un élément $x\in G(F)$ et ses images dans $T_\gamma(F) \back G(F)$ et $A_\gamma(F) \back G(F)$. Soit $v_\ggo$ une valuation sur $\ggo(F)$. Comme le morphisme $T_\gamma\back  G \to \ggo$ défini par $x\mapsto \Ad(x^{-1})\gamma$ est une immersion fermée, on définit une valuation sur $T_\gamma(F) \back G(F)$ en posant 
$$v_{T_\gamma\back  G}(x)= v_\ggo (\Ad(x^{-1})\gamma)$$
(cf. \cite{K-LTF} proposition 18.1). Comme $A_\gamma$ est le sous-tore déployé maximal de $T_\gamma$, l'application $v_{A_\gamma\back  G}$ définie par
$$v_{A_\gamma \back  G}(x)=v_{T_\gamma\back  G}(x)= v_\ggo (\Ad(x^{-1})\gamma)$$
est une valuation sur $A_\gamma(F) \back G(F)$ (cf.  \cite{K-LTF} corollaire 18.10). Comme $A_\gamma$ est déployé, 
$$v(x)=\max_{t\in A_\gamma(F)} v_G(tx)$$
est aussi une valuation sur $A_\gamma(F) \back G(F)$ (cf.  \cite{K-LTF} proposition 18.3)).  Comme la valuation $v_\ggo$ est minorée sur $\ggo(\of)$, on en déduit qu'il existe $c\geq 0$ tel que si $x$ vérifie
$$\Ad(x^{-1})\gamma \in \ggo(\of)$$
alors $x$ appartient à l'ensemble 
$$\{t y \ |\  t\in A_\gamma(F) \text{ et } v_G(y)\geq -c \}.$$
Pour conclure, il suffit de remarquer que $A_\gamma(F)=\eps^{X_\ast(A_\gamma)}A_\gamma(\of)$ et qu'il existe $c\geq 0$ tel que $t\in A_\gamma(\of)$ et $v_G(y)\geq -c$ implique $v_G(ty)\geq -c_1$.
\end{preuve}

Pour tout $P \in \fc^{H}(T)$, soit $P_G\in \fc^G(T)$ défini par
\begin{equation}
  \label{eq:defPG}
  \Phi^{P_G}=\{\al\in \Phi^{P}\ | \ \al^\vee\in \Phi^{G,\vee}\}.
\end{equation}
Pour tout $M\in \lc^{H}(T)$ et tout $x\in G(F)$, la famille $(H_{P_G}(x))_{P \in \pc^{H}(M)}$ est une $(H,M)$-famille. \\

On choisit un produit scalaire sur $a_{T}$ invariant par le groupe de Weyl $W^H(T)$. On note $\| \cdot \|$ la norme euclidienne associée sur $a_{T}$. Ce produit scalaire sur $a_{T}$ permet d'identifier $a_{T}^\ast$ à $a_{T}$. De cette manière,  $a_{T}^\ast$ est aussi muni de la norme euclidienne $\|\cdot \|$.

\begin{lemme} \label{lem:cts} Soit $v_G$ une valuation sur $G(F)$. Pour tout $c_1\geq 0$, il existe $c_2\geq 0$ et $c_3\geq 0$ tels que pour tout $y$ et tout $P\in \fc^H(T)$ l'assertion suivante soit vraie : si $v_G(y)\geq -c_1$ alors la décomposition d'Iwasawa 
$$y=nmk$$
suivant $G(F)=N_{P_G}(F) M_{P_G}(F) K$ vérifie 
\begin{enumerate}
\item $\| \nu\| \leq c_2$ pour tout $\nu \in X_\ast(T)$ tel que 
$$m\in (K\cap M_{P_G}(F)) \eps^\nu (K\cap M_{P_G}(F)) \ ;$$ 
\item Pour tout  $\chi\in X_\ast(T)$ tel que $ \al(\chi)> c_3$ pour tout $\al \in \Phi^{N_P}$, on a 
$$\eps^\chi n\eps^{-\chi}\in K.$$
\end{enumerate}

\end{lemme}

\begin{preuve} Comme $P$ est isomorphe comme schéma au produit $N \times M$, il est clair que les valuations  $v_N(n)$ et $v_M(m)$ sont minorées pour tous  $n\in N(F)$ et $m\in M(F)$ pour lesquels il existe $k\in K$ tel que $v_G(nmk)\geq -c_1$. L'assertion 2 s'en déduit immédiatement. 

Soit $X_0$ une base du $\ZZ$-module $X^\ast(T)$.  Pour tout $\nu\in X_\ast(T)$ dominant (pour un certain sous-groupe de Borel de $M$), on pose
$$v(k_1\eps^\nu k_2)=\min_{\chi\in X_0} (0,\chi(\nu))$$
pour tout $k_1,k_2\in K\cap M_{P_G}(F)$. Alors $v$ est une valuation sur $M(F)$ (\cite{K-LTF} lemme 18.11). L'assertion 1 s'en déduit.

\end{preuve}

\end{paragr}

\begin{paragr}[Démonstration de la proposition \ref{prop:tronc-ell}.] --- \label{S:tronc-ell}Soit $v_G$ une valuation sur $G(F)$ et $c_1\geq 0$ qui vérifie les conditions du lemme \ref{lem:TG}. On en déduit une constante $c_2$ qui vérifie les conditions du lemme \ref{lem:cts}.

Montrons alors que la proposition \ref{prop:tronc-ell} vaut pour la constante 
$$c=c_2\max_{\al\in \Phi^H}(\|\al\|).$$
Soit  $\la\in X_\ast(T)$ tel que $d(\la)\geq c$. Il s'agit de montrer l'inclusion
$$\Xgo_\gamma^G(\Sigma',D'_{\la})\subset \Xgo_\gamma^G(\Sigma,D_\la)$$
car l'inclusion inverse est triviale. 

Soit $x\in G(F)$ dont la classe modulo $K$ appartient à $\Xgo_\gamma^G(\Sigma',D'_{\la})$. D'après le lemme  \ref{lem:TG} ci-dessus, 
$$x=\eps^\chi y$$ 
avec 
$$\chi\in X_\ast(A_\gamma) \ \ \text{   et   } \ \ v_G(y)\geq -c_1.$$
D'après l'assertion 1 du lemme \ref{lem:cts}, pour tous $B\in \pc^H(T)$ et $\al\in \Delta_B$ on a 
$$\al(H_{B_G}(y))\leq \|\al\| c_2 \leq c.$$
Comme $d(\la)\geq c$, la $(H,T)$-famille $(\la_B-H_{B_G}(y))_{B\in \pc^{H}(T)}$ est très positive. Soit $D$ le diviseur sur $Y$ qui lui est associé par l'application (\ref{eq:appliGM}). Il résulte des lemmes \ref{lem:cvx} et \ref{lem:famille-positive} que le projeté $\chi^H$ de $\chi$ sur $a_{T}^H$  appartient à $\mathfrak{P}(D)$ si et seulement la classe $xK$ appartient à  $\Xgo_\gamma^G(\Sigma,D_{\la})$. 

Soit $P_0\in \pc^H(M_0)$ tel que $\chi$ soit $P_0$-dominant (i.e. pour tout $\al\in \Phi^{P_0}$, $\al(\chi)\geq 0$). Soit $B\in \pc^H(T)$ tel que $B\subset P_0$. Soit $P$ le sous-groupe parabolique qui contient $P_0$ et qui admet comme facteur de Lévi $M$ défini par
$$\Delta_{M\cap B}=\{\al\in \Delta_B \ | \ \al(\chi)=0\}.$$
Comme la famille $(\la_{B'}-H_{B'_G}(y))_{B'\in \pc^H(T)}$ est très positive et que $\chi$ est $B$-dominant, on sait (cf. \cite{localtrace} lemme 3.1) que $\chi^H$ appartient à $\mathfrak{P}(D)$ si et seulement si 
\begin{equation}\label{eq:ineg-P}
\varpi_\sigma(\chi) \leq \varpi_\sigma(\la_B -H_{B_G}(y))
\end{equation}
pour tout $\sigma\in \Sigma(1)$, $\sigma\subset a_B^{H,+}$.

Puisque $xK$ appartient  à $\Xgo_\gamma^G(\Sigma',D'_{\la})$, on sait que l'inégalité (\ref{eq:ineg-P}) vaut pour les rayons $\sigma \in \Sigma(1)$ tels que $\sigma\subset a_{P_0}^{H,+}$. Comme $ a_{P}^{H,+}\subset a_{P_0}^{H,+}$, il reste à prouver l'inégalité (\ref{eq:ineg-P}) pour les rayons $\sigma\in \Sigma(1)$ tels que $\sigma\subset a_B^{H,+}$ et $\sigma\not\subset a_P^{H,+}$. Or ces rayons correspondent aux racines dans $\Delta_{M\cap B}$ (cf. la discussion du \S \ref{S:GMfamille}). D'après un lemme de Langlands (cf. \cite{dis_series} corollaire 2.2), ces inégalités résultent du fait que pour tout $\al\in \Delta_{M\cap B}$ on a 
\begin{equation}\label{eq:ineg-P2}
\al(\chi) \leq \al(\la_B -H_{B_G}(y))
\end{equation}
Cette dernière inégalité est vraie puisque  la famille $(\la_{B'}-H_{B'_G}(y))_{B'\in \pc^H(T)}$ est très positive et que $\al(\chi)=0$  pour de tels $\al$.
\end{paragr}

\begin{paragr}[Démonstration de la proposition \ref{prop:cartan}.] --- \label{S:Cartan} Soit $v_G$ une valuation sur $G(F)$ et $c_1\geq 0$ qui vérifie les conditions du lemme \ref{lem:TG}. On en déduit des constantes positives $c_2$ et $c_3$ qui vérifie les conditions du lemme \ref{lem:cts}.

Pour toute base de $\Phi^H$, soit $N$ la somme des coefficients des racines positives écrites dans cette base. Cet entier ne dépend du choix de la base. On pose
$$c_4=  \max_{\al\in \Phi^H}(\|\al \|) c_2+c_3$$ 
et 
$$c=N c_4+\max_{\al\in \Phi^H}(\|\al \|) c_2.$$

Montrons alors que la proposition \ref{prop:cartan} vaut pour la constante $c$. Soit $\la\in X_\ast(T)$ tel que $d(\la)\geq c$ et $x\in G(F)$ tel que $xK\in \Xgo_\gamma^G$. Comme précédemment on écrit 
$$x=\eps^\chi y$$ 
avec 
$$\chi\in X_\ast(A_\gamma) \ \ \text{   et   } \ \ v_G(y)\geq -c_1.$$
Soit $B\in \pc^H(T)$ tel que $\chi$ soit $B$-dominant. Alors $xK\in \Xgo_\gamma^G(D_\la)$ si et seulement si les inégalités (\ref{eq:ineg-P}) sont vraies
pour tout $\sigma\in \Sigma(1)$, $\sigma\subset a_B^{H,+}$.
Soit $P\in \fc^H(T)$ qui contient $B$ et de Lévi $M=M_P$ défini par
$$\Delta_{M\cap B}=\{\al\in \Delta_B \ | \ \al(\chi)\leq c_4\}.$$
En utilisant l'assertion 1 du lemme \ref{lem:cts}, on voit que l'inégalité (\ref{eq:ineg-P2}) vaut pour tout $\al\in \Delta_{M\cap B}$. D'après le lemme de Langlands déjà cité (cf. \cite{dis_series} corollaire 2.2), l'inégalité (\ref{eq:ineg-P}) vaut pour les cônes $\sigma$ associés aux éléments de $\Delta_{M\cap B}$.  On vient de prouver que $xK\in \Xgo_\gamma^G(D_\la)$ si et seulement si
\begin{equation}\label{eq:ineg-P3}
\varpi_\sigma(\chi) \leq \varpi_\sigma(\la_B -H_{B_G}(y))
\end{equation}
vaut pour tout $\sigma\in \Sigma(1)$, $\sigma\subset a_P^{G,+}$.
On a donc le lemme :

\begin{lemme}\label{lem:inter1}
La classe $xK$ appartient à $\Xgo_\gamma^G(D_\la)$ si et seulement si
\begin{equation}\label{eq:ineg-P4}
\varpi_\sigma(H_{P_G}(x)) \leq \varpi_\sigma(\la_P).
\end{equation}
pour tout rayon $\sigma\in \Sigma(1)$, $\sigma\subset a_P^{G,+}$.
\end{lemme}

Écrivons 
\begin{equation}
  \label{eq:def1}
  y=nmk
\end{equation}
suivant la décomposition d'Iwasawa $G(F)=N_{P_G}(F)M_{P_G}(F) K$. La décomposition de Cartan de $M(F)$ donne l'écriture
\begin{equation}
  \label{eq:def2}
m=k_1 \eps^{\nu} k_2
\end{equation}
où $k_1$ et $k_2$ appartiennent à $K^M=K \cap M_{P_G}(F)$ et $\nu\in X_\ast(T)$ est $M_{P_G} \cap B_G$-dominant. Soit $\mu\in X_\ast(T)$ l'unique cocaractère $M_{P_G} \cap B_G$-dominant qui vérifie
\begin{equation}
  \label{eq:def3}
\eps^{\chi}k_1 \eps^{\nu}k_2 \in K^M \eps^{\mu} K^M.
\end{equation}
On a alors 
\begin{equation}
  \label{eq:maj-de-la-diff}
\|\chi-\mu \| \leq \|\nu  \| \leq c_2.
\end{equation}
La première inégalité est une conséquence du lemme 21.2 de \cite{K-LTF} et la seconde provient de l'assertion 1 du lemme \ref{lem:cts}.

\begin{lemme} Pour tout $\al\in \Phi^{N_P}$, on a 
\begin{equation}
  \label{eq:ineg-mu}
 \al(\mu)\geq 0. 
\end{equation}
\end{lemme}

\begin{preuve}
Une racine dans $N_P$ s'écrit 
$$\al=\sum_{\beta \in \Delta_B} n_\beta \beta,$$
où les entiers $n_\beta$ sont tous positifs et il existe $\gamma\in \Delta_B-\Delta_{M\cap B}$ tel que $ n_\gamma >0$. Puisque $\chi$ est $B$-dominant, l'inégalité 
$$\al(\chi)\geq \gamma(\chi) >c_4$$
vaut pour tout racine $\al\in \Phi^{N_P}$. En combinant avec la majoration (\ref{eq:maj-de-la-diff}), on obtient
$$\al(\mu)=\al(\chi)+\al(\chi-\mu)\geq c_4 - \| \al \| \| \nu \| \geq c_4 -c_2 \max_{\al\in \Phi^H}(\|\al \|)\geq 0 .$$
\end{preuve}

Soit $w\in W^M(T)$ qui envoie $ \Delta_{M\cap B}$ dans l'ensemble $\{\al \in \Phi^H\  | \ \al(\mu)\geq 0\}$. Soit $B'$ le sous-groupe de Borel $wB$. Il résulte de l'inégalité (\ref{eq:ineg-mu}) que $\mu$ est $B'$-dominant. Donc $\mu^H$ appartient à $\mathfrak{P}(D_\la)$ si et seulement pour tout $\sigma\in \Sigma(1)$ tel que $\sigma\subset a_{B'}^{H,+}$ on a
\begin{equation}
  \label{eq:ineg-mu2}
  \varpi_\sigma(\mu) \leq \varpi_\sigma( \la_{B'}).
\end{equation}
Toute racine $\al\in\Phi^M$ vérifie la majoration grossière
$$\al(\chi) \leq Nc_4.$$
En utilisant cette majoration et l'inégalité (\ref{eq:maj-de-la-diff}), on obtient pour tout $\al \in w\Delta_{M\cap B}$ la majoration
$$\al(\mu) \leq \al(\chi) +\al(\mu-\chi)\leq N c_4 + \max_{\al\in \Phi^H}(\|\al \|) c_2 \leq c.$$
Or $\la$ vérifie $d(\la)\geq c$. On a  donc
$$\al(\mu) \leq \al(\la_{B'}).$$
D'après le lemme de Langlands déjà cité, on a donc l'inégalité (\ref{eq:ineg-mu2}) pour les cônes $\sigma\in \Sigma(1)$ tels que  $\sigma\subset a_{B'}^{H,+}$ et  $\sigma\not\subset a_P^{H,+}$. 

\begin{lemme}\label{lem:inter2} Soit $\mu_1\in X_\ast(T)$ tel que  
$$xK \in \Xgo_\gamma^G \cap K\eps^{\mu_1} K / K.$$
Le projeté $\mu_1^H$ appartient à $\mathfrak{P}(D_\la)$ si et seulement si 
\begin{equation}
  \label{eq:ineg-mu3}
  \varpi_\sigma(\mu) \leq \varpi_\sigma( \la_{P}).
\end{equation}
pour tout $\sigma\in \Sigma(1)$ tel que  $\sigma\subset a_P^{H,+}$.
\end{lemme}

\begin{preuve} D'après l'assertion 2 du lemme \ref{lem:cts} et les égalités  (\ref{eq:def1}), (\ref{eq:def2}) et  (\ref{eq:def3}), on a 
$$x=\eps^\chi y= \eps^\chi nmk =\eps^\chi n \eps^{-\chi} \eps^{\chi}k_1\eps^\nu k_2k\in K\eps^\mu K.$$
Donc si $x\in K\eps^{\mu_1}K$ avec $\mu_1\in X_\ast(T)$ alors $\mu_1$ est conjugué à $\mu$ par un élément de $W^G(T)$. Par conséquent, $\mu^H_1\in \mathfrak{P}(D_\la)$ si et seulement si  $\mu^H\in \mathfrak{P}(D_\la)$. Or on vient de voir que $\mu^H$ appartient à  $\mathfrak{P}(D_\la)$ si et seulement si l'inégalité  (\ref{eq:ineg-mu2}) vaut pour les rayons  $\sigma\in \Sigma(1)$ tels que $\sigma\subset a_P^{H,+}$. Pour de tels rayons, l'inégalité (\ref{eq:ineg-mu2}) se récrit comme en (\ref{eq:ineg-mu3}).
\end{preuve}

On peut maintenant donner la preuve de la proposition \ref{prop:cartan}. D'après  (\ref{eq:def1}), (\ref{eq:def2}) et  (\ref{eq:def3}), on a 
$$H_{P_G}(x)=H_{P_G}(\eps^\chi m)= H_{P_G}(\eps^\mu).$$
Donc, pour tout rayon  $\sigma\in \Sigma(1)$ tel que  $\sigma\subset a_P^{H,+}$ on a
$$ \varpi_\sigma(\mu)= \varpi_\sigma(H_{P_G}(x)).$$
La proposition \ref{prop:cartan} résulte alors des lemmes \ref{lem:inter1} et \ref{lem:inter2}.
\end{paragr}

\begin{paragr} \label{S:3assertions} Pour tout entier $N$, soit $\NN[\frac{1}{N}]\subset \QQ$ le sous-monoïde additif engendré par $\frac{1}{N}$. Soit $M\subset L$ deux sous-groupes de Lévi dans $\lc^H(T)$. Soit $B$ un sous-groupe de Borel de $H$ qui contient $T$ et $P$ un sous-groupe parabolique qui contient $B$ tels que $M_P=M$. Soit $p_M^L$ la   projection sur $a_{M}^L$ relativement à la décomposition 
$$a_{T}=a_{T}^{M}\oplus a_{M}^{L}\oplus a_L.$$
Soit $\Delta_B^\vee$ l'ensemble des coracines simples de $T$ dans $B$ et $\Delta_P^\vee=\Delta_B^\vee \cap \Phi^{N_P,\vee}$. Soit $A_M$ le tore déployé central maximal de $M$. Le groupe des cocaractères $X_\ast(A_M\cap L_{\der})$ et le projeté $p_M^L(\Delta_P^\vee)$ forment deux parties génératrices de $a_M^L$. Il existe donc des entiers $N_0$, $N$ et $N'$ et une constante $c_0>0$ tels que les assertions suivantes soient vraies.

\begin{enumerate}
\item on a les inclusions
\begin{equation}
  \label{eq:defN}
  p_{M}^L(X_\ast(T))\subset  \ZZ[\frac{1}{N'}] X_\ast(A_{M}\cap L_{\der})\subset  \ZZ[\frac{1}{N}] p_M^L(\Delta_P^\vee) \ ;
\end{equation}
 
\item la boule dans $a_T$ de centre 0 et de rayon $c_0$ contient une base du $\ZZ$-module  
$$X_\ast(A_{M_Q}\cap L_{\der})\ ;$$

\item pour tout $x\in a_{T}$ tel que $\|x\|\leq c_0 \dim(a_T)$, la décomposition 
$$p_M^L(x)=\sum_{\al\in \Delta_P^\vee} m_\al \, p_M^L(\al^\vee)$$ 
vérifie
\begin{equation}
  \label{eq:nN0}
  m_\al\leq N_0.
\end{equation}

\end{enumerate}

Soit $N_0$, $N$, $N'$ et $c_0>0$ qui vérifient les assertions 1 à 3 pour tous les sous-groupes $M$, $L$, $B$ et $P$ comme ci-dessus. Soit $B_0\in \pc^H(T)$ un sous-groupe de Borel et 
\begin{equation}
  \label{eq:defla0}
  \la_0=N_0 \sum_{\al\in \Phi^{N_{B_0}}} \al^\vee.
\end{equation}

\end{paragr}

\begin{paragr}
Soit 
$$n\ :\ \Phi^H \to \NN[\frac{1}{N}]$$
une application. Soit $M\in \lc^H(T)$. Pour toute famille $(H,M)$-orthogonale  $\mu$, on définit une famille $\mu_n=(\mu_{n,P})_{P\in \pc^H(M)}$ par
$$\mu_{n,P}=\mu_P + p_M^H(\sum_{\al\in \Phi^{N_P}} n_\al \al^\vee).$$
La famille $\mu_n$ est une famille $(H,M)$-orthogonale. En outre, on a l'inégalité
\begin{equation}
  \label{eq:ineg-delta}
  \delta(\mu_n)\geq \delta(\mu).
\end{equation}
\end{paragr}

\begin{paragr}[Pureté des fibres de Springer affines tronquées dans le cas équivalué.] ---  Soit $(G,T)$ formé d'un groupe réductif connexe $G$ et d'un sous-tore maximal $T$ tous deux définis sur $k$. Soit $\gamma\in \ggo(F)$ un élément semi-simple régulier. Soit $T_\gamma$ le centralisateur de $\gamma$ dans $G$. Soit $F_1$ une extension finie de $F$ qui déploie le tore $T_\gamma$. Il existe une unique valuation sur $F_1$ qui induit la valuation $\val$ sur $F$. On la note encore $\val$.  À la suite de Goresky-Kottwitz-MacPherson \cite{GKM-purete}, on introduit la définition suivante.

\begin{definition}\label{def:equivalue}
On dit que $\gamma$ est \emph{équivalué} s'il existe un rationnel $s\in \QQ$ tel que pour tout $\chi\in X^\ast(T_\gamma)$ on ait $\val(\chi(\gamma))\geq s$ et si $\val(\al(\gamma))=s$ pour toute racine $\al$ de $T_\gamma$ dans $G$.
\end{definition}

On reprend les hypothèses et les notations du \S\ref{S:GH}. Soit $M\in \lc^G(T)$ le centralisateur connexe de $A_\gamma$ dans $G$.

\begin{theoreme} \label{thm:purete} Soit $\gamma\in \ggo(F)$ un élément semi-simple régulier et équivalué. Il existe $c\geq 0$ tel que pour tout $\la\in X_\ast(T)$ qui vérifie  $d(\la)\geq c$ et toute application 
$n\ :\ \Phi^H \to \NN[\frac{1}{N}]$, les composantes connexes de la fibre de Springer affine tronquée $\Xgo_\gamma^G(\Sigma',D'_{\la_n})$ sont des $k$-variétés projectives et pures. 
\end{theoreme}

\begin{remarque} Pour la preuve de ce théorème, on fait appel au théorème 1.1 de \cite{GKM-purete}. Pour que ce dernier soit valable, il faut que la caractéristique de $k$ soit assez grande (pour une borne précise, cf. \emph{ibid.}).
\end{remarque}\\

\begin{preuve} En raisonnant par récurrence, on suppose que le théorème vaut pour une constante $c_1$, tout entier $N$ qui vérifie (\ref{eq:defN}) et tout sous-groupe de Lévi propre $L$ de $H$ qui contient $M_0$. Le cas $L=M_0$ amorce la récurrence. En effet, dans cette situation, il n'y a pas de différence entre une fibre tronquée et la fibre de Springer affine non tronquée $\Xgo^{M_0}_\gamma$. D'après Kazhdan-Lusztig (cf. \cite{KL} proposition 1), on sait que $\Xgo^{M_0}_\gamma$ est une variété projective. Le résultat de pureté est dû à Goresky-Kottwitz-MacPherson (cf. \cite{GKM-purete}, corollaire 1.3).  

On raisonne ensuite par récurrence sur l'entier $|n|=N \sum_{\al\in \Phi^H} n_\al$. \\

Le cas $|n|=0$ débute la récurrence. D'après Goresky-Kottwitz-MacPherson, il existe un sous-groupe parahorique, disons $J$, tel que toute intersection d'une $J$-orbite dans la Grassmannienne affine $\Xgo^G$ avec la fibre de Springer affine $\Xgo^G_\gamma$ soit pure (cf. \cite{GKM-purete}, théorème 1.1). Ce sous-groupe parahorique $J$ est le fixateur d'un certain point de ``l'appartement'' $a_{T}$. Il résulte de la construction de ce point (cf. \cite{GKM-purete}, \S 5.3, 5.4) qu'on peut le supposer dans $a_{T\cap M_{\der}}$. Il existe donc $m \in \Norm_{M_{\der}(F)}(T)$ tel que $mJm^{-1}\subset K$. La translation à gauche par $m$ envoie $\Xgo_\gamma^G(\Sigma',D'_\la)$ sur $\Xgo_{\gamma'}^G(\Sigma',D'_\la)$ avec $\gamma'=\Ad(m)\gamma$. Soit $c_2$ une constante qui vérifie les conclusions des propositions  \ref{prop:tronc-ell} et \ref{prop:cartan} pour $\gamma'$. Alors pour tout $\la\in X_\ast(T)$ tel que $d(\la)\geq c_2$, on a 
$$\Xgo^G_{\gamma'}(\Sigma',D'_\la)=\Xgo^G_{\gamma'}(\Sigma,D_\la).$$
D'après la proposition \ref{prop:cartan}, la fibre tronquée  $\Xgo^G_{\gamma'}(\Sigma',D'_\la)$ est une réunion de $K$-orbites dans la grassmannienne affine $\Xgo^G$ intersectée avec la fibre de Springer affine $\Xgo^G_{\gamma'}$. Par conséquent, la fibre tronquée $\Xgo_\gamma^G(\Sigma',D'_\la)$ est une réunion de $J$-orbites dans $\Xgo^G$ intersectées avec $\Xgo^G_{\gamma}$. Les parties de $\Xgo^G$ où la fonction $H_G$ est constante sont ouvertes et fermées. L'intersection d'une telle partie avec $\Xgo_\gamma^G(\Sigma',D'_\la)$ est une réunion finie de $J$-orbites dans $\Xgo^G$ intersectées avec $\Xgo^G_{\gamma}$. C'est donc un schéma projectif qui est pur puisque pavé par les intersections (pures) des $J$-orbites dans $\Xgo^G$  avec $\Xgo^G_{\gamma}$.

Désormais, on omet le symbole $\Sigma'$ dans les notations.
Rappelons qu'on a défini en (\ref{eq:defla0}) un cocaractère $\la_0$. On pose
\begin{equation}
  \label{eq:defc3}
c_3=\|\la_0\| \max_{\al\in \Phi^H} \|\al\|.
\end{equation}

Soit $c_4 \geq 0$ tel que tout $x\in G(F)$ qui vérifie $xK\in \Xgo_\gamma^G$ s'écrive
\begin{equation}
  \label{eq:decomp=chiy}
  x=\eps^\chi y
\end{equation}
avec $\chi\in X_\ast(T)$ et $y\in G(F)$ tel que $v_G(y)\geq -c_4$. Soit $c_5$ qui vérifie l'assertion 1 du lemme \ref{lem:cts} pour tout $y$ tel que  $v_G(y)\geq -c_4$. Posons
\begin{equation}
  \label{eq:defc}
  c=\max(c_1,c_2,c_5\max_{\al\in \Phi^H} \|\al\|)+c_3.
\end{equation}
Montrons par récurrence que le théorème vaut pour la constante $c$ et toute application 
$$n\ : \ \Phi^H \to \NN[\frac{1}{N}].$$
On a déjà vu le cas $|n|=0$. Soit $n$ tel que le théorème vaut pour $c$ et $n$.
Soit $\al\in \Phi^H$ et $n'$ l'application de $\Phi^H$ dans $\NN[\frac{1}{N}]$ définie par
$$n'_\beta=\left\lbrace\begin{array}{l} n_\beta \text{ si } \beta\not=\al \ ; \\ n_\beta+\frac{1}{N} \text{ si }  \beta=\al. \end{array} \right.$$
On va montrer que le théorème vaut alors pour $c$ et $n'$. Si $\al\in \Phi^{M_0}$, on a $\la_n=\la_{n'}$ et le résultat est acquis. On suppose désormais $\al\notin \Phi^{M_0}$. 
Soit $\la\in X_\ast(T)$ tel que $d(\la)\geq c$. Soit $x\in \Xgo_\gamma^G$, $\chi$ et $y$ donnés par (\ref{eq:decomp=chiy}). La famille $\mu$ définie pour tout $P\in \pc^H(M_0)$ par $\mu_P=\la_P-H_{P_G}(y)$ est une $(H,M_0)$-famille positive --- où le sous-groupe parabolique $P_G$ est défini en (\ref{eq:defPG})---. Il en est donc de même de la famille $\mu_{n}$. Il résulte des lemmes \ref{lem:cvx} et \ref{lem:famille-positive} et que $x$ appartient à la différence
$$\Xgo_\gamma^G(D'_{\la_{n'}})-\Xgo_\gamma^G(D'_{\la_{n}})$$
si et seulement si le projeté $p_{M_0}^G(\chi)$ appartient à la différence
$$\mathfrak{P}(D'_{\mu_{n'}})- \mathfrak{P}(D'_{\mu_{n}}).$$
Si cette dernière condition vaut alors il existe $Q\in \fc^H(M_0)$ un sous-groupe parabolique maximal propre tel que $\al\in \Phi^{N_Q}$  et
$$\varpi_Q(\la_{n,Q}-H_{Q_G}(y)) < \varpi_Q(\chi) \leq \varpi_Q(\la_{n',Q}-H_{Q_G}(y))$$
soit encore
\begin{equation}
  \label{eq:ineg-piQ}
  0< \varpi_Q(\chi-\la_{n,Q}+H_{Q_G}(y))\leq \frac{1}{N}\varpi_Q(\al^\vee).
\end{equation}
D'après l'inclusion (\ref{eq:defN}), il existe un entier $i$ tel que 
$$\varpi_Q(\chi-\la_{n,Q}+H_{Q_G}(y))= \frac{i}{N}\varpi_Q(\al^\vee)$$
et l'inégalité (\ref{eq:ineg-piQ}) est équivalente à l'égalité 
$$ \varpi_Q(\chi) = \varpi_Q(\la_{n',Q}-H_{Q_G}(y)).$$
Chaque sous-groupe parabolique $Q\in \fc^H(M_0)$ définit une face $\mathfrak{f}_Q$ du polygone $\mathfrak{P}(D'_{\mu_{n'}})$ et on obtient ainsi une bijection croissante entre $ \fc^H(M_0)$ et l'ensemble des faces de $\mathfrak{P}(D'_{\mu_{n'}})$. Les calculs précédents montrent que $x$ appartient à la différence
\begin{equation}
  \label{eq:diff-X}
 \Xgo_\gamma^G(D'_{\la_{n'}})-\Xgo_\gamma^G(D'_{\la_{n}}) 
\end{equation}
si et seulement si le projeté $p_{M_0}^G(\chi)$ appartient à une face $\mathfrak{f}_Q$ telle que $\al\in \Phi^{N_Q}$. On peut donc écrire
$$\Xgo_\gamma^G(D'_{\la_{n'}})-\Xgo_\gamma^G(D'_{\la_{n}})= \bigcup_{\{Q\in \fc^H(M_0) \ | \ \al\in \Phi^{N_Q}\}} \Xgo_\gamma^G(\mathfrak{f}_Q)$$
où l'on posé
$$ \Xgo^G_\gamma(\mathfrak{f}_Q)=\{x\in \Xgo_\gamma^G \ | \ p_{M_0}^H(\chi)\in \mathfrak{f}_Q \}.$$

 Les fibres de l'application $H_G$ définissent des parties ouvertes et fermées de la grassmannienne affine $\Xgo^G$. L'intersection d'une telle fibre $H_G^{-1}(X_0)$ avec $\Xgo_\gamma^G(D'_{\la_{n'}})$ définit un schéma projectif. Il s'agit de voir que ce schéma est pur. Comme on raisonne par récurrence, le résultat est connu pour $\Xgo_\gamma^G(D'_{\la_{n}}) \cap H_G^{-1}(X_0)$ et il suffit d'analyser la différence 
\begin{equation}
  \label{eq:diff-X2}\bigcup_{\{Q\in \fc^H(M_0) \ | \ \al\in \Phi^{N_Q}\}} \Xgo^G_\gamma(\mathfrak{f}_Q) \cap H_G^{-1}(X_0).
\end{equation}
On vérifie que $\Xgo^G_\gamma(\mathfrak{f}_Q)$ est un ouvert de la différence (\ref{eq:diff-X}). De plus, une intersection quelconque de tels ouverts est soit vide soit de la forme $\Xgo^G_\gamma(\mathfrak{f}_Q)$. Il en résulte que $\Xgo^G_\gamma(\mathfrak{f}_Q)\cap H_G^{-1}(X_0)$ est un ouvert de (\ref{eq:diff-X2}) et qu'une intersection de tels ouverts  est soit vide soit de la forme $\Xgo^G_\gamma(\mathfrak{f}_Q)  \cap H_G^{-1}(X_0)$. Par conséquent, la pureté des ouverts $\Xgo^G_\gamma(\mathfrak{f}_Q)  \cap H_G^{-1}(X_0)$ entraîne celle de (\ref{eq:diff-X2}) et par récurrence celle de $\Xgo_\gamma^G(D'_{\la_{n'}})\cap H_G^{-1}(X_0)$.

Soit $Q\in \fc^H(M_0)$ tel que $\al\in \Phi^{N_Q}$ et $X_0\in a_G$. Nous allons  prouver la pureté de $\Xgo^G_\gamma(\mathfrak{f}_Q)\cap H_G^{-1}(X_0) $. Commençons par un lemme.

\begin{lemme}Soit $X_0\in a_G$  tel que l'intersection 
$$\Xgo^G_\gamma(\mathfrak{f}_Q)\cap H_G^{-1}(X_0)$$
soit non vide. L'application $H_{Q_G}$ ne prend qu'un nombre fini de valeurs sur cette intersection.
\end{lemme}

\begin{preuve} Soit $x\in \Xgo^G_\gamma(\mathfrak{f}_Q)\cap H_G^{-1}(X_0)$. Puisque $x\in\Xgo^G_\gamma(\mathfrak{f}_Q)$ , il existe $X_1\in a_H$, $X_2\in a_T^{M_0}$ et une famille $(x_P)_{P\in \pc^H(M_0), P\subset Q}$ de réels positifs de somme égale à 1 tels que
  \begin{equation}
    \label{eq:egalite-cvxite}
    \sum_{P\in \pc^H(M_0), P\subset Q} x_P H_{P_G}(x)=\sum_{P\in \pc^H(M_0), P\subset Q} x_P \,\la_{n',P} +X_1 +X_2.
  \end{equation}
Soit $p_{M_{Q_G}}$ le projecteur sur $a_{M_{Q_G}}$ parrallèlement à $a_T^{M_{Q_G}}$. En appliquant  $p_{M_{Q_G}}$ l'égalité ci-dessus, on trouve
$$H_{Q_G}(x)= p_{M_{Q_G}}(\sum_{P\in \pc^H(M_0), P\subset Q}  x_P\,\la_{n',P}) + X_1$$
puisque $X_2\in a_T^{M_0}\subset a_T^{M}\subset a_T^{M_{Q_G}}$ et $X_1\in a_H \subset a_G \subset a_{M_{Q_G}}$. En projetant sur $a_G$, on obtient
$$X_0= p_{G}(\sum_{P\in \pc^H(M_0), P\subset Q}  x_P\,\la_{n',P}) + X_1$$
et l'on voit que $X_1$ est astreint à être dans un compact. Par conséquent, $H_{Q_G}(x)$, qui est à valeurs dans un ensemble discret, ne peut prendre qu'un nombre fini de valeurs.
\end{preuve}

 Notons $\{X_1,\ldots, X_i\}$ l'ensemble fini des valeurs de $H_{Q_G}$ sur $\Xgo^G_\gamma(\mathfrak{f}_Q)\cap H_G^{-1}(X_0)$. Quitte à changer l'indexation, on peut supposer que l'adhérence de $\Xgo^G_\gamma(\mathfrak{f}_Q)\cap H_{Q_G}^{-1}(X_j)$ dans 
$\Xgo^G_\gamma(\mathfrak{f}_Q)\cap H_G^{-1}(X_0)$ est incluse dans
$$\bigcup_{1\leq p \leq j} \Xgo^G_\gamma(\mathfrak{f}_Q)\cap H_{Q_G}^{-1}(X_p).$$
Si chaque intersection $\Xgo^G_\gamma(\mathfrak{f}_Q)\cap H_{Q_G}^{-1}(X_p)$ est pure, on voit, par récurrence sur $j$, que la réunion ci-dessus est pure. En particulier pour $j=i$, on obtient que $\Xgo^G_\gamma(\mathfrak{f}_Q)\cap H_G^{-1}(X_0)$ est pure.

On est donc ramené à prouver la pureté de toute fibre de l'application $H_{Q_G}$ intersectée avec $\Xgo^G_\gamma(\mathfrak{f}_Q)$.

Introduisons l'éventail $\Sigma^{M_Q}_{M_0}$ et la variété torique $Y^{M_Q}_{M_0}$ associée. L'application $P\mapsto P\cap M$ définit une bijection de l'ensemble des sous-groupes paraboliques $P\in \pc^H(M_0)$ tels que $P\subset Q$ sur $\pc^{M_Q}(M_0)$. On définit donc une $(M_Q,M_0)$-famille $\nu$ en posant 
$$\nu_{P\cap M_Q}=p_{M_0}^{M_Q}(\la_{n',P}).$$
Soit $D_\nu$ le diviseur sur $Y^{M_Q}_{M_0}$ associé à $\nu$. 

\begin{lemme} L'hypothèse de  récurrence entraîne que les composantes connexes de la fibre de Springer affine tronquée $\Xgo_\gamma^{M_{Q_G}}(\Sigma^{M_Q}_{M_0},D_\nu)$ sont des variétés projectives et pures. 
\end{lemme}

\begin{preuve} D'après les assertions 1 et 2 du paragraphe \ref{S:3assertions}, on peut écrire
$$p_{M_0}^{M_Q}(\sum_{\beta \in \Phi^{N_Q}} n_\al \beta^\vee)= \chi_1 +\chi_2$$
avec $\chi_1\in \ZZ[\frac{1}{N'}]X_\ast(A_{M_Q}\cap L_{\der})$ et  $\chi_2\in X_\ast(A_{M_Q}\cap L_{\der})$ tels que $\|\chi_1\|\leq c_0 \dim(a_{T})$.

Soit $P\in \pc^H(M_0)$ tel que $P\subset Q$ et $P'=P\cap M_Q$. On a 
\begin{eqnarray}\label{eq:nuP'}
\nu_{P'}&=&p_{M_0}^{M_Q}( \la_P + \sum_{\beta\in \Phi^{N_{P}}} n'_\beta \beta^\vee)\nonumber \\
&=& p_{M_0}^{M_Q}( \la_P + \sum_{\beta\in \Phi^{N_{P'}}} n'_\beta \beta^\vee) +\chi_1+\chi_2
\end{eqnarray}
En utilisant la translation à gauche par $\eps^{\chi_2}$ sur $\Xgo_\gamma^{M_{Q_G}}$, on voit qu'il suffit de considérer le cas où $\chi_2=0$, cas que nous considérons désormais. 

Quitte à changer $\la$ dans son orbite sous le groupe de Weyl $W^H(T)$, on peut et on va supposer que $\la \in c_{B_0}^+ +a_G$ (cf. \S\ref{S:dla}), où $B_0$ est le sous-groupe de Borel qui intervient dans la définition (\ref{eq:defla0}) de $\la_0$. En utilisant la définition (\ref{eq:defc}) de $c$ et l'inégalité $d(\la) \geq c$, on voit 
$$\la-\la_0 \in c_{B_0}^+ +a_G .$$
On en déduit plus généralement que pour tout sous-groupe de Borel $B\in \pc^H(T)$ on a
$$\la_B=(\la-\la_0)_B + N_0 \sum_{\beta\in \Phi^{N_B}}  \beta^\vee.$$  
Pour tout $\la\in X_\ast(T)$, on note $\la^{M_Q}=(\la^{M_Q}_R)_{R\in \pc^{M_Q}(M_0)}$ la famille $(M_Q,M_0)$-orthogonale définie comme au \S \ref{S:dla}. Soit $B\subset P$. On observe 
$$p_{M_0}^{M_Q}(\sum_{\beta\in \Phi^{N_B}}  \beta^\vee)= \sum_{\beta\in \Phi^{N_{P'}}}  \beta^\vee$$
En reprenant la ligne (\ref{eq:nuP'}), on peut écrire
$$\nu_{P'}=(\la-\la_0)_{P'}^{M_Q} + p_{M_0}^{M_Q}(\sum_{\beta\in \Phi^{N_{P'}}}(N_0+n_\beta')  \beta^\vee) +\chi_1.$$

Soit $\Delta_{P'}=\Delta_B\cap \Phi^{N_{P'}}$. D'après (\ref{eq:defN}) , on a une écriture unique
$$\chi_1=p_{M_0}^{M_Q}(\sum_{\beta \in \Delta_{P'}} m_\beta \beta^\vee)$$
avec $m_\beta\in \ZZ[\frac{1}{N}]$. On prolonge la fonction $m$ à tout $\Phi^{M_Q}-\Delta_{P'}$ par
$$m_\beta= \left\lbrace\begin{array}{l} -m_\beta \text{ si } -\beta \in \Delta_{P'} \\ 0  \text{ sinon.}
  \end{array}\right.
  $$

On définit alors une application $n''\ :\ \Phi^{M_Q} \to  \ZZ[\frac{1}{N}]$
par
$$n''_\beta= \left\lbrace\begin{array}{l} 0 \text{ si } \beta \in \Phi^{M_0}\\ n'_\beta + N_0 +m_\beta  \text{ si } \beta \in \Phi^{M_Q}-\Phi^{M_0} \end{array}\right.
  $$
L'inégalité  (\ref{eq:nN0}) assure que cette application est en fait à valeurs dans $\NN[\frac{1}{N}]$. Il s'ensuit que 
$$\nu_{P'}=(\la-\la_0)_{n'',P'}.$$
Puisque $d(\la)\geq c$ , on a  $d(\la-\la_0)\geq c_1$ (cf. (\ref{eq:defc3}) et (\ref{eq:defc})), on peut utiliser l'hypothèse de récurrence qui donne la pureté voulue.
\end{preuve}

Soit $X\in a_{M_{Q_G}}$. Le lemme précédent entraîne que $\Xgo_\gamma^{M_{Q_G}}(D_\nu)\cap H_{M_{Q_G}}^{-1}(X)$ est une variété projective et pure. La pureté de l'intersection $\Xgo^G_\gamma(\mathfrak{f}_Q)\cap H_{Q_G}^{-1}(X)$ est alors une conséquence du lemme suivant.

\begin{lemme}  Le morphisme
$$\pi_{Q_G} \ : \ \Xgo^G_\gamma(\mathfrak{f}_Q)\cap H_{Q_G}^{-1}(X) \to \Xgo_\gamma^{M_{Q_G}}(D_\nu)\cap H_{M_{Q_G}}^{-1}(X) $$
fait de $\Xgo^G_\gamma(\mathfrak{f}_Q)\cap H_G^{-1}(X)$ un fibré vectoriel (itéré) sur $\Xgo_\gamma^{M_{Q_G}}(D_\nu)\cap H_{M_{Q_G}}^{-1}(X)$.
\end{lemme}

\begin{preuve}
En appliquant le projecteur $p_{M_0}^{M_Q}$ à l'égalité (\ref{eq:egalite-cvxite}), on voit que la rétraction $\pi_{Q_G}$ envoie $\Xgo^G_\gamma(\mathfrak{f}_Q)$ dans la grassmannienne affine tronquée $\Xgo^{M_{Q_G}}(D_\nu)$. Soit $xK\in  \Xgo^G$ et $x=mnk$ suivant la décomposition d'Iwasawa 
$$G(F)=M_{Q_G}(F) N_{Q_G}(F) K.$$
Alors $\Ad(x^{-1})\gamma\in \ggo(\of)$ si et seulement si $\Ad(m^{-1})\gamma \in \mgo_{Q_G}(\of)$ et
\begin{equation}
  \label{eq:nmgamma}
  \Ad(n^{-1})((\Ad m^{-1})\gamma)- \Ad m^{-1}\gamma \in \ngo_{Q_G}(\of).
\end{equation}
En particulier, $\pi_{Q_G}$ envoie $\Xgo^G_\gamma(\mathfrak{f}_Q)$ dans $\Xgo^{M_{Q_G}}_\gamma(D_\nu)$. Quitte à translater à gauche par un élément de $M(F)$, on peut et on va supposer que $X=0$. Soit $N$ le radical unipotent de $Q_G$, $M'$ le groupe dérivé de $M_{Q_G}$ et $\tilde{M}'$ l'ensemble des $m\in M'(F)$ tels que $mM'(\of) \in \Xgo_\gamma^{M_{Q_G}}(D_\nu)$. Soit $\tilde{M}'\times_{M'(\of)} N(F)/N(\of)$ le quotient de $\tilde{M}'\times N(F)/N(\of)$ par $M'(\of)$ où le groupe $M'(\of)$ agit sur le premier facteur par translation à droite et par conjugaison sur le second. 

Notons $V$ la partie de $\tilde{M}'\times_{M'(\of)} N(F)/N(\of)$ formée des couples $(m,n)$ tels que la condition (\ref{eq:nmgamma}) soit satisfaite. Alors $\pi_{Q_G}^{-1}(\Xgo_\gamma^{M_{Q_G}}(D_\nu)\cap H_{M_{Q_G}}^{-1}(0))$ s'identifie à  $\tilde{M}'\times_{M(\of)} N(F)/N(\of)$ et $V$ s'identifie à $\Xgo^G_\gamma(\mathfrak{f}_Q)\cap H_{Q_G}^{-1}(0)$. De plus, le morphisme $\pi_{Q_G}$ s'identifie à la projection évidente 
$$\tilde{M}'\times_{M'(\of)} N(F)/N(\of) \to  \Xgo_\gamma^{M_{Q_G}}(D_\nu)\cap H_{M_{Q_G}}^{-1}(0).$$

Introduisons la suite centrale décroissante
$$N_0=N \supset N_1=[N,N] \supset \ldots \supset N_{i+1}=[N,N_i] \ldots \supset N_r=\{1\}.$$

Soit $V_i$ la partie de 
$$\tilde{M}'\times_{M(\of)} N(F)/(N_i(F) + N(\of))$$
formé des couples $(m,n)$ tels que
\begin{equation}
  \label{eq:nmgamma2}
  \Ad(n^{-1})((\Ad m^{-1})\gamma)- \Ad m^{-1}\gamma \in \ngo(\of)+\ngo_i(F).
\end{equation}
Remarquons que $V_0=\Xgo_\gamma^{M_{Q_G}}(D_\nu)\cap H_{M_{Q_G}}^{-1}(0)$ et $V_r=V$.
Pour conclure, il suffit de voir que le morphisme évident
$$\tilde{M}'\times_{M(\of)} N(F)/(N_{i+1}(F) + N(\of)) \to \tilde{M}'\times_{M(\of)} N(F)/(N_i(F) + N(\of))$$
se restreint en un morphisme
$$V_{i+1} \to V_i$$
qui fait de $V_{i+1}$ un fibré vectoriel au-dessus de $V_i$. C'est une conséquence facile des calculs de Kazhdan-Lusztig (cf. \cite{KL} \S 5).

\end{preuve}

\end{preuve}

\end{paragr}

\begin{paragr} Soit $L \in \lc^H(T)$ tel que $T_\gamma\subset L$. On considère l'éventail $\Sigma_L^H$. Soit $\la\in X_\ast(T)$. Posons 
$$D_{L,\la}=D^{\Sigma_L^H}_\la,$$
où le diviseur à droite est défini en (\ref{eq:appliGM}).  Introduisons la définition suivante qui est commode pour la suite.

\begin{definition} \label{def:purete} On dit que la fibre de Springer affine tronquée
$$\Xgo_\gamma^G(\Sigma_L^H,D_{L,\la})$$ 
est pure si elle est une réunion de variétés projectives, pures et stables par $T$.
\end{definition}

On a alors le corollaire suivant au théorème \ref{thm:purete}.

\begin{corollaire}  \label{cor:purete} Soit $\gamma\in \ggo(F)$ un élément semi-simple régulier et équivalué. Il existe $c\geq 0$ tel que pour tout $\la\in X_\ast(T)$ qui vérifie $d(\la)\geq c$ la fibre de Springer affine tronquée
$$\Xgo_\gamma^G(\Sigma_L^H,D_{L,\la})$$ 
est pure.
\end{corollaire}

\begin{preuve} Il est clair que $\Xgo_\gamma^G(\Sigma_L^H,D_{L,\la})$ est la réunion croissante des fibres tronquées $\Xgo_\gamma^G(\Sigma',D'_{\la_{n_i}})$ pour $i\in \NN$ où $n_i$ est l'application 
$$n_i\ :\ \Phi^H \to \NN[\frac{1}{N}]$$
définie par $n_\al=0$ sauf si $\al\in \Phi^L-\Phi^{M_0}$ auquel cas $n_\al=i$. Le théorème  \ref{thm:purete} permet de conclure.

\end{preuve}

\end{paragr}

\section{Homologie des fibres de Springer affines tronquées pour $\SL(2)$}
\label{sec:SL2}

\begin{paragr} \label{par:SL2.1}Dans cette section, on considère le groupe $G=\SL(2)$ défini sur $k$, $T$ le sous-tore maximal diagonal. On note $B$ le sous-groupe de Borel des matrices triangulaires supérieures et $\bar{B}$ le sous-groupe de Borel opposé. Soit $\al$ l'unique racine de $T$ dans $B$ et $\al^\vee$ la coracine associée. Le groupe $X_\ast(T)$ est engendré par $\al^\vee$. Pour tout $x\in \SL(2,F)$, on définit les entiers $h_B(x)$ et $h_{\bar{B}}(x)$ par
$$H_B(x)=h_B(x)\al^\vee \ \ \text{et} \ \ H_{\bar{B}}(x)=-h_{\bar{B}}(x)\al^\vee.$$
Soit  $\Sigma$ un éventail $(G,T)$-adapté. Par conséquent, $\Sigma$ est inclus dans l'éventail 
$$\Sigma_T^G=\{(0), a_B^{G,+}, a_{\bar{B}}^{G,+}\}.$$
Soit $D=\sum_{\sigma \in \Sigma(1)}n_\sigma D_\sigma$ un diviseur dans $\Div_{\Tc}(Y_\Sigma)$. On pose $$n_B=n_{a_B^+}$$ si le cône $a_B^+$ est dans l'éventail $\Sigma$ et $n_B=+\infty$ sinon. On définit de même $n_{\bar{B}}$. Soit $t\in F$ et 
$$\gamma=\left( \begin{array}{cc} t & 0\\ 0 & -t \end{array}\right)\in \tgo(F).$$
On considère alors la grassmannienne affine tronquée 
$$\Xgo(D)=\{x\in \Xgo_\gamma \ |\ h_B(x)\leq n_B \ \text{et}\ h_{\bar{B}}(x)\leq n_{\bar{B}}\}.$$
et la fibre de Springer affine tronquée
$$\Xgo_\gamma(D)=\{x\in \Xgo(D) \ | \ x^{-1}\gamma x \in \mathfrak{sl}(2,\oc)\}$$
associées à ces données. Lorsque $\Sigma$ ne contient que le cône nul, on obtient la fibre de Springer affine usuelle $\Xgo_\gamma$.

\begin{lemme} \label{lem:cond-val} Soit $x=\left( \begin{array}{cc} a & b\\ c & d \end{array}\right)\in \SL(2,F)$. Alors $h_B(x)=-\min(\val(c),\val(d))$ et $h_{\bar{B}}(x)=-\min(\val(a),\val(b))$.
En particulier, $xK \in \Xgo(D)$, resp. $ xK \in \Xgo_\gamma$ si et seulement si les conditions suivantes 1 et 2 (resp. 3 et 4) sont réalisées :
\begin{enumerate}
\item $\val(a)\geq -n_{\bar{B}}$ et $\val(b)\geq -n_{\bar{B}}$ ;
\item $\val(c) \geq -n_B$ et $\val(d) \geq -n_B$ ;
\item $\val(ad+bc)\geq -\val(t)$ ;
\item $\val(bd)\geq -\val(2t)$ et $\val(ac)\geq -\val(2t)$.  
\end{enumerate}
\end{lemme}

\begin{preuve} Par la décomposition d'Iwasawa, il existe $s\in F^\times$ et $u\in F$ de sorte que 
$$x=\left( \begin{array}{cc} s & u\\ 0 & s^{-1}\end{array}\right) K.$$ 
On a alors $h_B(x)=\val(s)$. Considérons l'action standard de $\SL(2,F)$ sur l'espace $F^2$.  La projection sur la deuxième coordonnée du sous-$\oc$-module  $x\oc^2$ de $F^2$ est le sous-$\oc$-module de $F$
$$s^{-1}\oc=c\oc+d\oc.$$
De là, on en déduit l'égalité $\val(s^{-1})=\min(\val(c),\val(d))$ et le calcul de $h_B(x)$. On procède de même pour $h_{\bar{B}}(x)$. Les conditions 1 et 2 sont alors évidentes. Un calcul élémentaire montre que $x\in \Xgo_\gamma$ si et seulement si 
$$\left( \begin{array}{cc} t(ad+bc)  & 2tbd \\  -2tac &  -t(ad+bc) \end{array}\right)\in \ggo(\oc),$$
d'où les conditions 3 et 4.
\end{preuve}

\end{paragr}

\begin{paragr} Soit $I$ le sous-groupe d'Iwahori standard de $G(F)$.  Les intersections des cellules de Bruhat avec la fibre de Springer affine tronquée définissent, pour tout $n\in \ZZ$, des ``cellules tronquées'' 
$$C_n=I \eps^{n\al^\vee} K \cap \Xgo_{\gamma}(D).$$

\begin{proposition} Pour $n\in \ZZ$, les cellules tronquées $C_n$ sont soit vides soit isomorphes à des espaces affines standard.
\end{proposition}

\begin{preuve} Supposons tout d'abord $n\in \NN$. Alors l'application 
  \begin{equation}
    \label{eq:app-cellule}
    u\mapsto  x_u=\left( \begin{array}{cc} 1 & u\\ 0 & 1\end{array}\right)   \left( \begin{array}{cc} \eps^n & 0\\ 0 & \eps^{-n}\end{array}\right)K 
  \end{equation}
induit un isomorphisme de $\oc / \eps^{2n}\oc$ sur la cellule de Bruhat $I\eps^{n\al^\vee}K$. Pour l'élément $x_u$, les conditions 1 à 4 du lemme se traduisent par les inégalités suivantes
\begin{itemize}
\item $n\geq -n_{\bar{B}}$ et $\val(u)\geq n-n_{\bar{B}}$ ;
\item $n \leq n_B$ ;
\item $\val(t)\geq 0$ ;
\item $\val(u)\geq 2n-\val(2t)$.
\end{itemize}
Supposons $-n_{\bar{B}}\leq n\leq n_B$ et $\val(t)\geq 0$ sans quoi $C_n$ est vide. On pose 
$$m=\max(0, n-n_{\bar{B}},2n-\val(2t)).$$
Alors l'application (\ref{eq:app-cellule}) induit un isomorphisme de $\eps^m\oc/\eps^{2n}\oc$ sur $C_n$. Si $n\in \ZZ-\NN$, on procède de même en remplaçant l'application  (\ref{eq:app-cellule}) par
 $$u\mapsto  \left( \begin{array}{cc} 1 & 0\\  u & 1\end{array}\right)   \left( \begin{array}{cc} \eps^n &0 \\  0& \eps^{-n}\end{array}\right)K $$
qui induit un isomorphisme de $\eps\oc /\eps^{-2n}\oc$ sur $I\eps^{n\al^\vee}K$.
\end{preuve}

 \begin{remarque} On voit que $\Xgo_\gamma\not=\emptyset$ si et seulement si $\gamma\in \tgo(\oc)$. On suppose dorénavant que cette dernière condition est satisfaite.
 \end{remarque}
\end{paragr}

\begin{paragr}\label{S:Xn}Pour $n\in \NN$, on pose 
$$C_{\leq n}=\bigcup_{\{n\in \ZZ,\ -n\leq m\leq n\}} C_{m}.$$
Ce sont des schémas projectifs sur $k$.

\begin{corollaire} \label{cor:Xnpur}La fibre de Springer affine tronquée $\Xgo_\gamma(D)$ est la réunion croissante des schémas projectifs $C_{\leq n}$, $n\in\NN$. L'homologie de $C_{\leq n}$  est pure et s'annule en degrés impairs.
\end{corollaire}

\begin{preuve} La première assertion est claire. Pour $n\in \NN$, la  différence $C_{\leq n+1}-C_{\leq n}$ est la réunion des deux cellules tronquées $C_{\pm (n+1)}$. On vient de voir dans le lemme précédent que ces cellules sont vides ou isomorphes à des espaces affines. Par récurrence sur $n$, on en déduit que  l'homologie de $C_{\leq n}$ est pure et s'annule en degrés impairs. 
\end{preuve}

\end{paragr}

\begin{paragr} On associe à chaque $\lambda\in k^\times$ un automorphisme $\sigma_\lambda$ du corps $F$ défini par $\sigma_\lambda(\eps^m)=\lambda^m\eps^m$ pour tout $m\in \ZZ$. On en déduit alors une action de $k^\times$ sur $F$ qui préserve la valuation donc une action de $k^\times$ sur $\SL(2,F)$ qui préserve $\SL(2,\oc)$. On obtient ainsi une action de $\Gm$ sur $\Xgo$ qui commute à l'action de $T$. Suivant \cite{GKM} \S5.5, on appelle tore pivotant le groupe $\Gm$ muni de cette action sur la grassmannienne affine $\Xgo$ et tore étendu le produit $\tilde{T}=T\times \Gm$. Le tore étendu $\tilde{T}$ préserve les fibres tronquées $\Xgo_\gamma(D)$ (puisqu'il préserve clairement les conditions 1 à 4 du lemme \ref{lem:cond-val}), les cellules tronquées $C_n$ et leurs réunions $C_{\leq n}$. Le caractère $\alpha$ de $T$ donne par composition avec la projection canonique $\tilde{T}\to T$ un caractère de $\tilde{T}$ encore noté $\alpha$. Pour $n\in \ZZ$, on note $\lambda_n$ le caractère de $\tilde{T}$ obtenu en composant  la projection $\tilde{T}\to \Gm$ avec le caractère $t\mapsto t^n$ de $\Gm$. 

\begin{lemme}\label{lem:ptsfixes-orbitesSL2} Soit $\gamma\in \tgo(\oc)$ et  $n\in \NN$. Dans  $C_{\leq n}$, il y a un nombre fini de points fixes sous l'action de $\tilde{T}$ : ce sont les points $t_l=\alpha^{\vee}(\eps^l)K$ où $l$ est un entier qui vérifie les deux conditions suivantes :
  \begin{enumerate}
  \item $-n_{\bar{B}}\leq l \leq n_B$ ;
  \item $-n\leq l \leq n.$
  \end{enumerate}
Deux tels points $t_l$ et $t_m$, avec $l\not=m$, sont dans l'adhérence d'une orbite de $\tilde{T}$ dans $C_{\leq n}$ de dimension $1$ si et seulement si 
$$\val(2t)\geq |l-m|.$$ 
Sous cette dernière condition, l'orbite est unique et un point de cette orbite a pour centralisateur le noyau du caractère  $\alpha\lambda_{l+m}$. Il n'y a donc qu'un nombre fini d'orbites de dimension $1$ de $\tilde{T}$ dans  $C_{\leq n}$.
\end{lemme}

\begin{preuve} D'après le théorème \ref{thm:orbites}, les points de $\Xgo$ fixes sous $T$ sont les points de $\Xgo^T$ à savoir l'ensemble des $t_l$, $l\in \ZZ$. Ces points sont aussi fixes sous $\tilde{T}$ et ils sont forcément dans $\Xgo_\gamma$ puisqu'on a supposé que $\gamma \in \tgo(\oc)$. On vérifie immédiatement que les conditions 1 et 2 expriment l'appartenance de $t_l$ respectivement à $\Xgo(D)$ et à $C_{\leq n}$. 

Prenons deux entiers $l$ et $m$ tels que $m<l$. D'après le lemme 6.4 de \cite{GKM}, il existe une unique $\tilde{T}$-orbite de dimension $1$ dans $\Xgo$ qui contient dans son adhérence $t_l$ et $t_m$ à savoir l'orbite de 
$$x_{l,m}=\left( \begin{array}{cc} \eps^l &  \eps^m \\ 0  & \eps^{-l}\end{array}\right) K.$$ 
Supposons que $t_l$ et $t_m$ appartiennent à $C_{\leq n}$. À l'aide du lemme \ref{lem:cond-val}, on voit que $x_{l,m}$, et donc sa $\tilde{T}$-orbite, est dans $\Xgo_\gamma(D)$ si et seulement si $\val(2t)\geq l-m$. Il reste à voir que, sous cette dernière condition,  $x_{l,m}\in C_{\leq n}$. 

Si $l+m\geq 0$, en écrivant 
$$ \left( \begin{array}{cc} \eps^l &  \eps^m \\ 0  & \eps^{-l}\end{array}\right)=\left( \begin{array}{cc} 1 &  \eps^{m+l} \\ 0  & 1 \end{array}\right)\left( \begin{array}{cc} \eps^l & 0 \\ 0  & \eps^{-l}\end{array}\right),$$
on voit que $x_{l,m}\in C_l$. Si $l+m\leq -1$, on utilise l'égalité
$$\left( \begin{array}{cc} \eps^l &  \eps^m \\ 0  & \eps^{-l}\end{array}\right)=\left( \begin{array}{cc} 1 &  0 \\ \eps^{-(m+l)}  & 1\end{array}\right)\left( \begin{array}{cc} \eps^m & 0 \\ 0  & \eps^{-m}\end{array}\right)\left( \begin{array}{cc} \eps^{l-m} & 1 \\ -1  & 0\end{array}\right)$$
pour obtenir $x_{l,m}\in C_m$. 

Finalement, le calcul du centralisateur de $x_{l,m}$ dans $\tilde{T}$ est immédiat (cf.  la preuve du lemme 6.4 de \cite{GKM}).
\end{preuve}
 \end{paragr}

\begin{paragr} \label{par:SL2-4}On note $\mathcal{S}$ et $\mathcal{D}$ les algèbres définies aux \S\ref{S:coh-equiv1} et  \S\ref{S:coh-equiv3} relatives au tore $T$. On note $\tilde{\Sc}$ et $\tilde{\dc} $ les algèbres analogues relatives à $\tilde{T}$.

Soit  $\Qlb[\Tc]=\Qlb[X_\ast(T)]$ l'algèbre des fonctions régulières sur $\Tc$. Soit $\oc(D)$ le faisceau cohérent sur la variété torique $Y_\Sigma$ associé au diviseur $\Tc$-équivariant $D$. Soit $\Gamma(D)$ le $\Qlb$-espace vectoriel des sections globales du faisceau $\oc(D)$. C'est le sous-$\Qlb$-espace de $\Qlb[\Tc]$ engendré par $(\al^\vee)^n$ pour $-n_B\leq n \leq n_{\bar{B}}$. 

L'ensemble  $\Xgo^T(D)$ des points de la grassmannienne affine tronquée du tore est l'ensemble
$$\{t_n \ | \  -n_{\bar{B}}\leq n\leq  n_B\}.$$
En identifiant $(\al^\vee)^n$ au point fixe $t_{-n}$, on définit un isomorphisme de $\dc$-module 
$$\Gamma(D)\otimes_{\Qlb}\Sc \to H^T_\bullet( \Xgo^T(D)).$$
Pour tout $d\in \NN$, $d\geq 1$, soit $R_d$ le sous-$\dc$-module de $\Gamma(D)\otimes_{\Qlb}\Sc$  engendré par 
$$(1-\al^\vee)^d (\al^\vee)^n \otimes \ker(\partial_\al^d),$$
pour $n\in \ZZ$ vérifiant $-n_{B}\leq n \leq n+d \leq  n_{\bar{B}}$.

\begin{proposition} \label{prop:hom-1orb-SL2} On a la suite exacte 
$$\begin{CD}
0 @>>> H_\bullet^T(\Xgo_\gamma(D),\Xgo^T(D)) @>>> H_\bullet^T(\Xgo^T(D))@>>> H_\bullet^T(\Xgo_\gamma(D)) @>>> 0
\end{CD}$$
et l'image du premier morphisme dans $H_\bullet^T(\Xgo^T(D))$ s'identifie au sous-$\dc$-module
$$\sum_{d=1}^{\val(2t)} R_d.$$
\end{proposition}

\begin{preuve}Soit $n\in \NN$. On a défini au \S\ref{S:Xn} un schéma projectif $C_{\leq n}$, qui est muni d'une action du tore $\tilde{T}$. Soit  $C_{\leq n}^0$ le sous-ensemble des points fixes de $\tilde{T}$ et  $C_{\leq n}^1$ la réunion des orbites de $\tilde{T}$ de dimension inférieure ou égale à $1$. Puisque l'homologie de $C_{\leq n}$ est pure, le lemme de Chang-Skjelbred (cf. \S\ref{S:coh-equiv4}) donne une suite exacte
$$ H^{T}_{\bullet}(C_{\leq n}^1,C_{\leq n}^0) \rightarrow H^{T}_{\bullet}(C_{\leq n}^0)\rightarrow H^{T}_{\bullet}(C_{\leq n})\rightarrow 0.$$
Or les orbites de $\tilde{T}$ dans $C_{\leq n}$ de dimension inférieure ou égale à $1$ sont en nombre fini (cf. lemme \ref{lem:ptsfixes-orbitesSL2}). Alors la suite exacte s'explicite de la manière suivante (cf. \S\ref{S:coh-equiv4} et lemme  \ref{lem:ptsfixes-orbitesSL2})

\begin{equation}
    \label{eq:suite-sl2-Ttilde-n-explicite}
 \bigoplus_{\begin{array}{c}(l,m)\in I^2\\ 0< l-m \leq \val(2t)\end{array}} \!\!\!\!\!\!\!\! \ker(\partial_\alpha +(m+l)\partial_{\lambda}) \longrightarrow \sum_{l\in I}  (\al^\vee)^{-l}\otimes \tilde{\Sc} \longrightarrow H^{\tilde{T}}_\bullet(C_{\leq n})\longrightarrow 0,
\end{equation}
 où $I$ est la partie de $\ZZ$ formée des entiers $l$ qui vérifient les égalités 1 et 2  du lemme \ref{lem:ptsfixes-orbitesSL2} et la  flèche de gauche envoie $f\in \ker(\partial_\alpha +(m+l)\partial_{\lambda})$ sur $((\al^\vee)^{-l}-(\al^\vee)^{-m})\otimes f$.  Soit $R_\gamma(D)$ l'image de cette flèche : c'est le sous-$\dc$-module de $\Gamma(D)\otimes\tilde{\Sc}$ engendré par 
$$((\al^\vee)^{-l}-(\al^\vee)^{-m})\otimes \ker(\partial_\alpha +(m+l)\partial_{\lambda})$$
pour des entiers $l$ et $m$ qui vérifient $ 0< l-m \leq \val(2t)$ et  $-n_{\bar{B}}\leq m <l \leq n_B$. Puisque la limite inductive sur $n$ de $H^{\tilde{T}}_\bullet(C_{\leq n})$ est l'homologie $H^{\tilde{T}}_\bullet(\Xgo_\gamma(D))$, on obtient une suite exacte courte
\begin{equation}
    \label{eq:suite-sl2-Ttilde}
0\longrightarrow R_\gamma(D) \longrightarrow \Gamma(D)\otimes\tilde{\Sc} \longrightarrow H^{\tilde{T}}_\bullet(\Xgo_\gamma(D))\longrightarrow 0.
\end{equation}
Comme l'homologie de $\Xgo_\gamma(D)$ est pure, on a la relation 
$$H^{T}_\bullet(\Xgo_\gamma(D))=H^{\tilde{T}}_\bullet(\Xgo_\gamma(D))\cap \ker \partial_\lambda.$$
On vérifie que la dérivation $\partial_\lambda$ induit une surjection de $R_\gamma(D)$ sur lui-même. On déduit alors de (\ref{eq:suite-sl2-Ttilde}) la suite exacte
\begin{equation}
    \label{eq:suite-sl2-T}
0\longrightarrow R_\gamma(D) \cap \ker \partial_\lambda \longrightarrow \Gamma(D) \otimes\Sc \longrightarrow H^{T}_\bullet(\Xgo_\gamma(D))\longrightarrow 0.
\end{equation}
Les sous-$\dc$-modules $R_\gamma(D) \cap \ker \partial_\lambda$ et $\sum_{d=1}^{\val(2t)} R_d$ sont égaux. En effet, le lemme 12.2 de \cite{GKM} donne une inclusion $R_d \subset R_\gamma(D) \cap \ker \partial_\lambda$. Pour montrer l'inclusion inverse, on procède comme dans la preuve de la proposition 12.7 de \cite{GKM}.

Finalement, la  suite exacte courte
\begin{equation*}
0\longrightarrow H^{T}_\bullet(\Xgo_\gamma(D),\Xgo^T(D))  \longrightarrow   H^{T}_\bullet(\Xgo^T(D))  \longrightarrow H^{T}_\bullet(\Xgo_\gamma(D))\longrightarrow 0,
\end{equation*}
qui provient du lemme de Chang-Skjelbred et de l'injectivité de la flèche de gauche due à l'annulation de l'homologie de $\Xgo^T(D)$ en degré impair, s'identifie à la suite (\ref{eq:suite-sl2-T}).
\end{preuve}

\end{paragr}

\section{Homologie des orbites de dimension $0$ et $1$ des fibres de Springer tronquées.}\label{sect:Hom_orb}

\begin{paragr} \label{S:moduleR} Soit $G$ un groupe réductif connexe, $T$ un sous-tore maximal et $M$ un sous-groupe de Lévi qui contient $T$, tous ces groupes étant définis sur $k$. Soit  $\Sigma$ un éventail $(G,M)$-adapté et $D$ un diviseur $\Tc$-invariant sur la variété torique $Y_\Sigma$. Soit $\Gamma(D)$ le $\Qlb$-espace vectoriel des sections globales du faisceau cohérent $\oc(D)$ : c'est le sous-$\Qlb$-espace de $\Qlb[\Tc]$ engendré par 
$$\{\la \in X_\ast(T) \ | \ (\la) + D \geq 0  \}.$$
Soit $d\in \NN$, $d\geq 1$ et $\al\in \Phi^G$. Soit $R_{\al,d}$ le sous-$\dc$-module de $\Gamma(D)\otimes_{\Qlb} \Sc$ engendré par 
$$(1-\al^\vee)^d \lambda \otimes \ker(\partial_\al^d),$$
pour $\la\in X_\ast(T)$ tel que l'ensemble
$$\{\la, \al^\vee \la,\ldots,(\al^\vee)^d  \la\}$$
soit inclus dans $\Gamma(D)$. On notera que $R_{\al,d}=R_{\al^{-1},d}$.  
\end{paragr}

\begin{paragr} Pour $\al\in \Phi^G_+$, on note  $L_\alpha$ l'unique sous-groupe $L_\al\in \lc^G(T)$ de rang semi-simple $1$ défini par $\Phi^{L_\al}=\{\pm \al\}$. On a alors une inclusion naturelle 
$$\Xgo^{L_\al}_\gamma(D) \subset \Xgo^{G}_\gamma(D)$$
de la fibre de Springer affine tronquée de $L_\al$ dans celle de $G$. Afin d'alléger les notations, on pose $\Xgo^\al_\gamma(D)=\Xgo^{L_\al}_\gamma(D)$ et $\Xgo_\gamma(D)=\Xgo^{G}_\gamma(D)$. Les points  de $\Xgo_\gamma(D)$ fixes sous $T$ s'identifient aux points de $\Xgo^T(D)$, la grassmannienne affine tronquée de $T$ (cf. théorème \ref{thm:orbites}). On note enfin $\Xgo_\gamma(D)_1$ la réunion des orbites de $T$ dans $\Xgo_\gamma(D)$ de dimension inférieure ou égale à $1$.

\begin{proposition} \label{prop:hom-orb-0-1} On a les isomorphismes de $\dc$-modules suivants :
  \begin{enumerate}
  \item $H_\bullet^T(\Xgo^T(D))\simeq \Gamma(D)\otimes_{\Qlb}\Sc$  ;
\item $H_\bullet^T(\Xgo_\gamma(D)_1,\Xgo^T(D))\simeq\displaystyle \bigoplus_{\al \in \Phi^G_+} H_\bullet^T(\Xgo^\al_\gamma(D),\Xgo^T(D))$ ;
\item $H_\bullet^T(\Xgo^\al_\gamma(D),\Xgo^T(D))\simeq \displaystyle \sum_{d=1}^{\val(\alpha(\gamma))} R_{\al,d}.$
\end{enumerate}
\end{proposition}
\end{paragr}

\begin{paragr} Ce paragraphe est consacré à la preuve de la proposition précédente.

Un point de $\Xgo^T$ est de la forme $\eps^{-\la} K$ pour $\la\in X_\ast(T)$. Ce point est dans la grassmannienne tronquée $\Xgo^T(D)$ si et seulement si $D_M^G(\eps^{-\la})\leq D$ c'est-à-dire $\la\in \Gamma(D)$ puisque $D_M^G(\eps^{-\la})=-(\la)$. L'isomorphisme 1 résulte alors de cette observation et de l'isomorphisme de Chern-Weil.

Ensuite, le théorème \ref{thm:orbites} implique d'une part l'égalité 
$$\Xgo_\gamma(D)_1= \bigcup_{\al\in \Phi^G_+} \Xgo_\gamma^\al(D)$$
et d'autre part que, pour $\al$ et $\beta$ deux éléments distincts de $\Phi^G_+$, on a 
$$ \Xgo_\gamma^\al(D)\cap \Xgo_\gamma^\beta(D)=\Xgo^T(D).$$
L'isomorphisme $2$ s'en déduit immédiatement.

Il reste à prouver le troisième isomorphisme. Il suffit de considérer le cas où $G$ est de rang semi-simple $1$. Dans ce cas, on a soit $M=G$ soit $M=T$. Notons $G_{\scnx}$ le revêtement simplement connexe du groupe dérivé de $G$. C'est un groupe isomorphe à $\SL(2)$. Soit $T_{\scnx}$ le sous-tore maximal de $G_{\scnx}$ qui est l'image réciproque de $T$ par le morphisme canonique $G_{\scnx}\to G$. On note $\tgo_{\scnx}$ l'algèbre de Lie de $T_{\scnx}$. Soit $T_0$ la composante neutre du noyau de l'une des deux racines de $T$ dans $G$. On note $\tgo_0$ son algèbre de Lie. On a alors une décomposition $\tgo=\tgo_{\scnx}\oplus \tgo_0$ suivant laquelle on écrit $\gamma=\gamma'+\gamma_0$.

On pose simplement $\Xgo^{\scnx}=\Xgo^{G_{\scnx}}$. Le morphisme canonique $G_{\scnx}\to G$ induit une injection 
\begin{equation}
  \label{eq:inj-Gsc}
  \phi:\  \Xgo^{\scnx}_{\gamma'} \to \Xgo_\gamma.
\end{equation}
La translation par $\eps^{\la}$, pour $\la\in X_\ast(T)$, composée avec $\phi$ est notée $\phi_\la$.  On note $\Lambda$ un système de représentants du quotient $X_\ast(T) / X_\ast(T_{\scnx})$. Le lemme 8.4 de \cite{GKM} décrit $\Xgo_\gamma$ comme réunion disjointe d'ouverts
\begin{equation}
  \label{eq:reunion-dis}
  \Xgo_\gamma=\bigcup_{\la \in \Lambda} \phi_\la(\Xgo_{\gamma'}^{\scnx}).
\end{equation}

Soit $\la\in \Lambda$. On va décrire $\Xgo_\gamma\cap \phi_\la(\Xgo_{\gamma'}^{\scnx})$. On identifie $a_{T_{\scnx}}^\ast$ à $(a_T^G)^\ast$. Soit $B$ un sous-groupe de Borel de $G_{\scnx}$ . On note $\Sigma_B(1)$ le sous-ensemble de $\Sigma(1)$ formé des cônes $\sigma$ de sorte que la  condition suivante soit remplie
$$\sigma\subset a_B^+ + a_G^\ast \text{ et } \sigma\not\subset a_G^\ast.$$
Soit $\al^\vee_B$ l'unique coracine de $T_{\scnx}$ positive pour $B$. Si $\sigma \in \Sigma_B(1)$, la condition ci-dessus implique que $\varpi_\sigma(\al^\vee_B)>0$.

On considère alors l'éventail $\Sigma'$ dans $X^\ast(T_{\scnx})$ formé du cône $\{0\}$ et des cônes $a_B^+$, pour les sous-groupes de Borel $B$ de $G_{\scnx}$ tels que $\Sigma_B(1)\not=\emptyset$. On pose 
\begin{equation}
  \label{eq:div'}
 n_{a_B^+}(\la)=\min_{\sigma \in \Sigma_B(1)}  E( \frac{n_\sigma- \varpi_\sigma(\la)}{\varpi_\sigma(\al^\vee_B)})
\end{equation}
où $D=\sum_{\sigma \in \Sigma(1)} n_\sigma D_\sigma$ et $E$ désigne la partie entière. On définit alors  un diviseur $\widehat{T_{\scnx}}$-invariant sur la variété torique $Y_{\Sigma'}$ par  $D'(\la)=\sum_{\sigma'\in  \Sigma'(1)}n_{\sigma'}(\la)D_{\sigma'}$. On introduit
\begin{equation}
  \label{eq:LamdaD}
\Lambda(D)=\{\la\in \Lambda \ | \ (\sigma\in \Sigma(1) \ \mathrm{et}\ \sigma\subset a_G^\ast) \Rightarrow  \varpi_\sigma(\la)\leq n_\sigma\}.
\end{equation}
À l'aide du lemme suivant, on vérifie   
\begin{equation}
  \label{eq:inter-phi}
  \Xgo_\gamma\cap \phi_\la(\Xgo_{\gamma'}^{\scnx})= \left\lbrace\begin{array}{l} \phi_\la(\Xgo_{\gamma'}^{\scnx}(D'(\la))) \text{ si }  \la\in \Lambda(D); \\ \emptyset \text{ sinon}. \end{array}\right.
\end{equation}

\begin{lemme} \label{lem:calc-DMG}Soit $x\in G_{\scnx}(F)$ et $\la\in X_\ast(T)$. Par abus, on ne distingue pas dans les notations $x$ et son image dans $G(F)$. Alors
$$D_M^G(\eps^\la x)=(\la) + \sum_{B\in \pc^{G_{\scnx}}(T_{\scnx})}  \varpi_{a_B^+} (H_B(x)) \sum_{\sigma\in \Sigma_B(1)} \varpi_\sigma(\al^\vee_B)D_\sigma.$$
\end{lemme}

\begin{preuve}Par définition,
$$D_M^G(\eps^\la x)=\sum_{\sigma \in \Sigma(1)}\varpi_{\sigma}(H_{P_\sigma}(\eps^\la x))D_\sigma,$$
où $P_\sigma \in \pc^G(T)$ vérifie $\sigma \subset a_{P_\sigma}^+ + a_G^\ast$ (un tel sous-groupe parabolique existe puisque l'éventail $\Sigma$ est $(G,M)$-adapté). De deux choses l'une : soit $\sigma\not\subset a_G^\ast$ auquel  cas $P_\sigma$ est un sous-groupe de Borel. On considère alors $B$ le sous-groupe de Borel de $G_{\scnx}$ image réciproque de $P_\sigma$ et $\al^\vee$ la coracine de $T_{\scnx}$ positive pour $B$. Alors 
$$\varpi_{\sigma}(H_{P_\sigma}(\eps^\la x))=\varpi_{\sigma}(\la)+\varpi_\sigma(\al^\vee) \varpi_{a_B^+} (H_B(x)).$$

Soit $\sigma \subset a_G^\ast$ et dans ce cas $P_\sigma=G$ et $\varpi_{\sigma}(H_{G}(\eps^\la x))=\varpi_\sigma(\la)$.
\end{preuve}

De (\ref{eq:reunion-dis}) et (\ref{eq:inter-phi}), on déduit la réunion disjointe d'ouverts
\begin{equation}
  \label{eq:reunion-dis2}
  \Xgo_\gamma=\bigcup_{\la \in \Lambda(D)} \phi_\la(\Xgo_{\gamma'}^{\scnx}(D'(\la)).
\end{equation}
On pose $\dc_{\scnx}= \dc(T_{\scnx})$ et $\Sc_{\scnx}= \Sc(T_{\scnx})$. Le tore $T_{\scnx}$ opère sur $\Xgo_\gamma(D)$ via le morphisme canonique $T_{\scnx}\to T$. À  partir  de (\ref{eq:reunion-dis2}) et de la proposition \ref{prop:hom-1orb-SL2}, on peut calculer le $\dc_{\scnx}$-module $H^{T_{\scnx}}_\bullet(\Xgo_\gamma(D),\Xgo^T(D))$. On trouve qu'il est isomorphe au sous-$\dc_{\scnx}$-module de $\Gamma(D)\otimes \Sc_{\scnx}$ engendré par
$$(1-\al^\vee)^d (\al^\vee)^n\la^{-1} \otimes \ker(\partial_\al^d),$$
où $\al$ est l'unique élément de $\Phi_+^{G_{ \scnx}}$, $\la\in \Lambda(D)$, $d$ et $n$ sont des entiers qui vérifient les inégalités
$$1\leq d \leq \val(\al(\gamma'))$$
 et 
$$-  n_B(\la)\leq n \leq n+d \leq  n_{\bar{B}}(\la).$$
C'est le sous-module engendré par 
$$(1-\al^\vee)^d \la \otimes \ker(\partial_\al^d),$$
pour $\la\in X_\ast(T)$ tel que l'ensemble
$$\{\la, \al^\vee \la,\ldots,(\al^\vee)^d  \la\}$$
soit inclus dans $\Gamma(D)$. 
On a introduit au début de ce paragraphe un tore $T_0$. On a la décomposition $\Sc=\Sc_{\scnx}\otimes \Sc(T_0)$.  Comme le tore $T_0$ agit trivialement sur $\Xgo_\gamma(D)$, on a
$$H^T_\bullet(\Xgo_\gamma(D),\Xgo^T(D))=H^{T_{\scnx}}_\bullet(\Xgo_\gamma(D),\Xgo^T(D))\otimes \Sc_\bullet(T_0)$$
ce qui donne le résultat cherché. 
\end{paragr}

\begin{paragr} Indiquons un corollaire à la proposition précédente.

  \begin{proposition}\label{prop:descrip-hom}
    Supposons que $\Xgo_\gamma(D)$ soit pure. Alors on a une suite exacte de $\dc$-module 
$$\begin{CD} 
0@>>> \displaystyle \sum_{\al\in \Phi_+^G} \sum_{d=1}^{\val(\al(\gamma))} R_{\al,d} @>>> \Gamma(D)\otimes \Sc_\bullet @>>> H_\bullet^T(\Xgo_\gamma(D)) @>>> 0.
\end{CD}$$
En particulier, l'homologie $T$-équivariante de $\Xgo_\gamma(D)$ s'annule en degré impair.
  \end{proposition}

  \begin{preuve} L'hypothèse de pureté entraîne la validité du lemme de Chang-Skjelbred

$$\begin{CD} 
H_\bullet^T(\Xgo_\gamma(D)_1,\Xgo^T(D))@>>> H_\bullet^T(\Xgo^T(D))@>>>H_\bullet^T(\Xgo_\gamma(D)) @>>> 0.
\end{CD}$$
Le module gradué $\Gamma(D)\otimes \Sc_\bullet$ est isomorphe à $H_\bullet^T(\Xgo^T(D))$ par l'isomorphisme de Chern-Weil qui double les degrés. On en déduit la dernière assertion. L'image de $H_\bullet^T(\Xgo_\gamma(D)_1,\Xgo^T(D))$ s'identifie précisément au sous-module
$$\sum_{\al\in \Phi_+^G} \sum_{d=1}^{\val(\al(\gamma))} R_{\al,d}$$
(cf. proposition \ref{prop:hom-orb-0-1}).
  \end{preuve}
\end{paragr}

\section{Complexes de faisceaux sur les variétés toriques.}\label{sec:cplx-fx}

\begin{paragr} Soit $(G,M,T)$ un triplet formé de $G$ un groupe défini sur $k$, réductif et connexe , de $M$ un sous-groupe de Lévi de $G$ et de $T$ un sous-tore maximal vérifiant $T\subset M\subset G$. Soit  $\gamma\in \tgo(F)$. Soit $\Sigma=\Sigma_M^G$ l'éventail  défini à la ligne (\ref{eq:SGM}) du \S \ref{S:GMfamille} et  $Y=Y_\Sigma$ la variété torique associée (cf. section \ref{sec:vartorique}). Soit  $A$  l'anneau gradué des coordonnées homogènes défini \S\ref{S:coord-hom} l.(\ref{eq:coord-hom}). Rappelons qu'on fait le choix dans ce paragraphe d'une projection $\la \to \la_S$ de $X_\ast(T)$ sur $X_\ast(S)$ (cf. (\ref{eq:projecteur})).

Pour tout $\la\in X_\ast(T)$, on introduit les éléments suivants de $A$ :
\begin{equation}
  \label{eq:monomes}
  y^\lambda_+=\la_S \prod_{\sigma \in \Sigma(1), \varpi_\sigma(\la)\geq 0} y_\sigma^{\varpi_\sigma(\lambda)}\ \ \text{et}\ \ \ y^\lambda_-=\prod_{\sigma \in  \Sigma(1), \varpi_\sigma(\la)< 0} y_\sigma^{-\varpi_\sigma(\lambda)}
\end{equation}
On notera que $y^\lambda_+$ et $y^\lambda_-$ sont homogènes de même degré (pour la graduation de $A$). 

On a défini au \S \ref{S:coh-equiv3} des algèbres $\Sc=\Sc(T)$ et $\dc=\dc(T)$. L'algèbre $\dc$ s'identifie à l'algèbre des opérateurs différentiels à coefficients constants sur l'algèbre de polynômes $\Sc$. Rappelons que tout caractère $\al\in X^\ast(T)$ définit un élément $\partial_\al \in \dc$. On introduit alors $R_\al$ le sous-$A\otimes \dc$-module gradué de $A\otimes \Sc$ 
$$R_\al= \sum_{d=1}^{\val(\al(\gamma))} (y^{\al^\vee}_+-y^{\al^\vee}_-)^d A\otimes \ker(\partial_\al^d).$$

Soit $H$ un groupe réductif connexe défini sur $k$ dont le groupe dual $\Hc$ est un sous-groupe de $\Gc$ qui contient $\Tc$. On définit alors le $A\otimes\dc$-module gradué  $L^H_T$ par
$$L^H_T = A\otimes \Sc / \sum_{\al\in \Phi^H_+} R_\al.$$
Soit $L^H=L^H_T\{\dc^+\}$ le sous-$A$-module de $L^H_T$ annulé par les éléments de l'idéal d'augmentation $\dc^+$ de $\dc$. 
La graduation naturelle sur $\Sc$ induit des graduations $L^H_T=\oplus_{k\in \NN} L^H_{T,k}$ et $L^H=\oplus_{k\in \NN} L^H_k$.
La proposition suivante justifie l'introduction du $A$-module $L^H$.

\begin{proposition} \label{prop:module-hom} Soit $ D\in \Div_{\Tc}(Y)$ tel que la fibre de Springer affine tronquée $\Xgo_\gamma^H(D)$ soit pure au sens de la définition \ref{def:purete}.  On a alors un isomorphisme naturel 
$$H_{2n}(\Xgo_\gamma^H(D),\Qlb) \simeq L^H_n[D],$$
pour tout $n\in \NN$. Par ailleurs, pour tout entier $n$,
$$H_{2n+1}(\Xgo_\gamma^H(D),\Qlb)=0.$$
\end{proposition}

\begin{preuve} Puisque $\Xgo_\gamma^H(D)$ est pure, on a 
$$H_{\bullet}(\Xgo_\gamma^H(D),\Qlb)=H_{\bullet}^T(\Xgo_\gamma^H(D),\Qlb)\{\dc^+\}$$
(cf. \S\ref{S:coh-equiv3}). Comme on a, par ailleurs,  
$$ L^H_\bullet[D]=L^H_{T,\bullet}[D]\{\dc^+\}$$
il suffit d'exhiber un isomorphisme
$$H_{2\bullet}^T(\Xgo_\gamma^H(D),\Qlb)\to L^H_{T,\bullet}[D].$$
On a défini au  \S \ref{S:coord-hom} un sous-tore $S$ de $T$. En particulier, l'anneau $\Qlb[X_\ast(S)]$ est un sous-anneau de $\Qlb[X_\ast(T)]$ et le morphisme injectif de $\Qlb[X_\ast(T)]$ dans $A'$, défini l.(\ref{eq:mor}) de ce paragraphe, permet d'identifier $\Qlb[X_\ast(S)]$ à son image $A[0]$. L'application $\la \mapsto y^D y^\la$ induit un isomorphisme de $A[0]$-module entre $\Gamma(D)=\Gamma(Y,\oc(D))$  et $A[D]$.  On en déduit un isomorphisme de $A[0]\otimes \dc$-module entre $\Gamma(D)\otimes \Sc$ et $A[D]\otimes \Sc$. De plus, on vérifie que cet isomorphisme et son inverse échange, pour chaque racine $\al\in \Phi_+^H$, les sous-modules $\sum_{d=1}^{\val(\al(\gamma))} R_{\al,d}\subset \Gamma(D)\otimes \Sc$ (définis au \S\ref{S:moduleR}) et $ L_\al[D]\subset A[D]\otimes\Sc$. Cet isomorphisme induit donc l'isomorphisme cherché comme on le voit sur le calcul explicite de  $H_{2 \bullet}^T(\Xgo_\gamma^H(D),\Qlb)$ donné à la proposition \ref{prop:descrip-hom}. 
Finalement, l'annulation de l'homologie de $\Xgo_\gamma^H(D)$ en degré impair résulte de l'annulation analogue en homologie $T$-équivariante donnée dans cette même proposition.
\end{preuve}
\end{paragr}

\begin{paragr} Soit $\Hc'\subset \Hc$ des sous-groupes réductifs de $\Gc$ qui contiennent $\Tc$. Soit $H'$ et $H$ des groupes duaux sur $k$ munis d'un plongement de $T$. On pose 

$$\nabla_{H'}^{H}= \prod_{\al^\vee \in \Phi^{\Hc}_+-\Phi^{\Hc'}_+} (y^{\alpha^\vee}_+- y^{\alpha^\vee}_-).$$
C'est un élément homogène de $A$ qui ne dépend du choix d'un système de représentants $\Phi_+^{\Gc}$ qu'à un élément inversible de $A$ près. Soit $V_{H'}^{H}$ l'ouvert de $V \times \hat{S}$ (cf. \S\ref{S:coord-hom}) où $\nabla_{H'}^{H}$ ne s'annule pas. Cet ouvert est stable par $\widehat{\mathrm{Cl}}(Y)$ et le quotient est un ouvert de $Y$ noté  $U_{H'}^{H}$.

À la suite de Goresky-Kottwitz-MacPherson, on introduit le facteur de transfert homologique  
$$\Delta_{H'}^{H}= \prod_{\al^\vee \in \Phi^{\Hc}_+-\Phi^{\Hc'}_+} \partial_\alpha^{\val(\alpha(\gamma))}$$
qui est un élément de $\mathcal{D}$.

Soit $\tilde{L}^H(D)$ le faisceau quasi-cohérent associé à $L^H(D)$ (cf. la fin du \S\ref{S:coord-hom}).

\begin{lemme} \label{lem:iso-lg1} Le facteur de transfert $\Delta_{H'}^H$ induit un isomorphisme de faisceaux 
$$ \tilde{L}^H(D) \longrightarrow  \tilde{L}^{H'}(D)$$
sur l'ouvert $U_{H'}^{H}$.
\end{lemme}

\begin{preuve} On fait appel aux notations des paragraphes \ref{S:vartor-intro} et \ref{S:coord-hom}. Le lemme est de nature locale : il suffit de le prouver sur chaque ouvert $U_{\sigma}\cap U_{H'}^{H}$ pour $\sigma\in \Sigma$. Soit $\sigma\in \Sigma$ et $B=A_\sigma[0]$. Sur l'ouvert affine $U_{\sigma}\simeq \Spec(B)$, le faisceau  $\tilde{L}^H(D)$ est le faisceau quasi-cohérent associé au $B$-module $(L^H(D)\otimes A_\sigma)[0]$. 
Pour chaque racine $\al^\vee\in \Phi^{\Gc}_+$, on a la propriété suivante : les entiers $\varpi_\tau(\al^\vee)$ sont tous de même signe lorsque $\tau$  vérifie les conditions $\tau\in \Sigma(1)$ et $\tau \subset \sigma$. En effet, il suffit de se rappeler que $\sigma=a_P^{G,+}$ pour un élément $P\in\fc^G(M)$ et que $\al^\vee$ ou $-\al^\vee$ est une coracine dans $P$. Quitte à changer le système de représentants $\Phi^{\Gc}_+$ et à prendre $-\al^\vee\in \Phi^{\Gc}_+$, on suppose ces entiers positifs. On en déduit que pour toute  racine $\al^\vee\in \Phi^{\Gc}_+$, le monôme $y_-^{\al^\vee}$ est un élément inversible de $A_\sigma$, que les éléments $y^{\al^\vee}-1$ et 
$$f=\prod_{\al^\vee \in \Phi^{\Hc}_+-\Phi^{\Hc'}_+} (y^{\al^\vee}-1)$$
appartiennent à $B$. Par ailleurs, $f$ et $\nabla^H_{H'}$ sont égaux à un élément inversible de $A_\sigma$ près. L'ouvert $U_{\sigma}\cap U_{H'}^{H}$ est alors isomorphe à $\Spec(B_f)$ où $B_f$ est l'anneau localisé de $B$ par rapport à la partie multiplicative engendré par $f$.
Il s'agit donc de vérifier que le morphisme de $B$-module donné par le facteur de transfert $\Delta_{H'}^H$ 
$$(L^H(D)\otimes A_\sigma)[0]\to (L^{H'}(D)\otimes A_\sigma)[0]$$
devient un isomorphisme après tensorisation par $B_f$. Or ce morphisme s'explicite de la manière suivante 
$$A[D]\otimes \dc / \sum_{\al \in \Phi^H_+} \sum_{d=1}^{\val(\al(\gamma))} (y^{\al^\vee}-1)^d  
A[D]\otimes \ker(\partial_\al^d) \to A[D]\otimes \dc / \sum_{\al \in \Phi^{H'}_+} \sum_{d=1}^{\val(\al(\gamma))} (y^{\al^\vee}-1)^d  
A[D]\otimes \ker(\partial_\al^d).$$
On reprend alors brièvement le raisonnement de la preuve du théorème 10.2 de \cite{GKM}. Le morphisme ci-dessous est clairement surjectif. Il s'agit de prouver que son noyau s'annule après tensorisation par $B_f$. Soit $x\in A[D]\otimes \dc$ un élément qui s'envoie dans le noyau de ce morphisme c'est-à-dire $x$ vérifie
$$\Delta_{H'}^H(x)\in   \sum_{\al \in \Phi^{H'}_+} \sum_{d=1}^{\val(\al(\gamma))} (y^{\al^\vee}-1)^d  
A[D]\otimes \ker(\partial_\al^d).$$
Le facteur de transfert  $\Delta_{H'}^H$ induit une surjection de $\ker(\partial_\al^d)$ dans lui-même pour $\al \in \Phi^{H'}_+$. Ainsi, on peut supposer que $x \in \ker(\Delta_{H'}^H)$. Comme on a
$$\ker(\Delta_{H'}^H)=\sum_{ \al \in \Phi^H_+-\Phi^{H'}_+ }  \ker(\partial_\al^{\val(\al(\gamma))})$$
on voit qu'il existe $n \in \NN$ tel que 
$$f^n x\in \sum_{\al \in \Phi^H_+-\Phi^{H'}_+} \sum_{d=1}^{\val(\al(\gamma))} (y^{\al^\vee}-1)^d A[D]\otimes \ker(\partial_\al^d)$$
ce qui prouve bien que le noyau s'annule après tensorisation par $B_f$.
\end{preuve}
\end{paragr}

\begin{paragr}\label{S:ec} On a fixé au \S \ref{S:coord-hom} en (\ref{eq:projecteur}) une projection $X_\ast(T) \to X_\ast(S)$. Ce choix donne une décomposition non canonique 
$$Y\simeq V//\widehat{\mathrm{Cl}}(Y) \times \hat{S}.$$
On a donc une projection $Y \to \hat{S}$ qui est un morphisme propre puisque le quotient  $V//\widehat{\mathrm{Cl}}(Y)$ est propre. On vérifie que ce morphisme ne dépend pas du choix du projecteur $X_\ast(T) \to X_\ast(S)$. Soit $s_M\in \hat{S}$ et $Y_{s_M}$ la fibre du morphisme 
 $$Y \to \hat{S}$$
au-dessus de $s_M$. 

L'ensemble 
$$\{\Gc_s \ | \ s \text{ appartient à la fibre de } \Tc \to  \hat{S} \text{ au-dessus de } s_M \}$$
est fini. On note $\hat{\ec}=\hat{\ec}(s_M)$ le sous-ensemble de ses éléments maximaux pour l'inclusion. 

Pour $\Hc\in \hat{\ec}$, soit $U_H=U_H^G$ l'ouvert défini au début du paragraphe précédent.

\begin{lemme}\label{lem:recouvt} La réunion des ouverts $U_{H}$ pour $\Hc\in \hat{\ec}$ contient $Y_{s_M}$ tout entier.
\end{lemme}

\begin{preuve} Le choix d'une projection $X_\ast(T) \to X_\ast(S)$ donne une section du morphisme canonique $\Tc \to \hat{S}$. Par abus, on note encore $s_M$ l'élément de $\Tc$ qui est l'image de $s_M$ par cette section. 

La projection sur le second facteur $V\times \hat{S}\to  \hat{S}$ se factorise par le morphisme $Y \to \hat{S}$. Il suffit de prouver que $V\times \{s_M\}$ est inclus dans la réunion des ouverts $V_H$ pour $\Hc\in \hat{\ec}$.  On va prouver la même assertion au niveau des points sur $\Qlb$. Soit $x=(x_\sigma)_{\sigma\in \Sigma(1)}$ un point de $V(\Qlb)$.

Soit  $P\in \pc(M)$ tel que que $x_\sigma\not=0$ pour $\sigma\not\subset a_P^{G,+}$. Soit $\Sigma(1)'$ l'ensemble des $\sigma\in \Sigma(1)$ pour lesquels $x_\sigma\not=0$. Soit $z\in \Tc$ défini par 
$$z=\prod_{\sigma\in\Sigma(1)'} x_\sigma^{\varpi_\sigma}.$$
Pour tout $\chi \in X^\ast(\hat{S})=X_\ast(S)$, on a $\chi(z)=1$. Donc $z\in \ker(\Tc \to \hat{S})$. Soit $s=s_M z$. Par définition de $\hat{\ec}$, il existe $\Hc\in  \hat{\ec}$ tel que $\Gc_s\subset \Hc$. Pour un tel $\Hc$, nous allons montrer que $\nabla_{H}(x)\not=0$ c'est-à-dire que pour toute  racine $\al^\vee$ hors de $\Phi^{\Hc}$, 
 \begin{equation}
   \label{eq:ineg}
  (y^{\alpha^\vee}_+- y^{\alpha^\vee}_-)(x,s_M)\not=0. 
 \end{equation}
Fixons une telle racine $\al^\vee$ : on a donc  
$$\al^\vee(s)\not=1.$$
Quitte à prendre $-\al^\vee$, on peut supposer que $\al \in \Phi^P$. Dans ce cas, si $\varpi_\sigma(\al^\vee) <0$, on a $\sigma\not\subset a_P^{G,+}$ et $x_\sigma\not=0$. On en déduit
$$y_{-}^{\al^\vee}(x,s_M)=\prod_{\sigma \in  \Sigma(1), \varpi_\sigma(\al^{\vee})< 0} x_\sigma^{-\varpi_\sigma(\al^\vee)}\not=0.$$
D'autre part,
$$y_{+}^{\al^\vee}(x,s_M)=\al^\vee(s_M)\prod_{\sigma \in  \Sigma(1), \varpi_\sigma(\al^{\vee})\geq 0} x_\sigma^{\varpi_\sigma(\al^\vee)}.$$
Si $y_{+}^{\al^\vee}(x,s_M)=0$, l'inégalité (\ref{eq:ineg}) est évidente. Sinon, pour tout $\sigma\in \Sigma(1)- \Sigma(1)'$, on a $ \varpi_\sigma(\al^{\vee})=0$ et 
$$(y^{\alpha^\vee}_+- y^{\alpha^\vee}_-)(x,s_M)= y_{-}^{\al^\vee}(x,s_M)^{-1} (\al^\vee(s) -1).$$
ce qui montre que cette expression n'est pas nulle.
\end{preuve}
\end{paragr}

\begin{paragr} Pour tout entier $n$, soit $[n]$ l'ensemble $\{1,\ldots,n\}$. Soit $s_M\in \Tc$ et  $n$ le cardinal de l'ensemble $\hat{\ec}$ (cf. \S \ref{S:ec}). On indexe les éléments de $\hat{\ec}$ par $[n]$ ainsi $\hat{\ec}=\{\Gc_1,\ldots,\Gc_n\}$. Pour toute partie $I\subset[n]$, soit $G_I$ un groupe réductif sur $k$ dont le dual est donné par
$$\Gc_I= \big(\bigcap_{i\in I }\Gc_i \big)^0.$$
Par convention, $G_\emptyset=G$. On pose 
$$\Delta_J^I=\Delta_{G_J}^{G_I}.$$
On remarque que si $I$ est le complémentaire dans $J$ d'un singleton alors le groupe $G_J$ est un groupe endoscopique de $G_I$.

Soit $I\subset [n]$ et $C_I^k$ l'ensemble des parties à $k$ éléments de $I$. Soit  $\Sc^I$ le $\Sc$-module libre de base $(e_i)_{i\in I}$. L'algèbre extérieure $\bigwedge \Sc^I$ est un $\Sc$-module libre de base $(e_K)_{K\subset I}$ avec pour $K\subset I$ 
$$e_K=e_{\phi(1)}\wedge \ldots \wedge e_{\phi(|K|)},$$
où $\phi$ est l'unique bijection croissante de $[|K|]$ sur $K$. On le munit d'une structure de $\mathcal{D}$-module de la façon suivante : pour tout $\partial\in \mathcal{D}$ et tout $a\in \Sc$, $\partial( a e_K)=\partial(a)e_K$.

Pour toute partie $J\subset[n]$ telle que $I\cap J=\emptyset$, on définit une différentielle $d_J$ (de degré $1$ pour la graduation naturelle de $\bigwedge \Sc^I$) par 
$$d_J(a e_K)= \sum_{i\in I} \Delta^{J\cup K}_{J\cup K\cup\{i\}}(a) e_i \wedge e_K,$$
pour toute partie $K\subset I$ et tout $a\in \Sc$.

On a alors un complexe de $\dc$-module de longueur $|I|$ 
$$\begin{CD}
0 @>>> \bigwedge^{0}\Sc^I @>>>  \bigwedge^{1}\Sc^I @>>> \ldots  @>>>  \bigwedge^{|I|}\Sc^I @>>> 0
\end{CD}
$$
En tensorisant par l'anneau $A$, on obtient un nouveau complexe. Pour toute partie $K\subset[n]$, on introduit le $A\otimes\dc$-module  
$$R_K=\sum_{\al\in \Phi_+^{G_K}} R_\al\ ;$$
on a donc 
$$L^{G_K}_T= A\otimes \Sc / R_K.$$
Alors, pour $K\subset I$, la différentielle $1\otimes d_J$ envoie le $A\otimes\dc$-module $L_{K\cup J} e_K$ dans le module $\sum_{i\in I} L_{K\cup J\cup\{i \}} e_i\wedge e_K$. On en déduit donc un complexe de $A\otimes \dc$-modules :

$$\begin{CD}
0 @>>> L_T^{G_J} @>>> \bigoplus_{K\in C^1_I} L_T^{G_{J\cup K}}   @>>> \ldots  
\end{CD}
$$
$$\begin{CD}@>>> \bigoplus_{K\in C^{|I|-1}_I} L_T^{G_{J\cup K}}    @>>>  L_T^{G_{J\cup I}}  @>>> 0
\end{CD}
$$
En prenant le noyau de l'idéal d'augmentation, on obtient un complexe pour les  $A\otimes \dc$-modules $L^{G_{J\cup K}}$. 

On pose 
$$\Fgo^H(D)=\tilde{L}^H(D)$$
où  $\tilde{L}^H(D)$ est le faisceau quasi-cohérent associé à $L^H(D)$ (cf. la fin du \S\ref{S:coord-hom}). On note $\hat{\Fgo}^H(D)$ le complété de  $\Fgo^H(D)$ le long de $Y_{s_M}$. On déduit du complexe précédent un complexe noté $\mathcal{K}(J,I,D)$

$$\begin{CD}
0 @>>>  \hat{\Fgo}^{G_J}(D) @>>> \bigoplus_{K\in C^1_I} \hat{\Fgo}^{G_{J\cup K}}(D)  @>>> \ldots     @>>>   \hat{\Fgo}^{G_{J\cup I}}(D)  @>>> 0
\end{CD}
$$

\begin{theoreme} \label{thm:cplx-fx} Le complexe $\mathcal{K}(\emptyset,[n],D)$
$$\begin{CD}
0 @>>>  \hat{\Fgo}^{G}(D) @>>> \bigoplus_{K\in C^1_{[n]}} \hat{\Fgo}^{G_{K}}(D)  @>>> \ldots     @>>>   \hat{\Fgo}^{G_{[n]}}(D)  @>>> 0
\end{CD}
$$
est exact.
\end{theoreme}

\begin{preuve} Puisque les ouverts $U_{G_i}$ pour $i\in [n]$ recouvrent $Y_{s_M}$ (cf. lemme \ref{lem:recouvt}), il suffit de prouver que la suite ci-dessus est exacte sur chaque ouvert $Y_{s_M} \cap U_{G_i}$, $i\in[n]$. Fixons donc un tel $i$. 

Nous allons en fait montrer que pour toutes parties $J$ et $I$ de $[n]$ telles que $i\in I$ et $I\cap J=\emptyset$, le complexe de faisceaux $\mathcal{K}(J,I,D)$ est exact sur l'ouvert $Y_{s_M} \cap U_{G_i}$. Pour cela, on raisonne par récurrence sur le cardinal de $I$. Lorsque $I=\{i\}$, le complexe $\mathcal{K}(J,\{i\},D)$ est le complexe de longueur $1$ suivant
$$\begin{CD}
  \hat{\Fgo}^{G_J}(D) @>{\Delta_{J\cup \{i\}}^{J}}>> \hat{\Fgo}^{G_{J\cup \{i\}}}(D).
\end{CD}
$$
On va montrer que ce morphisme est un isomorphisme sur l'ouvert $Y_{s_M} \cap U_{G_i}$.  Pour alléger les notations, on pose $H=G_J$ et $H'=G_{J\cup \{i\}}$. On a alors le diagramme suivant d'inclusions 
$$\begin{array}{ccc} \Hc & \subset  & \Gc \\ \cup & &\cup \\ \Hc' & \subset & \Gc_i \end{array}$$
et l'inclusion $\Phi^H_+-\Phi^{H'}_+\subset \Phi^G_+ -\Phi^{G_i}_+$. L'ouvert $U_{H'}^{H}$ contient donc $U_{G_i}$ et sur cet ouvert le lemme \ref{lem:iso-lg1} donne l'isomorphisme 
$$\begin{CD}
  \Fgo^{H}(D) @>>> \Fgo^{H'}(D)
\end{CD}
$$
d'où \emph{a fortiori} un isomorphisme après complétion.

Traitons maintenant le cas où le cardinal de $I$ plus grand que $1$. En particulier, il existe $j\in I$ distinct de $i$. Posons $I'=I-\{j\}$ et $J'=J\cup\{j\}$. On définit un endomorphisme $f$ du $\dc$-module $\bigwedge(\Sc^{I'})$ par 
$$f(a e_K)=\Delta_{J'\cup K}^{J\cup K}(a) e_K$$
pour toute partie $K\subset I'$. On vérifie que $d_{J'}\circ f = f\circ d_J$. En outre, $f$ induit un morphisme entre les complexes $\mathcal{K}(J,I',D)$ et $\mathcal{K}(J',I',D)$. Par hypothèse de récurrence ces complexes sont exacts. Il en donc de même du cône du morphisme induit par $f$. On vérifie que ce dernier est isomorphe à  $\mathcal{K}(J,I,D)$ ce qui achève la récurrence.
\end{preuve}

\end{paragr}

\begin{paragr} On a construit au paragraphe précédent un complexe de $A\otimes \dc$-modules 
$$\begin{CD}
0 @>>> L^{G_J} @>>> \bigoplus_{K\in C^1_I} L^{G_{J\cup K}}   @>>> \ldots     @>>>  L^{G_{J\cup I}}  @>>> 0
\end{CD}
$$
Pour tout $D\in \Div_{\Tc}(Y)$, on en déduit un complexe de $A[0]\otimes \dc$-modules
 
$$\begin{CD}
0 @>>> L^{G}[D]@>>> \bigoplus_{K\in C^1_{[n]}} L^{G_{K}}[D]  @>>> \ldots    @>>>  L^{G_{[n]}}[D] @>>> 0
\end{CD}$$
Supposons que, pour tout $I$, l'homologie de $\Xgo_\gamma^{G_I}(D)$ soit pure. Alors, d'après la proposition \ref{prop:module-hom}, ce complexe se récrit
\begin{equation}
  \label{eq:cplxH}
  \begin{CD}
0 @>>> H_\bullet(\Xgo_\gamma^{G}(D))  @>>> \bigoplus_{K\in C^1_{[n]}}  H_\bullet(\Xgo_\gamma^{G_K}(D))@>>> \ldots    
\end{CD}
\end{equation}
$$\begin{CD} @>>>  H_\bullet(\Xgo_\gamma^{G_{[n]}}(D))    @>>> 0
\end{CD}$$

 Pour tout $I\subset [n]$, le groupe  $S(F)$ agit  sur $\Xgo_\gamma^{G_I}(D)$. Via le morphisme $\la\in X_\ast(S)\mapsto \eps^\la$, on obtient une action du groupe $X_\ast(S)$ sur $\Xgo_\gamma^{G_I}(D)$ et donc une action de l'algèbre $\Qlb[\hat{S}]=\Qlb[X_\ast(S)]$ sur  $H_\bullet(\Xgo_\gamma^{G_I}(D))$. Soit $\mathcal{I}$ l'idéal maximal de $\Qlb[\hat{S}]$ défini par le point $s_M$. On note $\hat{H}_\bullet(\Xgo_\gamma^{G_I}(D))$ le complété $\mathcal{I}$-adique de $H_\bullet(\Xgo_\gamma^{G_I}(D))$.

\begin{theoreme}\label{thm:suite-exacte}  Pour tout $\la\in X_\ast(T)$ soit $D_\la$ le diviseur défini au \S\ref{S:cohosup}. Soit $D'  \in \Div_{\Tc}(Y)$. Il existe un  entier $d$ tel que pour tout $\la\in X_\ast(T)$ tel que $d(\la)\geq d$ on ait : si pour tout $I\subset [n]$, la fibre de Springer affine $\Xgo_\gamma^{G_I}(D'+ D_\la)$ est pure au sens de la définition \ref{def:purete} alors le complexe suivant pour $D=D'+D_\la$
\begin{equation}
  \label{eq:cplxHc}
  \begin{CD}
0 @>>> \Hc_\bullet(\Xgo_\gamma^{G}(D))  @>>> \bigoplus_{K\in C^1_{[n]}}  \Hc_\bullet(\Xgo_\gamma^{G_K}(D))@>>> \ldots 
\end{CD}
\end{equation}
$$\begin{CD} @>>>  \Hc_\bullet(\Xgo_\gamma^{G_{[n]}}(D))    @>>> 0
\end{CD}
$$
obtenu par complétion $\mathcal{I}$-adique à partir du complexe (\ref{eq:cplxH}) est exact.
\end{theoreme}

\begin{preuve} Comme l'homologie est graduée et que les facteurs de transfert sont homogènes, nous allons prouver l'exactitude du complexe degré par degré. Vu l'hypothèse de pureté, on sait que l'homologie de $\Xgo_\gamma^{G_I}(D)$ s'annule en degré impair. On se limite donc aux degrés paris. Si l'on pose 
$$d_K=\deg(\Delta_K),$$
on a en degré $i$ (pair) un complexe
$$ \begin{CD}
0 @>>> \Hc_i(\Xgo_\gamma^{G}(D))  @>>> \bigoplus_{K\in C^1_{[n]}}  \Hc_{i-2d_K}(\Xgo_\gamma^{G_K}(D))  @>>>\ldots  
\end{CD}$$
$$ \begin{CD}
@>>>  \Hc_{i-2d_{[n]}}(\Xgo_\gamma^{G_{[n]}}(D))    @>>> 0
\end{CD}$$

Kazhdan-Lusztig ont prouvé la formule suivante pour la dimension des fibres de Springer 
$$\dim(\Xgo_\gamma^{G_I})=\sum_{\al \in \Phi_+^{G_I}} \val(\al(\gamma)).$$
En particulier, pour tout $D$, la dimension de $\Xgo_\gamma^{G_I}(D))$ est majorée par cette quantité. On en déduit qu'en degré $i$ tel que
$$i>2\sum_{\al \in \Phi_+^{G}} \val(\al(\gamma))$$
tous les termes du complexe sont nuls. Il suffit donc de considérer l'ensemble fini des complexes indexés par $2i$ avec $0\leq i \leq 2\dim(\Xgo_\gamma^{G})$. Fixons un tel $i$. 

Pour tout $I$, la graduation  $L^{G_I}=\oplus_{j\in \NN} L^{G_I}_j$ entraîne des graduations naturelles sur les faisceaux $\Fgo^{G_I}(D)$ et $\hat{\Fgo}^{G_I}(D)$. Le théorème  \ref{thm:cplx-fx} entraîne l'existence d'une suite exacte

\begin{equation}\label{eq:cplx-grad}
\begin{CD}
0 @>>>  \hat{\Fgo}^{G}_i(D) @>>> \bigoplus_{K\in C^1_{[n]}} \hat{\Fgo}^{G_{K}}_{i-d_K}(D)  @>>> \ldots     @>>>   \hat{\Fgo}^{G_{[n]}}_{i-d_{[n]}}(D)  @>>> 0
\end{CD}
\end{equation}
Soit $\Fgo$ un module cohérent sur $Y$ et $\hat{\Fgo}$ son complété le long de $Y_{s_M}$, qui est, rappelons-le, la fibre du morphisme propre 
$$Y\to \hat{S}$$
au-dessus de $s_M$. D'après \cite{EGAIII} \S4 corollaire 4.1.7, il existe un isomorphisme canonique 
\begin{equation}
  \label{eq:iso-formel}
  H^\bullet(\hat{Y}, \hat{\Fgo}) \to  \hat{H}^\bullet(Y,\Fgo)
\end{equation}
où le second membre désigne le complété $\mathcal{I}$-adique de $H^\bullet(Y,\Fgo)$.

Pour chaque $I\subset [n]$ et chaque entier $i$, le théorème  \ref{thm:cohosup} appliqué au $A$-module de type fini $L_i^{G_I}(D')$ donne un entier $\delta(i,I)$ qui vérifie: pour tous $j>0$ et $\la\in X_\ast(T)$ tel que $d(\la)\geq \delta(i,I)$ on a  
$$H^j(Y,\Fgo^{G_I}_i(D'+D_\la))=0$$
et 
$$H^0(Y,\Fgo^{G_I}_i(D'+D_\la))=L^{G_I}_i[D'+D_\la].$$
Puisque $\Xgo_\gamma^{G_I}(D'+D_\la)$ est pure, la proposition \ref{prop:module-hom} donne un isomorphisme
$$H^0(Y,\Fgo^{G_I}_i(D'+D_\la))\simeq H_{2i}(\Xgo_\gamma^{G_I}(D'+D_\la),\Qlb).$$
On en déduit
$$H^j(\hat{Y},\hat{\Fgo}^{G_I}_i(D'+D_\la))= \hat{H}^j(Y,\Fgo^{G_I}_i(D'+D_\la))=0$$
et 
\begin{equation}
  \label{eq:iso-complete}
 H^0(\hat{Y},\hat{\Fgo}^{G_I}_i(D'+D_\la))= \hat{H}^0(Y,\Fgo^{G_I}_i(D'+D_\la))\simeq \Hc_{2i}(\Xgo_\gamma^{G_I}(D'+D_\la),\Qlb).\end{equation}
On pose alors 
$$\delta=\max_{\left\lbrace\begin{array}{l} 0\leq i \leq 2\dim(\Xgo_\gamma^{G})\\ I\subset [n]\end{array}\right\rbrace} \delta(i,I)$$

Pour tout $\la$ tel que $d(\la)\geq \delta$, le foncteur des sections globales transforme donc la suite exacte (\ref{eq:cplx-grad}) en une suite exacte 
$$\begin{CD}
0 @>>>  H^0(\hat{Y},\hat{\Fgo}^{G}_i(D'+D_\la)) @>>>  \ldots     @>>>   H^0(\hat{Y}, \hat{\Fgo}_{i-d_{[n]}}^{G_{[n]}}(D'+D_\la)) @>>> 0
\end{CD}
$$
qui s'identifie à la suite exacte recherchée grâce aux isomorphismes (\ref{eq:iso-complete}).

\end{preuve}

\end{paragr}

\section{Intégrales orbitales pondérées}\label{sec:IOP}

\begin{paragr} Soit $\Fq$ un corps fini et $k$ une clôture algébrique de $\Fq$. On pose $F=\Fq((\eps))$, $L=k((\eps))$, $\of=\Fq[[\eps]]$ et $\of_L=k[[\eps]]$.  Soit $\tau$ l'automorphisme de Frobenius de $\bar{k}$ sur $k$. Soit $\overline{L}$ une clôture séparable de $L$ et $\overline{F}$ la clôture séparable de $F$ dans  $\overline{L}$. On note encore $\tau$ l'automorphisme de $L$ défini par
$$\tau(\sum_{n=n_0}^\infty a_n \eps^n)=\sum_{i=n}^\infty \tau(a_n) \eps^n.$$
Soit $\Gamma=\Gal(\bar{F}/F)$ le groupe de Galois de $\bar{F}$ sur $F$ et $\Gamma_L$ celui de $\bar{L}$ sur $L$.
\end{paragr}

\begin{paragr}[Mesures de Haar.] --- \label{S:Haar} Le groupe $G(F)$ est muni de la mesure de Haar qui donne la mesure $1$ au sous-groupe ouvert compact $K=G(\of)$. 

Pour tout tore défini sur $F$, on note $T^{\nr}$ le sous-tore non ramifié maximal de $T$. Il vérifie $X_\ast(T^{\nr})=X_\ast(T)^{\Gamma_L}$. Via le morphisme $\la\mapsto\eps^\la$, on identifie $X_\ast(T^{\nr})$ à un sous-groupe de $T(L)$. En particulier, $X_\ast(T^{\nr})^\tau$ s'identifie à un sous-groupe discret et cocompact de $T(F)$. On munit $T(F)$ de l'unique mesure de Haar qui vérifie (cf. \cite{GKM} \S15.2, 15.3)
$$\mu(T)=\mes(X_\ast(T^{\nr})^\tau \back T(F))= \frac{|\Cok[ X_\ast(T^{\nr})_\Gamma \to X_\ast(T)_\Gamma   ]|}{|\Ker[ X_\ast(T^{\nr})_\Gamma \to X_\ast(T)_\Gamma   ]|}.$$
Si $T$ est défini sur $\Fq$, on a $\mu(T)=1$.
\end{paragr}

\begin{paragr}\label{S:def-IOP} Soit $G$ un groupe réductif connexe et $T\subset G$ un sous-tore maximal tous deux définis sur $\Fq$. Soit $M$ un facteur de Lévi d'un sous-groupe parabolique de $G$, tous deux définis sur $\Fq$. On suppose $T\subset M$. 

Soit $\Sigma$ un éventail $(G,M)$-adapté dans $X^\ast(T)$. Le groupe engendré par le Frobenius $\tau$  opère sur le $\ZZ$-module $X^\ast(T)$ ainsi que sur $a_T^\ast$.  On suppose que $\Sigma$ est $\tau$-stable c'est-à-dire que l'action de $\tau$ sur $a_{T}^\ast$ permute les cônes dans $\Sigma$. Soit $Y=Y_\Sigma$ la variété torique associée au tore $\Tc$ et à l'éventail $\Sigma$. Soit $D=\sum_{\sigma\in\Sigma(1)} n_\sigma  D_\sigma \in \Div_{\Tc}(Y)$ un diviseur $\tau$-stable au sens où, pour tout rayon $\sigma\in \Sigma(1)$, on a 
$$n_{\tau(\sigma)}=n_{\sigma}.$$

On suppose que la réunion des cônes dans $\Sigma$
$$|\Sigma|=\bigcup_{\sigma\in\Sigma} \sigma$$
forme un sous-$\RR$-espace vectoriel de $a_M^\ast$. 
Soit $S$ le sous-tore de $T$ défini au \S\ref{S:coord-hom} et $a_S=X_\ast(S)\otimes \RR$. Notons que $S$ est défini sur $\Fq$. Pour tout ensemble muni d'une action du Frobenius  $\tau$, on note par un exposant $\tau$ le sous-ensemble des points fixes par $\tau$.
 
Le noyau de l'application linéaire
\begin{equation}
  \label{eq:la-div}
 \la\in a_{T} \mapsto (\la)=\sum_{\sigma\in \Sigma(1)}\varpi_\sigma(\la)D_\sigma \in \Div_{\Tc}(Y)_\RR
\end{equation}
est le sous-espace $a_S $. 

Soit 
$$\Lambda_\Sigma= \{\la \in a_T^\tau/a_S^\tau \ | \ \exists t\in T(L) \ \ (\la)=(H_T(t))\}$$
On remarque que $\Lambda_\Sigma$ est un réseau dans le quotient $a_{T}^\Gamma / a_S^\Gamma$.

On définit alors une application poids sur $G(F)$ par
\begin{equation}
  \label{eq:poids}
v_D(x)=|\{\la\in \Lambda_\Sigma \ | \ (\la)+D_M^G(x)\leq D\}|.
\end{equation}

\begin{lemme}\label{lem:vD}  L'application $v_D$ est à valeurs dans $\NN$. Elle est invariante par translation à gauche par $M(F)$ et à droite par $K$.
\end{lemme}

\begin{preuve} Dans la définition du poids, l'ensemble qui apparaît est l'intersection d'un réseau avec un compact. Il est donc fini. Soit $m\in M(F)$, $x\in G(F)$ et $k\in K$. On a la formule
$$D_M^G(mxk)= \sum_{\sigma\in \Sigma(1)}\varpi_\sigma(H_M(m)) D_\sigma +D_M^G(xk)=(H_M(m))+D_M^G(x)$$
qui donne déjà l'invariance par $K$. Soit $B$ un sous-groupe de Borel de $M$ défini sur $k$ et qui contient $T$. En utilisant la décomposition d'Iwasawa $m=ntk_M$ suivant $M(L)=N_B(L)T(L)M(\of_L)$, on voit que $(H_M(m))=(H_T(t))$. Il en résulte que la classe de $H_M(m)$ dans $a_T^\tau/ a_S^\tau$ appartient à $\Lambda_\Sigma$. L'invariance par $M(F)$ est alors claire. 
\end{preuve}
\end{paragr}

\begin{paragr} Par abus, on note encore $G$ et $M$ les groupes sur $F$ obtenus par extension des scalaires. Soit $\gamma\in \mgo(F)$ un élément semi-simple et $G$-régulier ce qui signifie que son centralisateur $T_\gamma$ dans $G$ est un tore maximal. Notons que $T_\gamma\subset M$. On considère l'intégrale orbitale pondérée convergente 
$$J_D(\gamma)=J_D^G(\gamma)=\int_{T_\gamma(F)\back G(F)}\mathbf{1}_{\ggo(\of)}(\Ad(x^{-1})\gamma) \,v_D(x) \,\frac{dx}{dt}.$$
Elle dépend du choix des mesures de Haar sur $T_\gamma(F)$ et sur $G(F)$ (cf. \S\ref{S:Haar}).

On vérifie que, pour tout $m\in M(F)$,
$$J_D(\Ad(m)\gamma)=J_D(\gamma).$$

 À la suite de Kottwitz, on introduit $B(G)$ l'ensemble des classes de $\tau$-conjugaison dans $G(L)$. Pour tout $F$-tore $U$, soit 
$$\mathfrak{D}(U) =\ker(B(U)\to B(M)).$$

L'isomorphisme 
$$B(U)\simeq X_\ast(U)_\Gamma$$
dû à Kottwitz (cf. \cite{Ko1} ou \cite{Ko2}) permet de définir un accouplement
$$\bg . , . \bd : \ \hat{U}^\Gamma \times B(U) \to \Qlb,$$
où $\Uc$ est le tore sur $\Qlb$ dual de $U$.

Tout $t\in \mathfrak{D}(T_\gamma)$ est la classe de $\tau$-conjugaison  d'un élément de $T_\gamma(L)$ qui s'écrit $m_t^{-1}\tau(m_t)$ avec $m_t\in M(L)$. On pose alors 
$$\gamma_t=\Ad(m_t)\gamma.$$
C'est un élément de $\mgo(F)$ dont la classe de $M(F)$-conjugaison ne dépend pas du choix de $m_t$. L'expression $J_D(\gamma_t)$ est donc définie sans ambiguïté. De plus, l'application $t \mapsto \gamma_t$ définit une bijection de $\mathfrak{D}(T_\gamma)$ sur l'ensemble des classes de $M(F)$-conjugaison dans l'orbite sous $M(\overline{F})$ de $\gamma$.

Pour tout $s\in \Tc_\gamma^\Gamma$, on considère la $s$-intégrale orbitale pondérée 
$$J_D^s(\gamma)=\sum_{t\in \mathfrak{D}(T_\gamma)}  \bg s , t \bd J_D(\gamma_t).$$
Remarquons que si $s\in Z_{\Mc}^\Gamma\subset \Tc_\gamma^\Gamma$, l'accouplement $\bg s , t \bd$ vaut $1$ pour tout $t\in \mathfrak{D}(T_\gamma)$.
\end{paragr}

\begin{paragr}\label{S:tsafquot} On continue avec l'élément $\gamma$ défini au paragraphe précédent. On suppose de plus que $\gamma$ est entier au sens où pour toute racine de $T_\gamma$ dans $G$ la  valuation de $\al(\gamma)$ est positive.  Rappelons la définition de la fibre de Springer affine tronquée (cf. définition \ref{def:fsafftr})
$$\Xgo_\gamma^G(D)=\{x\in G(L)/G(\of_L) \ |\  \Ad(x^{-1}   )\gamma \in \ggo(\of_L) \text{ et }  D_M^G(x)\leq D \}.$$
Comme $D$ est $\tau$-stable, la fibre ci-dessus est stable par $\tau$ et l'ensemble des points $\tau$-fixes admet la description suivante
\begin{equation}
  \label{eq:fsat-tau}
  \Xgo_\gamma^G(D)^\tau=\{x\in G(F)/G(\of) \ |\  \Ad(x^{-1}   )\gamma \in \ggo(\of) \text{ et }  D_M^G(x)\leq D \}.
\end{equation}
Soit $S_\gamma$ le sous-tore de $T_\gamma$ défini par
\begin{equation}\label{eq:Sgamma}
X_\ast(S_\gamma)= \Ker\left( \begin{array}{ccc} X_\ast(T_\gamma) &\to &\Div_{\Tc}(Y) \\ \la &\mapsto & D_M^G(\eps^\la) 
  \end{array}\right).
\end{equation}
Le groupe $S_\gamma(L)$ opère par translation à gauche sur $\Xgo_\gamma^G(D)$. On en déduit une action libre de $X_\ast(S_\gamma)$ sur $\Xgo_\gamma^G(D)$, via le morphisme usuel de $X_\ast(S_\gamma)$ dans $S_\gamma(L)$. Le quotient $X_\ast(S_\gamma^{\nr}) \back \Xgo_\gamma^G(D)$ est un schéma projectif défini sur $\Fq$.

\end{paragr}

\begin{paragr} \label{S:var-IOP} Le morphisme de  restriction $X^\ast(M)\to X^\ast(T_\gamma)$ induit l'inclusion canonique $a_M^\ast \subset a_{T_\gamma}^\ast$. En particulier, pour tout rayon $\sigma\in \Sigma(1)$, le vecteur $\varpi_\sigma\in a_M^\ast$ définit une forme linéaire sur $a_{T_\gamma}$. On a donc une application linéaire 
$$\la\in a_{T_\gamma} \mapsto \sum_{\sigma\in \Sigma(1)}\varpi_\sigma(\la)D_\sigma \in \Div_{\Tc}(Y)_\RR.$$

Rappelons qu'on dispose d'un morphisme 
$$H_{T_\gamma }\ : \ T_\gamma(L) \to a_{T_\gamma}.$$

Soit $\Lambda_\gamma$ le réseau  de $a_T^\tau/a_S^\tau$ défini par
 $$\Lambda_\gamma= \{\la\in a_T^\tau/a_S^\tau \ | \ \exists t\in T_\gamma(F), \ \  (\la)=(H_{T_\gamma}(t))\}.$$
Pour tout $x\in G(F)$, on introduit la variante suivante du poids $v_D$ 
$$v_{\gamma,D}(x)=|\{\la\in \Lambda_\gamma \ | \ (\la)+D_M^G(x)\leq D\}|$$

\begin{lemme}\label{lem:vgamma}
Le poids $v_{\gamma,D}$ est une application sur $G(F)$ à valeurs dans $\NN$, invariante par translation à gauche par $T_\gamma(F)$ et à droite par $K$. 
\end{lemme}

\begin{preuve} Seule l'invariance par $T_\gamma(F)$ n'est pas évidente. Pour tous $t\in T_\gamma(F)$ et $x\in G(F)$, on a la formule
$$D_M^G(tx)=(H_M(t))+D_M^G(x).$$
L'invariance cherchée résulte de la relation $(H_M(t))=(H_{T_\gamma}(t))$. 
\end{preuve}

Contrairement au poids $v_D$, le poids $v_{\gamma,D}$ n'est pas invariant à gauche par $M(F)$. En fait pour $m\in M(F)$, on a
$$v_{\gamma,D}(m^{-1}x)=v_{\gamma,D+D_M^G(m)}(x).$$

\begin{lemme} On a 
$$v_{\gamma,D}(x)=|{\la\in X_\ast(T^{\nr}_\gamma)^\tau /X_\ast(S^{\nr}_\gamma)^\tau  \ | \ D_M^G(\eps^\la x)\leq D}|$$
\end{lemme}

\begin{preuve} Les morphismes $t  \in T_\gamma(F) \mapsto (H_{T_\gamma}(\gamma))$ et $\la\in X_\ast(T^{\nr}_\gamma)^\tau \mapsto D_M^G(\eps^\la)$ ont même image. Pour conclure, il suffit de remarquer que le noyau du second morphisme est  $X_\ast(S^{\nr}_\gamma)^\tau$ (cf.
(\ref{eq:Sgamma})).\\

\end{preuve}

 On définit alors la variante de l'intégrale orbitale pondérée
$$\tilde{J}_D(\gamma)=\int_{T_\gamma(F)\back G(F)}\mathbf{1}_{\ggo(\of)}(\Ad(x^{-1})\gamma) \,v_{\gamma,D}(x) \,\frac{dx}{dt}.$$
Cette intégrale admet une interprétation en terme de comptage.

\begin{lemme} \label{lem:Jtilde}On a l'égalité 
$$\tilde{J}_D(\gamma)= \mu(T_\gamma)^{-1} | X_\ast(S^{\nr}_\gamma)^\tau \back \Xgo_\gamma^G(D)^\tau |     .$$

\end{lemme}

\begin{preuve}  On écrit à l'aide du lemme précédent
  \begin{eqnarray*}
    \tilde{J}_D(\gamma)& =&  \mu(T_\gamma)^{-1}  \int_{X_\ast(T^{\nr}_\gamma)^\tau \back G(F)}\mathbf{1}_{\ggo(\of)}(\Ad(x^{-1})\gamma) \,v_{\gamma,D}(x) \,dx \\
&=&  \mu(T_\gamma)^{-1} |\{ x\in   X_\ast(S^{\nr}_\gamma)^\tau \back G(F) / G(\of) \ | \   \Ad(x^{-1}   )\gamma \in \ggo(\of) \text{ et }  D_M^G(x)\leq D    \}|\\
  \end{eqnarray*} 
On conclut à l'aide de (\ref{eq:fsat-tau}).
\end{preuve}

\end{paragr}

\begin{paragr} \label{S:tr-frob} Soit $\hat{S}^{\nr}_\gamma$ le tore sur $\Qlb$ dual de $S_\gamma^{\nr}$. Soit $s\in (\hat{S}^{\nr}_\gamma)^\Gamma$ d'ordre fini. Pour alléger les notations, on pose
$$X_\gamma=X_\ast(S^{\nr}_\gamma )$$

La représentation $\ell$-adique donnée par
$$\chi \in X_\gamma \mapsto \chi(s)\in \Qlb$$
définit un système local  sur $\Fq$, noté $\mathfrak{L}_s$, sur $X_\gamma\back \Xgo_\gamma^G(D)$. 

Comme $S_\gamma(L)$ agit sur $\Xgo_\gamma^G(D)$ par translation à gauche et que cette action commute à celle de  $X_\gamma$, le groupe $S_\gamma(L)$ agit sur les groupes d'homologie 
$$H_i( X_\gamma\back \Xgo_\gamma^G(D),\mathfrak{L}_s).$$
On montre comme dans \cite{GKM} 15.6 que cette action est semi-simple et qu'on a une décomposition en sous-espaces isotypiques 
\begin{equation}
  \label{eq:isotypique}
 H_i( X_\gamma\back \Xgo_\gamma^G(D),\mathfrak{L}_s)=\bigoplus_{s'} H_i( X_\gamma\back \Xgo_\gamma^G(D),\mathfrak{L}_s)_{s'}
\end{equation}
où $s'$ parcourt la fibre de $\hat{S}_\gamma^{\Gamma_L} \to   \hat{S}^{\nr}_\gamma$ au-dessus de $s$. Un tel $s'$ s'interprète comme un caractère de $S_\gamma(L)$ à valeurs dans $\Qlb$ via le morphisme de Kottwitz $S_\gamma(L)\to X_\ast(S_\gamma)_{\Gamma_L}$ (cf. \cite{Ko2} 7.6).

On pose
$$\mathrm{trace}(\tau, H_\bullet( X_\gamma\back \Xgo_\gamma^G(D), \mathfrak{L}_s))= \sum_{i=0}^{2\dim(\Xgo_\gamma^G(D))} (-1)^i \, \mathrm{trace}(\tau, H_i( X_\gamma\back \Xgo_\gamma^G(D), \mathfrak{L}_s)).$$
On a une définition semblable pour tout sous-espace $\tau$-stable de $H_\bullet( X_\gamma\back \Xgo_\gamma^G(D), \mathfrak{L}_s)$. La proposition suivante est l'analogue (en homologie) du théorème 15.8 de \cite{GKM} pour les intégrales orbitales pondérées.

\begin{proposition}\label{prop:traceFrob} Pour tout $t\in \mathfrak{D}(S_\gamma)$, soit  $m_t$ un élément de $M(L)$ tel que $m_t^{-1}\tau(m_t)$ appartienne à $S_\gamma(L)$ et dont la classe de $\tau$-conjugaison soit celle de $t$.

Soit $s\in \hat{S}_\gamma^\Gamma$ et $\bar{s}$ l'image de $s$ par le morphisme canonique $\hat{S}^\Gamma_\gamma \to (\hat{S}^{\nr}_\gamma)^\Gamma$. On suppose que $\bar{s}$ est d'ordre fini.

On a l'égalité
$$\mathrm{trace}(\tau,H_\bullet( X_\gamma\back \Xgo_\gamma^G(D), \mathfrak{L}_{\bar{s}})_s )= \frac{\mu(T_\gamma)}{\mu(S_\gamma)} \sum_{t\in \mathfrak{D}(S_\gamma)} \bg s ,t \bd^{-1} \tilde{J}_{D +D_M^G(m_t)}(\Ad(m_t)\gamma).$$
\end{proposition}

\begin{preuve} C'est une simple adaptation des arguments de Goresky-Kottwitz-MacPherson (\emph{ibid.}). Soit $t\in S_\gamma(L)$ et $x\in G(L)/G(\of_L)$ dont la classe modulo $X_\gamma$ est fixée par $t\tau$. Il existe alors $\la_x \in X_\gamma$ tel que
$$t\tau(x)=\la_x x$$
et la classe de $\la_x$ dans les co-invariants  $ (X_\gamma)_\Gamma\simeq B(S_\gamma^{\nr})$ ne dépend que de la classe de $x$ modulo $X_\gamma$.

La version homologique de la formule des traces de Grothendieck-Lefschetz donne  l'égalité 

\begin{eqnarray*}\mathrm{trace}(t\tau,H_\bullet(X_\gamma\back \Xgo_\gamma^G(D), \mathfrak{L}_{\bar{s}}))&=&
\sum_{x\in (X_\gamma \back \Xgo_\gamma^G(D))^{t\tau}} \bg \bar{s} ,\la_x \bd\\
&=& \sum_{\lambda \in (X_\gamma)_\Gamma}    \bg \bar{s} ,\lambda \bd^{-1} | (X_\gamma \back \Xgo_\gamma^G(D))^{\eps^\la t\tau}|.
\end{eqnarray*}
Par ailleurs, comme $\tau$ permute les composantes isotypiques de la décomposition (\ref{eq:isotypique}), on a 
$$\mathrm{trace}(t\tau,H_\bullet(X_\gamma\back \Xgo_\gamma^G(D), \mathfrak{L}_{\bar{s}}))= \sum_{s'} \mathrm{trace}(t\tau,H_\bullet(X_\gamma\back \Xgo_\gamma^G(D), \mathfrak{L}_{\bar{s}})_{s'})$$
où la somme est prise sur les $s'\in \hat{S}_\gamma^\Gamma$ qui s'envoient sur $\bar{s}$. Introduisons une constante $c$ 
$$c=|  \Cok[B(S_\gamma^{\nr}) \to B(S_\gamma)]  |= | \Cok[(X_\gamma)_\Gamma \to X_\ast(S_\gamma)_\Gamma].$$

On en déduit 
\begin{eqnarray*} 
 \mathrm{trace}(t\tau,H_\bullet(X_\gamma\back \Xgo_\gamma^G(D), \mathfrak{L}_{\bar{s}})_{s})&=& c^{-1} \! \!\!\!\!\!  \sum_{t \in \Cok[B(S_\gamma^{\nr}) \to B(S_\gamma)]}  \bg s ,t \bd^{-1}\mathrm{trace}(t\tau,H_\bullet(X_\gamma\back \Xgo_\gamma^G(D), \mathfrak{L}_{\bar{s}})) \\
&=& \mu(S_\gamma)^{-1} \sum_{t \in B(S_\gamma)}\bg s ,t \bd^{-1} | (X_\gamma \back \Xgo_\gamma^G(D))^{t\tau}|
\end{eqnarray*}

Si $(X_\gamma \back \Xgo_\gamma^G(D))^{t\tau}\not=\emptyset$, il existe $x\in G(L)$ tel que $x^{-1}t\tau(x)\in G(\of_L)$. En utilisant la décomposition d'Iwasawa $G(L)=P(L)G(\of_L)$ pour $x$, on voit qu'il existe $m\in M(L)$ tel que $m^{-1}t\tau(m)\in M(\of_L)$. On sait bien alors que la classe de $t$ est triviale dans $B(M)$. Dans la somme ci-dessus, on peut donc se restreindre aux  $t\in \mathfrak{D}(S_\gamma)$. Pour un tel $t$, il existe $m_t$ tel que $t=m_t^{-1}\tau(m_t)$. Soit $\gamma'= \Ad(m_t)\gamma$. La translation  $x\mapsto m_t x$ réalise alors une bijection entre les ensembles $(X_\gamma \back \Xgo_\gamma^G(D))^{t\tau}$ et $(X_{\gamma'} \back \Xgo_{\gamma'}^G(D+D_M^G(m_t)))^{\tau}$. On conclut en utilisant le lemme \ref{lem:Jtilde} et le fait que 
$$\mu(T_{\gamma'})=\mu(T_{\gamma}).$$
\end{preuve}
\end{paragr}

\begin{paragr} Soit $t\in T_\gamma(F)$. On a vu dans la démonstration du lemme  \ref{lem:vgamma} que $(H_{T_\gamma}(t))=(H_M(t))$. Dans la démonstration du lemme \ref{lem:vD}, on a dit que la  classe de $H_M(t)$ dans $a_T^\tau/a_S^\tau$ appartient à $\Lambda_\Sigma$. Il s'ensuit que
$$\Lambda_\gamma \subset \Lambda_\Sigma.$$
Soit $\Lambda_\Sigma^\gamma$ un système de représentant  du quotient fini $\Lambda_\Sigma / \Lambda_\gamma$. Pour tout $x\in G(F)$, on a la relation 
$$v_D(x)=\sum_{\la\in \Lambda_\Sigma^\gamma} v_{T,D+(\la)}(x)$$
et donc 
\begin{equation}\label{eq:J-Jtilde}
J_D(\gamma)=\sum_{\la\in \Lambda_\Sigma^T} \tilde{J}_{D+(\la)}(\gamma).
\end{equation}
Soit $m\in M(L)$ tel que $m^{-1}\tau(m)\in T_\gamma(L)$. La conjugaison par $m$ induit alors un isomorphisme entre $T_\gamma(F)$ et $(mT_\gamma m^{-1})(F)$ à partir duquel on voit que $\Lambda_{\Ad(m)\gamma }=\Lambda_\gamma$. Donc la relation ci-dessus est encore vraie pour le même ensemble  $\Lambda_\Sigma^\gamma$ si l'on remplace $\gamma$ par $\Ad(m)\gamma$.

\begin{theoreme} \label{thm:IOP=trace} Soit $\kappa\in \Tc_\gamma^\Gamma$. Soit $s$ et $\bar{s}$ les images de $\kappa^{-1}$ par les morphismes canoniques $\Tc_\gamma^\Gamma \to \hat{S}_\gamma^\Gamma$ et   $\Tc_\gamma^\Gamma \to (\hat{S}_\gamma^{\nr})^\Gamma$. On suppose $\bar{s}$ d'ordre fini.
On a la relation suivante

$$\sum_{\la\in   \Lambda_\Sigma^\gamma} \mathrm{trace}(\tau,H_\bullet( X_\gamma \back \Xgo_\gamma^G(D+(\la)), \mathfrak{L}_{\bar{s}})_s )=  \frac{\mu(T_\gamma)}{\mu(S_\gamma)} |  \ker(X_\ast(S_\gamma)_\Gamma \to X_\ast(T_\gamma)_\Gamma) |J^\kappa_D(\gamma).$$
\end{theoreme}

\begin{preuve} Soit $\la\in   \Lambda_\Sigma^\gamma$. La proposition  \ref{prop:traceFrob} donne l'égalité
$$\mathrm{trace}(\tau,H_\bullet( X_\gamma \back \Xgo_\gamma^G(D+(\la)), \mathfrak{L}_{\bar{s}})_s )= \frac{\mu(T_\gamma)}{\mu(S_\gamma)} \sum_{t\in \mathfrak{D}(S_\gamma)} \bg s ,t \bd^{-1} \tilde{J}_{D_{\la,t}}(\Ad(m_t)\gamma).$$
où l'on a posé $D_{\la,t}=D+(\la)+D_M^G(m_t)$.
En sommant le second membre sur $\la\in   \Lambda_\Sigma^\gamma$, on trouve 
$$\frac{\mu(T_\gamma)}{\mu(S_\gamma)} \sum_{t\in \mathfrak{D}(S_\gamma)} \bg s ,t \bd^{-1} \big[  \sum_{\la\in   \Lambda_\Sigma^\gamma}  \tilde{J}_{D_{\la,t}}(\Ad(m_t)\gamma)\big].$$

D'après l'égalité (\ref{eq:J-Jtilde}) appliquée à $\Ad(m_t)\gamma$, la somme entre crochets vaut
$$J_{D+D_M^G(m_t)}(\Ad(m_t)\gamma)$$
qui n'est autre que
$$J_{D}(\Ad(m_t)\gamma)$$
puisque pour $t\in  \mathfrak{D}(S_\gamma)$, on vérifie que $D_M^G(m_t)\in \Lambda_\Sigma$.

Pour conclure, il nous reste à prouver 
\begin{equation}
  \label{eq:D-D'}
  |  \ker(X_\ast(S_\gamma)_\Gamma \to X_\ast(T_\gamma)_\Gamma) |J^{\kappa}_D(\gamma)= \sum_{t\in \mathfrak{D}(S_\gamma)} \bg s ,t \bd^{-1} J_{D}(\Ad(m_t)\gamma).
\end{equation}
Soit $T_{\scnx}$ l'image réciproque de $T_\gamma$ dans le revêtement simplement connexe du groupe dérivé de $M$. Il est bien connu que
\begin{equation}\label{eq:Dsc}
\mathfrak{D}(T_{\scnx}) \to \mathfrak{D}(T_\gamma)
\end{equation}
est surjective. Or $X_\ast(T_{\scnx})$ est le sous-groupe de $X_\ast(T_\gamma)$ engendré par les coracines de $T_\gamma$ dans $M$. On en déduit l'inclusion
$$X_\ast(T_{\scnx})\subset X_\ast(S_\gamma).$$
Par conséquent, l'application (\ref{eq:Dsc}) se factorise par $\mathfrak{D}(S_\gamma)$ et l'application
$$\mathfrak{D}(S_\gamma) \to \mathfrak{D}(T_\gamma)$$
est surjective. Son noyau est fini et s'identifie à
$$\ker(X_\ast(S_\gamma)_\Gamma \to X_\ast(T_\gamma)_\Gamma).$$
Si $t\in \mathfrak{D}(S_\gamma)$ a pour image $\bar{t}$ dans $\mathfrak{D}(T_\gamma)$, on a $m_{\bar{t}}=m_t$ et $\bg s ,t \bd^{-1}=\bg \kappa ,\bar{t} \bd$.
L'égalité (\ref{eq:D-D'}) est alors claire.
\end{preuve}

\end{paragr}

\begin{paragr} \label{S:trad2} On a fixé au \S \ref{S:def-IOP}  un sous-tore maximal $T$ de $M$ défini sur $\Fq$. Soit $P$ un sous-groupe parabolique de $G$ défini sur $\Fq$ ayant $M$ comme facteur de Lévi. Soit $B$ un sous-groupe de Borel de $G$ défini sur $k$ qui contient $T$ et qui est inclus dans $M$.

Le Frobenius $\tau$ opère sur la donnée radicielle associée à $(G,T)$
$$(X^\ast(T),\Phi^G,X_\ast(T),\Phi^{G,\vee})$$
et il préserve la donnée radicielle associée à $(M,T)$. 

Soit $s\in \Tc^\tau$. Soit ${H^\ast}$ un groupe déployé sur $\Fq$ et $T^{H^\ast}$ un  sous-tore maximal de ${H^\ast}$  défini sur $\Fq$ de sorte qu'on ait
$$(X^\ast(T^{H^\ast}),\Phi^{H^\ast},X_\ast(T^{H^\ast}),\Phi^{{H^\ast},\vee})=(X^\ast(T),\Phi^M_{s},X_\ast(T),\Phi^{M,\vee}_{s})$$
où l'on a posé
$$\Phi^{M,\vee}_{s}=\{ \al^\vee \in \Phi^{M,\vee}\ |\ \al^{\vee}(s)=1\}$$
et
$$\Phi^{M}_{s}=\{ \al\in \Phi^M\ |\ \al^\vee\in   \Phi^{M,\vee}_{s}\}.$$
Puisque $s$ est fixé par le Frobenius $\tau$, ce dernier stabilise $\Phi^M_{s}$ et 
$\Phi^{M,\vee}_{s}$. On en déduit une action du Frobénius $\tau$ sur la donnée radicielle de ${H^\ast}$. 

Soit $\Delta^{H^\ast}\subset \Phi^{H^\ast}$ la base formée du sous-ensemble des racines simples dans $\Phi^{H^\ast}\cap\Phi^{B}$. Le Frobenius $\tau$ ne préserve pas nécessairement $\Delta^{H^\ast}$. Cependant il envoie $\Delta^{H^\ast}$ sur la base formée de l'ensemble des racines simples dans $\Phi^{{H^\ast}}\cap \Phi^{\tau({B})}$. Soit $w_\tau$ l'unique élément du groupe de Weyl du système de racines $\Phi^{H^\ast}$ qui vérifie 
$$(w_\tau \circ\tau)\Delta^{H^\ast}=\Delta^{H^\ast}.$$
Soit $B^{H^\ast}$ le sous-groupe de Borel de $H^\ast$ qui contient $T^{H^\ast}$ et qui est défini par la base $\Delta^{H^\ast}$. Pour tout $\al\in \Delta^{H^\ast}$, soit $E_\al^{H^\ast}$   un vecteur radiciel de poids $\al$. Soit $\Psi^{H^\ast}$ l'unique automorphisme de ${H^\ast}$ qui stabilise l'épinglage $(B^{H^\ast},T^{H^\ast},\{E_\al^{H^\ast}\}_{\al\in \Delta^{H^\ast}})$ et qui induit sur la donnée radicielle de ${H^\ast}$ l'automorphisme $w_\tau\circ\tau$.

Soit $H$ le groupe sur $\Fq$ défini comme la forme extérieure de $H^\ast$ obtenue par torsion par l'automorphisme $\Psi^{H^\ast}$. L'isomorphisme $T\to T^{H^\ast}$ sur $k$ induit un plongement $\iota : T\to H$. Il résulte des définitions que 
$$\iota\circ\tau=w_\tau^{-1} \circ \tau \circ\iota.$$
D'après le théorème de Kottwitz-Steinberg (cf. \cite{Ko}), quitte à conjuguer $\iota$ par un élément de $H$, on peut et on va supposer qu'on a

\begin{equation}
  \label{eq:tau-equiv}
\iota\circ\tau=\tau \circ \iota.
\end{equation}

Rappelons qu'on a a défini au \S\ref{S:def-IOP} un sous-tore $S$ de $T$ défini sur $\Fq$. On a donc un morphisme $\tau$-équivariant $\Tc\to \hat{S}$. Pour tout $s'\in s\Ker(\Tc \to \hat{S})$, on pose
$$\Phi^\vee_{s'}=\{ \al^\vee \in \Phi^{\vee}\ |\ \al^{\vee}(s')=1\}.$$
L'ensemble
$$\{\Phi^\vee_{s'}\ | \ s'\in s \Ker(\Tc \to \hat{S})\}$$
est fini et non vide. Soit $\ec=\ec_s$ l'ensemble des éléments de cet ensemble maximaux pour l'inclusion. Soit $n$ le cardinal de $\ec$. On fixe 
$$i\in [n] \to \Phi^\vee_i\in \ec$$
une bijection de $[n]$ sur $\ec$. Puisque $s$ est fixé par le Frobenius $\tau$, ce dernier stabilise $\ec$ et définit une permutation de $[n]$, par abus encore notée $\tau$, de la manière suivante
$$\Phi^\vee_{\tau(i)}=\tau(\Phi^\vee_{i}).$$
Soit $I$ une partie non vide de $[n]$ et
$$\Phi^\vee_I=\bigcap_{i\in I} \Phi^\vee_i.$$
Soit $\Phi_I\subset \Phi^G$ d'image $\Phi^\vee_I$ par l'application qui à une racine associe sa coracine.\\

Soit $G_I^\ast$ un groupe réductif connexe et déployé sur $\Fq$ et $T_I^\ast\subset G_I^\ast$ un sous-tore maximal défini sur $\Fq$ tels que la donnée radicielle associée à $(G_I^\ast,T_I^\ast)$ soit égale à
$$(X^\ast(T),\Phi_I,X_\ast(T),\Phi^\vee_I).$$
Soit $\Delta_I$ la base du système de racines $\Phi_I$ formée de l'ensemble des racines simples dans $\Phi_I\cap\Phi^{B}$. Soit $B_I^\ast$ le sous-groupe de Borel de $G_I^\ast$ qui contient  $T_I^\ast$ et qui est défini par $\Delta_I$. Soit $(B_I^\ast,T_I^\ast,\{E_\al^I\}_{\al\in \Delta_I})$ un épinglage .

\begin{lemme} On a la relation
\begin{equation}
  \label{eq:w-delta}
w_\tau\circ \tau (\Delta_I)=\Delta_{\tau(I)}.
\end{equation}
\end{lemme}

\begin{preuve} L'élément  $w_\tau$, vu comme élément  du groupe de Weyl  $W^{\Mc}(\Tc)$, fixe tous les points de $s\Ker(\Tc \to \hat{S})$. En effet, par définition, $w_\tau$ fixe $s$. D'autre part, on vérifie que $\Ker(\Tc \to \hat{S})\subset Z_{\Mc}$. 

On en déduit qu'on a  
$$w_\tau\circ\tau(\Phi_I)=\Phi_{\tau(I)}.$$
On a $\Phi^{B}=\Phi^{M\cap {B}}\cup \Phi^P$. Donc on a aussi
$$\Phi_I\cap \Phi^B= (\Phi_I\cap\Phi^{M\cap {B}}) \cup (\Phi_I\cap\Phi^P).$$
Pour tout $I$, on a 
$$\Phi_I\cap\Phi^{M\cap {B}}=\Phi^{H^\ast} \cap \Phi^B.$$
Par conséquent, par définition de $w_\tau$, on sait que $w_\tau\circ\tau$ stabilise 
$$\Phi_I\cap\Phi^{M\cap {B}}=  \Phi_{\tau(I)}\cap\Phi^{M\cap {B}}. $$
Comme $w_\tau$ se décompose en symétrie élémentaires associées à des racines dans $M$, le composé $w_\tau\circ\tau$ stabilise aussi $\Phi^{P}$ et envoie  $\Phi_I\cap \Phi^P$ sur  $\Phi_{\tau(I)}\cap \Phi^P$. La relation (\ref{eq:w-delta}) est alors évidente. 
\end{preuve}

Du lemme précédent, on déduit que $w_\tau \circ \tau$ induit un isomorphisme de la donnée radicielle de $G^\ast_I$ sur celle de $G^\ast_{\tau(I)}$ qui envoie $\Delta_I$ sur $\Delta_{\tau(I)}$. Il existe donc un unique isomorphisme
$$\Psi_I \ :\ G^\ast_I \to G^\ast_{\tau(I)}$$
qui envoie l'épinglage  $(B_I^\ast,T_I^\ast,\{E_\al^I\}_{\al\in \Delta_I})$ sur celui associé à $G^\ast_{\tau(I)}$ et qui induit sur les données radicielles correspondantes le même isomorphisme que $w_\tau\circ\tau$. 

Soit
$$\mathbf{G}^\ast_I=\coprod_{J} G_J^\ast$$
où $J$ parcourt l'ensemble des parties de $[n]$ qui sont dans l'orbite de $I$ sous l'action de $\tau$. On munit $\mathbf{G}^\ast_I$ de l'automorphisme $\mathbf{\Psi}_I$ défini  sur la composante d'indice $J$ par $\Psi_J$. On note alors $G_I$ la forme extérieure de $\mathbf{G}^\ast_I$ tordue par $\mathbf{\Psi}_I$. Si $\tau(I)=I$ alors $G_I$ est un groupe réductif connexe sur $\Fq$ et c'est la forme extérieure de $G_I^\ast$ tordue par $\Psi_I$.

Soit $M_I^\ast$ l'unique sous-groupe de Lévi de $G^\ast_I$ qui contient $T_I^\ast$ et dont la donnée radicielle est égale à
$$(X^\ast(T),\Phi^{H^\ast},X_\ast(T),\Phi^{{H^\ast},\vee}).$$
Remarquons l'égalité $\Psi_I(M_I^\ast)=M_{\tau(I)}^\ast$. On  définit alors de manière évidente  des sous-variétés $\mathbf{M}^\ast_I\subset \mathbf{G}^\ast_I$ et $M_I\subset G_I$. Le morphisme évident $H^\ast \to M^\ast_I$ permet de définir une action (par translation à gauche) définie sur $\Fq$ de $H$ sur $M_I$. Si $\tau(I)=I$, ce morphisme induit un $\Fq$-isomorphisme de $H$ sur $M_I$.

Rappelons qu'on a défini un plongement $\iota$ de $T$ dans $H$ défini sur $\Fq$. On en déduit une action de $T$ sur $G_I$ défini sur $\Fq$.
\end{paragr}

\begin{paragr}  Soit $\Sigma=\Sigma_M^G$ l'éventail défini en (\ref{eq:SGM}). C'est un éventail $\tau$-stable. Soit $D\in \Div_{\Tc}(Y)$ un diviseur $\tau$-stable. Soit $I\subset[n]$ non vide. On sait que l'éventail $\Sigma$ est $(G_I^\ast,M_I^\ast)$-adapté au sens de la définition \ref{def:adapte} (cf. proposition \ref{eventail-adapte}). Soit $\gamma\in \tgo(F)$ un élément semi-simple, $G$-régulier et entier (cf. début du \S \ref{S:tsafquot}). Par abus, on note encore $\gamma$ l'image de $\gamma$ dans $\hgo(F)$ par le plongement $\iota$ ou son image dans $\hgo_I^\ast$ obtenue par composition avec l'isomorphisme $\hgo \to \hgo^\ast_I$ sur $k$. \\

On introduit alors la grassmannienne affine tronquée 
$$\Xgo^{G_I}_\gamma(D)=\coprod_{J} \Xgo^{G^\ast_J}_\gamma(D)\ ;$$
dans la somme, $J$ parcourt l'orbite de $I$ sous $\tau$ et pour un tel $J$
$$ \Xgo_\gamma^{G^\ast_J}(D)=\{x\in G_J^\ast(L) / G_J^\ast(\of_L)  \ |\ \Ad(x^{-1})\gamma \in \ggo_J^\ast(\of_L) \ \text{ et } \ D_{M_J^\ast}^{G_J^\ast}(x)\leq D\},$$
L'action de $\tau$ sur $G_I$ induit une action de $\tau$ sur 
$$\coprod_{J} G_J^\ast(L) / G_J^\ast(\of_L)$$
et on vérifie que cette action stabilise $\Xgo^{G_I}_\gamma(D)$.

Soit $\bar{s}\in \hat{S}^\tau$ l'image de $s$ par le morphisme $\Tc\to \hat{S}$. On suppose que $\bar{s}$ est d'ordre fini. Comme au \S \ref{S:tr-frob}, la représentation  $\ell$-adique de $X_\ast(S)$ associée à $\bar{s}$ fournit un système local $\mathfrak{L}_{\bar{s}}^I$ sur le quotient $X_\ast(S)\back \Xgo_\gamma^{G_I}(D)$ défini sur $\Fq$.

Soit $\al\in \Phi^G$. On dit que cette racine est symétrique si $\tau(\al)=-\al$. 
On définit alors  un signe 
$$\Delta_{s}(\gamma)=\prod_{\al} (-1)^{\val(\al(\gamma))}$$
où $\al$ parcourt un système de représentants des orbites symétriques de $\tau$ dans $\Phi^M$ tel que $\al^\vee(s)\not=1$.

On pose 
$$d_I=\sum_{\al\in \Phi_+^G-\Phi_+^{G_I^{\ast}}}\val(\al(\gamma)).$$
On remarque que $d_I$ ne dépend que de l'orbite sous $\tau$ de $I$. On note $C_{[n]}^i/\tau$ l'ensemble des orbites de $\tau$ dans $C_{[n]}^i$.

\begin{theoreme} \label{thm:LFP} Soit $\kappa\in \Tc^\tau$ dont l'image dans $\hat{S}$ est d'ordre fini et $D$ un diviseur $\Gamma$-stable. On suppose que pour tout $\la$ dans l'ensemble fini $\Lambda_\Sigma^\gamma$ les assertions suivantes sont vérifiées :
  \begin{enumerate}
  \item pour tout $I\subset \ec_{\kappa}$, la fibre de Springer affine tronquée $\Xgo_\gamma^{G^\ast_I}(D+(\la))$ est pure au sens de la définition \ref{def:purete} ;
  \item le diviseur  $D+(\la)$ satisfait les conclusions du théorème  \ref{thm:suite-exacte}.
   \end{enumerate}
Alors la relation suivante est vraie
$$\Delta_\kappa(\gamma) J_D^{G,\kappa}(\gamma)= \sum_{I}(-1)^{|I|-1} q^{d_I} J_D^{G_I,\kappa}(\gamma)$$
où la somme est prise sur les parties $I$ non vides de $\ec_{\kappa}$ telles que $\tau(I)=I$.
\end{theoreme}

\begin{preuve} Soit $D$ un diviseur $\tau$-stable qui satisfait les hypothèses du théorème  \ref{thm:suite-exacte} et qui vérifie : pour tout $I\subset \ec_{\kappa}$, la fibre de Springer affine tronquée $\Xgo_\gamma^{G^\ast_I}(D)$ est pure. Soit l'anneau $\mathcal{A}=\Qlb[\hat{S}]$ et son idéal maximal $\mathcal{I}$ défini par le point $\bar{s}\in \hat{S}$, image de $s=\kappa^{-1}$ par le morphisme canonique $\Tc\to\hat{S}$.  Pour tout $\mathcal{A}$-module $V$, on note $\hat{V}$ le complété $\mathcal{I}$-adique de $V$. On dispose de l'ensemble fini $\ec_s=\ec_{\kappa^{-1}}$ qu'on met en bijection avec $[n]$. Rappelons qu'on a défini au paragraphe \ref{S:ec} un ensemble fini $\hat{\ec}(\bar{s})$. L'application qui à un sous-groupe de $\Gc$ qui contient $\Tc$ associe l'ensemble des racines de $\Tc$ dans ce sous-groupe envoie $\hat{\ec}(\bar{s})$ bijectivement sur $\ec_{s}$. \\
Pour alléger les notations, on pose pour tous $i\in \NN$ et $I\subset [n]$ (avec la convention $G_\emptyset=G$)
$$H_i^{G_I}=H_i(\Xgo_\gamma^{G_I}(D)).$$

 On a la décomposition
$$ H_i^{G_I}=\bigoplus_J H_i(\Xgo^{G^\ast_J}_\gamma(D)),$$
où $J$ parcourt l'orbite de $I$ sous l'action de $\tau$. Le Frobenius $\tau$ envoie  $H_i(\Xgo^{G^\ast_J}_\gamma(D))$ sur $H_i(\Xgo^{G^\ast_{\tau(J)}}_\gamma(D))$. La trace de $\tau$ sur $ H_i^{G_I}$ est donc nulle sauf si $\tau(I)=I$. 

On suppose désormais que $\tau(I)=I$. On fait l'hypothèse que $\Xgo^{G^\ast_I}_\gamma(D)$ est pure. Il résulte du calcul explicite de $H_i^{G_I}$ donné à la section \ref{sect:Hom_orb} que $H_i^{G_I}$ est soit nul si $i$ est impair soit un sous-espace d'un certain quotient de $\Qlb[X_\ast(T)]\otimes \Sym^{i/2} (X_\ast(T) \otimes \Qlb(1))$ si $i$ est pair. Soit $\chi\otimes X \in \Qlb[X_\ast(T)]\otimes \Sym^{i/2} (X_\ast(T) \otimes \Qlb(1))$. L'action de $\tau$ sur  $H_i^{G_I}$ résulte de l'action de $\tau$ sur $\Qlb[X_\ast(T)]\otimes \Sym^{i/2} (X_\ast(T) \otimes \Qlb(1))$ donnée par 
$$\tau(\chi\otimes X)=\tau(\chi) \otimes q^{i/2} X.$$

 Pour tout $\mathcal{A}$-module $V$, l'isomorphisme
$$\hat{V}\otimes \hat{\mathcal{A}}/ \mathcal{I}\hat{\mathcal{A}}\simeq V\otimes \mathcal{A}/\mathcal{I},$$
induit pour tout entier $p$ un isomorphisme
 $$\Tor^{\hat{\mathcal{A}}}_{p}(\hat{V},\hat{\mathcal{A}}/ \mathcal{I}\hat{\mathcal{A}})\simeq  \Tor^{\mathcal{A}}_p(V,\mathcal{A}/\mathcal{I}).$$
On en déduit pour tous entiers $i$ et $p$ un isomorphisme $\tau$-équivariant 
\begin{equation}
  \label{eq:isotor-torc}
\Tor^{\hat{\mathcal{A}}}_{p}(\hat{H}^{G_I}_i,\hat{\mathcal{A}}/ \mathcal{I}\hat{\mathcal{A}})\simeq  \Tor^{\mathcal{A}}_p(H_i^{G_I},\mathcal{A}/\mathcal{I}).
\end{equation}

La conclusion du théorème \ref{thm:suite-exacte} pour le diviseur $D$ est l'exactitude en tout degré $i$ du complexe de $\hat{\mathcal{A}}$-modules
$$  \begin{CD}
0 @>>> \Hc_i^{G}  @>>> \bigoplus_{I \in C^1_{[n]}/\tau}  \Hc_{i-2d_I}^{G_K}@>>> \ldots     @>>>  \Hc_{i-2d_{[n]}}^{G_{[n]}}    @>>> 0
\end{CD}$$
où les flèches sont données par les facteurs de transfert.
Pour $I\subset \ec_s$ tel que   $\tau(I)\not=I$, on a
$$\trace(\tau,\Tor^{\hat{\mathcal{A}}}_{p}(\Hc_{i}^{G_I},\hat{\mathcal{A}}/ \mathcal{I}\hat{\mathcal{A}}))=0.$$
Si $\tau(I)=I$ et $I\not=\emptyset$, on vérifie que
$$\tau(\Delta_{G^\ast_I}^G)= \Delta_s(\gamma) \Delta_{G^\ast_I}^G.$$
De ces deux remarques et de la suite exacte ci-dessus, on déduit la relation
$$\sum_{I\subset \ec_s,\ \tau(I)=I}(-1)^{|I|} e_I\,   q^{d_I} \sum_{p=0}^\infty (-1)^p   \trace(\tau, \Tor^{\hat{\mathcal{A}}}_{p}(\Hc_{i-2d_I}^{G_I},\hat{\mathcal{A}}/ \mathcal{I}\hat{\mathcal{A}})  )=0,$$
où $e_I= \Delta_s(\gamma)$ si $I\not=\emptyset$ et $e_\emptyset=1$.
Vu l'isomorphisme (\ref{eq:isotor-torc}), on a aussi 
\begin{equation}
  \label{eq:trace-tor}
\sum_{I\subset \ec_s,\ \tau(I)=I}(-1)^{|I|} e_I\,   q^{d_I} \sum_{p=0}^\infty (-1)^p   \trace(\tau, \Tor^{\mathcal{A}}_{p}(H_{i-2d_I}^{G_I},\mathcal{A}/ \mathcal{I})  )=0.
\end{equation}

La proposition \ref{prop:suitespec} implique l'existence d'une suite spectrale 
$$E_{pq}^{2}=\Tor^{\mathcal{A}}_{p}(H_{i}^{G_I},\mathcal{A}/ \mathcal{I}) \Rightarrow H_{p+q}(X_\ast(S) \back \Xgo_\gamma^{G_I}(D),\mathfrak{L}_{\bar{s}}^I).$$
De cette suite spectrale, on déduit la relation suivante
\begin{equation}
  \label{eq:trace-ss}
 \trace(\tau,H_\bullet(X_\ast(S)\back \Xgo_\gamma^{G_I}(D),\mathfrak{L}_{\bar{s}}^I))= \sum_{i=0}^\infty (-1)^i \sum_{p=0}^\infty (-1)^p \trace(\tau,\Tor^{\mathcal{A}}_{p}(H_{i}^{G_I},\mathcal{A}/ \mathcal{I})),
\end{equation}
où les sommes sont bien sûr à support fini.\\

En combinant (\ref{eq:trace-tor}) et (\ref{eq:trace-ss}), on obtient 
$$\sum_{I\subset \ec_s,\ \tau(I)=I}(-1)^{|I|} e_I \,  q^{d_I}  \trace(\tau,H_\bullet(X_\ast(S)\back \Xgo_\gamma^{G_I}(D),\mathfrak{L}_{\bar{s}}^I))=0.$$

On conclut en utilisant le théorème \ref{thm:IOP=trace}.
\end{preuve}

\end{paragr}

\begin{paragr} Soit $\la \in X_\ast(T)^\tau$. Le diviseur $D_\la$ (cf. la fin du \S\ref{S:div-GM})  est $\tau$-stable. On a également définit un entier $d(\la)$ (cf. \S\ref{S:dla}).

\begin{corollaire} \label{cor:LFP} Soit $\gamma\in\tgo(F)$ un élément semi-simple, $G$-régulier, entier et équivalué. Soit $\kappa\in \Tc^\tau$ dont l'image dans $\hat{S}$ est d'ordre fini. Il existe une constante $c>0$ tel que pour tout $\la\in X_\ast(T)^\tau$ tel que $d(\la)\geq c$ on ait
 $$\Delta_\kappa(\gamma) J_{D_\la}^{G,\kappa}(\gamma)= \sum_{I\subset \ec_\kappa, \ \tau(I)=I}(-1)^{|I|-1} q^{d_I} J_{D_\la}^{G_I,\kappa}(\gamma).$$
\end{corollaire}

\begin{preuve}  Soit $\la_1, \ldots,\la_r$ une famille finie d'éléments de $a_T^\tau$ dont l'ensemble des classes dans $a_T^\tau/a_S^\tau$ est $\Lambda_\Sigma^\gamma$. Pour tout indice $i$, soit $t_i\in T(L)$ tel que $(\la_i)=(H_T(t_i))$. Soit $I\subset \ec_\kappa$. D'après le corollaire \ref{cor:purete}, il existe $c_I$ tel que pour tout $\la\in X_\ast(T)$ qui vérifie $d(\la)\geq c_I$, la fibre de Springer affine tronquée $\Xgo_\gamma^{G^\ast_I}(D_\la)$ est pure. La translation à gauche par $t_i$ est un isomorphisme de $\Xgo_\gamma^{G^\ast_I}(D_\la)$ sur $\Xgo_\gamma^{G^\ast_I}(D_\la+(\la_i))$. Par conséquent, cette dernière est également pure. On voit que les hypothèses du théorème \ref{thm:suite-exacte} sont vérifiées. Donc ses conclusions valent et on peut appliquer le théorème \ref{thm:LFP}.
\end{preuve}

\end{paragr}

\bibliographystyle{plain}

\bibliography{troncature-biblio}

\begin{thebibliography}{10}

\bibitem{dis_series}
J.~Arthur.
\newblock The characters of discrete series as orbital integrals.
\newblock {\em Invent. Math.}, 32:205--261, 1976.

\bibitem{inv-trace}
J.~Arthur.
\newblock The invariant trace formula. {I}. {L}ocal theory.
\newblock {\em J. Amer. Math. Soc.}, 1(2):323--383, 1988.

\bibitem{localtrace}
J.~Arthur.
\newblock A local trace formula.
\newblock {\em Publ. Math., Inst. Hautes \'Etudes}, 73:5--96, 1991.

\bibitem{STF1}
J.~Arthur.
\newblock A stable trace formula. {I}. {G}eneral expansions.
\newblock {\em J. Inst. Math. Jussieu}, 1(2):175--277, 2002.

\bibitem{BB}
A.~Beilinson and V.~Drinfeld.
\newblock Quantization of {H}itchin's integrable system and {H}ecke
  eigensheaves.
\newblock Prépublication.

\bibitem{Chang}
T.~Chang and T.~Skjelbred.
\newblock The topological {S}chur lemma and related results.
\newblock {\em Ann. of Math. (2)}, 100:307--321, 1974.

\bibitem{Cox}
D.~Cox.
\newblock The homogeneous coordinate ring of a toric variety.
\newblock {\em J. Algebraic Geom.}, 4(1):17--50, 1995.

\bibitem{Fulton}
W.~Fulton.
\newblock {\em Introduction to toric varieties}, volume 131 of {\em Annals of
  Mathematics Studies}.
\newblock Princeton University Press, Princeton, NJ, 1993.
\newblock The William H. Roever Lectures in Geometry.

\bibitem{Gaitsgory}
D.~Gaitsgory.
\newblock Construction of central elements in the affine {H}ecke algebra via
  nearby cycles.
\newblock {\em Invent. Math.}, 144(2):253--280, 2001.

\bibitem{GKM2}
M.~Goresky, R.~Kottwitz, and R.~MacPherson.
\newblock Equivariant cohomology, {K}oszul duality, and the localization
  theorem.
\newblock {\em Invent. Math.}, 131(1):25--83, 1998.

\bibitem{GKM}
M.~Goresky, R.~Kottwitz, and R.~Macpherson.
\newblock Homology of affine {S}pringer fibers in the unramified case.
\newblock {\em Duke Math. J.}, 121(3):509--561, 2004.

\bibitem{GKM-purete}
M.~Goresky, R.~Kottwitz, and R.~MacPherson.
\newblock Purity of equivalued affine {S}pringer fibers.
\newblock {\em Represent. Theory}, 10:130--146 (electronic), 2006.

\bibitem{EGAIII}
A.~Grothendieck.
\newblock \'{E}l\'ements de g\'eom\'etrie alg\'ebrique. {III}. \'{E}tude
  cohomologique des faisceaux coh\'erents. {I}.
\newblock {\em Inst. Hautes \'Etudes Sci. Publ. Math.}, (11):167, 1961.

\bibitem{KL}
D.~Kazhdan and G.~Lusztig.
\newblock Fixed point varieties on affine flag manifolds.
\newblock {\em Israel J. Math.}, 62(2):129--168, 1988.

\bibitem{Ko}
R.~Kottwitz.
\newblock Rational conjugacy classes in reductive groups.
\newblock {\em Duke Math. J.}, 49(4):785--806, 1982.

\bibitem{Ko1}
R.~Kottwitz.
\newblock Isocrystals with additional structure.
\newblock {\em Compositio Math.}, 56(2):201--220, 1985.

\bibitem{Ko2}
R.~Kottwitz.
\newblock Isocrystals with additional structure. {II}.
\newblock {\em Compositio Math.}, 109(3):255--339, 1997.

\bibitem{K-LTF}
R.~Kottwitz.
\newblock Harmonic analysis on reductive {$p$}-adic groups and {L}ie algebras.
\newblock In {\em Harmonic analysis, the trace formula, and Shimura varieties},
  volume~4 of {\em Clay Math. Proc.}, pages 393--522. Amer. Math. Soc.,
  Providence, RI, 2005.

\bibitem{Kumar}
S.~Kumar.
\newblock Infinite {G}rassmannians and moduli spaces of {$G$}-bundles.
\newblock In {\em Vector bundles on curves---new directions (Cetraro, 1995)},
  volume 1649 of {\em Lecture Notes in Math.}, pages 1--49. Springer, Berlin,
  1997.

\bibitem{MP}
A.~Moy and G.~Prasad.
\newblock Unrefined minimal {$K$}-types for {$p$}-adic groups.
\newblock {\em Invent. Math.}, 116(1-3):393--408, 1994.

\bibitem{W2}
J.-L. Waldspurger.
\newblock L'endoscopie tordue n'est pas si tordue : intégrales orbitales.
\newblock Prépublication Mai 2006.

\bibitem{W1}
J.-L. Waldspurger.
\newblock Endoscopie et changement de caract\'eristique.
\newblock {\em J. Inst. Math. Jussieu}, 5(3):423--525, 2006.

\end{thebibliography}

\begin{flushleft}
CNRS et Universit\'{e} Paris-Sud \\
 UMR 8628 \\
 Math\'{e}matique, B\^{a}timent 425 \\
F-91405 Orsay Cedex \\
France \\
\bigskip

mél :\\
Pierre-Henri.Chaudouard@math.u-psud.fr \\
Gerard.Laumon@math.u-psud.fr\\
\end{flushleft}

\end{document}